\documentclass{ruthesis}

\usepackage{amsmath,amsthm,amssymb,latexsym,hyperref,url,graphicx,cite}

\usepackage{tikz,color,lscape}
\usetikzlibrary{arrows}

\newtheorem{theorem}{Theorem}[chapter]
\newtheorem{conjecture}[theorem]{Conjecture}
\newtheorem{corollary}[theorem]{Corollary}
\newtheorem{definition}[theorem]{Definition}
\newtheorem{example}[theorem]{Example}
\newtheorem{lemma}[theorem]{Lemma}
\newtheorem{proposition}[theorem]{Proposition}

\newcommand{\seqnum}[1]{\href{http://oeis.org/#1}{\underline{#1}}}%
\usepackage{enumitem}%
\usepackage{arydshln}%
\usepackage{float} %
\newcommand{\beq}{\begin{equation*}}
\newcommand{\eeq}{\end{equation*}}
\newcommand{\beqa}{\begin{align*}}
\newcommand{\eeqa}{\end{align*}}
\newcommand{\ben}{\begin{enumerate}}
\newcommand{\een}{\end{enumerate}}
\newcommand{\bed}{\begin{definition}}
\newcommand{\eed}{\end{definition}}
\newcommand{\bet}{\begin{theorem}}
\newcommand{\eet}{\end{theorem}}
\newcommand{\bel}{\begin{lemma}}
\newcommand{\eel}{\end{lemma}}
\newcommand{\bec}{\begin{corollary}}
\newcommand{\eec}{\end{corollary}}
\newcommand{\bep}{\begin{proof}}
\newcommand{\eep}{\end{proof}}

\newcommand{\tab}{\hspace*{2em}}
\newcommand{\tbl}{\textquotedblleft}
\newcommand{\tbr}{\textquotedblright}
\newcommand{\modu}[1]{\ensuremath{(\!\!\!\!\mod #1)}}
\newcommand{\floor}[1]{\ensuremath{\left\lfloor#1\right\rfloor}}
\newcommand{\ceil}[1]{\ensuremath{\left\lceil#1\right\rceil}}

\def\Z{{\mathbb Z}}

\def\N{{\mathbb{N}}}
\def\S{\mathcal{S}}

\makeatletter
\newtheorem*{rep@theorem}{\rep@title}
\newcommand{\newreptheorem}[2]{%
\newenvironment{rep#1}[1]{%
 \def\rep@title{#2 \ref{##1}}%
 \begin{rep@theorem}}%
 {\end{rep@theorem}}}
\makeatother

\newreptheorem{theorem}{Theorem}
\newreptheorem{lemma}{Lemma}

\begin{document}
\phd
\title{Unimodal Polynomials and Lattice Walk Enumeration with Experimental Mathematics}
\author{Bryan Ek}
\program{Mathematics}
\director{Doron Zeilberger}
\approvals{4} %
\submissionyear{2018} %
\submissionmonth{May} %
\abstract{The main theme of this dissertation is retooling methods to work for different situations.

I have taken the method derived by O'Hara and simplified by Zeilberger to prove unimodality of $q$-binomials and tweaked it. This allows us to create many more families of polynomials for which unimodality is not, a priori, given. I analyze how many of the tweaks affect the resulting polynomial.

Ayyer and Zeilberger proved a result about bounded lattice walks. I employ their generating function relation technique to analyze lattice walks with a general step set in bounded, semi-bounded, and unbounded planes. The method in which we do this is formulated to be highly algorithmic so that a computer can automate most, if not all, of the work. I easily recover many well-known results for simpler step sets and discover new results for more complex step sets.}

\beforepreface
\acknowledgements{I would like to thank my advisor Doron Zeilberger for his direction. He provided several fun projects to work on that resulted in this dissertation. I would also like to thank Cole Franks for his edits and suggestions. And thank you to Michael Saks for his discussion and comments on my work that provided many insights for further analysis. Thank you to my committee for serving as such. Thank you to my friends for the welcome diversions. And a big thank you to my whole family for supporting me along this journey!\\

This research was funded by a SMART Scholarship: USD/R\&E (The Under Secretary of Defense-Research and Engineering), National Defense Education Program (NDEP) / BA-1, Basic Research.}

\dedication{\textit{\ \vskip 60pt \noindent Dedicated to my wife Kaitlin.}}

\figurespage
\tablespage
\afterpreface

\begin{chapter}{Introduction}	%

	Humans have been counting for millennia, whether it was to keep track of how many berry bushes remained to gather from, how many people lived in their village, or how many days it had been since the water froze. A single person could simply count the bushes or people or days; the important jump for society was to communicate that number. One of the earliest methods was tally marks. 
	In ancient Babylon 5,000 years ago, the abacus arose as a method of subtracting and adding numbers in shorthand. 
	A classic tale of Gauss' intellect is how he solved the grade school busy-work problem of adding the numbers from 1 to 100, which can be generalized to adding any sequence of numbers with constant difference.
	
	Zeilberger used recurrences and properties of unimodal polynomials inspired by O'Hara to show that the $q$-binomials are unimodal \cite{Zeilberger1989}. 
	He and Ayyer also used generating function relations to prove a result about polymers bounded between plates \cite{oldtimebasketball}.\\

	In addition to the important specific results, the previous examples provide something more consequential: the system of obtaining the result. To the cavemen, the amount of food left was probably their top priority, but what they began paved the way for more complex analysis\footnote{And eventually Complex Analysis.} of our world. The methods allow for a solution to not just the original question or task, but with simple tweaks, can answer much more.
	
	For many early questions, particularly Gauss' busywork, we can have a computer brute-force a result. But while computers are faster, {\it much} faster, than us they still take nonzero time for even the simplest computations. Over the years, we have created machines that can compute faster and faster. But instead of computing faster, we want to compute smarter.
	
	We could go through the process of solving a system of 1,000 linear equations and variables by hand, but that would take forever. Computers are very good at doing simple things quickly. The important facet is the method and set-up process to allow the computer to do its thing.
	
	This is the methodology of experimental mathematics: to use computers efficiently for mathematical investigations. The more data we gather, the easier it is to draw insights.  Computing has come quite a long way for mathematics.\\
	
	The main theme of this dissertation is retooling methods to work for different situations. As a part of this retooling, we have made two major Maple packages for demonstrations of the various devices.\\
	
	We begin by delving into unimodal polynomials. The first chapter is all about how we can tweak a single recurrence that proves unimodality for the $q$-binomial coefficients into a source of many more unimodal polynomials. The recurrence provided by Zeilberger \cite{zeilberger1989kathy} and inspired by O'Hara \cite{Ohara} is given and several adjustments are introduced. How these adjustments affect the final polynomial are explored for most of the chapter. The unimodal polynomial chapter was produced in conjunction with the Maple package {\bf Gnk} available at
	\begin{center}{\bf math.rutgers.edu/$\sim$bte14/Code/Gnk/Gnk.txt}.\end{center}
	The package provides for quick enumeration of families of unimodal polynomials as well as a few functions for insight into the background work.\\
	
	The second chapter is on the subject of lattice walks. Given a certain set of allowable steps in a 2D lattice (none of which move to the left), how many ways are there to walk between points in the plane? To tackle this question, I use generating function relations motivated by Ayyer and Zeilberger \cite{oldtimebasketball}. For walks bounded above and below, we are guaranteed to find an exact expression for the generating function: a rational function. For more general walks, we are only guaranteed to find a polynomial for which the generating function is a root (in terms of formal power series). For some polynomials this means we can give an exact expression but, in general, the {\it minimal polynomial} is the best we, or anyone, can do. The lattice walk chapter was produced in conjunction with the Maple package {\bf ScoringPaths} available at
	\begin{center}{\bf math.rutgers.edu/$\sim$bte14/Code/ScoringPaths/ScoringPaths.txt}.\end{center}
	Within this package are all of the procedures needed to generate the mentioned results automatically. There are several different ways we can enumerate the walks; this chapter discusses them, their pros and cons, and the specific methods of enumeration are included in {\bf ScoringPaths}.
	
\end{chapter}%

\begin{chapter}{Unimodal Polynomials}

	\begin{section}{Introduction}
		This chapter is available as a stand-alone paper on arXiv.org, number 1711.11252. It has been submitted to the {\it Journal of Difference Equations and Applications}.\\
		
		The $q$-binomial coefficients were conjectured to be unimodal as early as the 1850s, but it remained unproven until Sylvester's 1878 proof using invariant theory. In 1982, Proctor gave an \tbl elementary\tbr\ proof using linear algebra. Finally, in 1989, Kathy O'Hara provided a {\it combinatorial} proof of the unimodality of the $q$-binomial coefficients. Very soon thereafter, Doron Zeilberger translated the argument into an elegant recurrence. We introduce several perturbations to the recurrence to create a larger family of unimodal polynomials. We analyze how these perturbations affect the final polynomial and analyze some specific cases.

		\begin{subsection}{Motivation}
			\begin{quote}
				\tbl The study of unimodality and log-concavity arise often in combinatorics, economics of uncertainty and information, and algebra, and have been the subject of considerable research.\tbr \cite{Alvarez2000}
			\end{quote}
			Intuitively, many sequences seem to be unimodal, but how does one prove that fact? We will review some methods of building unimodal sequences as well as a few lemmas that imply unimodality.
			
			Knowing a sequence is unimodal allows for guaranteed discovery of the global extremum using an easy search algorithm. Unimodality is also useful for probability applications. Identifying a probability distribution as unimodal allows certain approximations for how far a value will be from its mode (Gauss' inequality \cite{GaussInequality}) or mean (Vysochanskij-Petunin inequality \cite{VPinequality}).\\
			
			Suppose we start with a set of nice combinatorial objects that satisfy property $\mathcal{P}$. How can we make the set larger and still satisfy property $\mathcal{P}$? Or a slightly different property $\mathcal{P}'$? The goal of this project is to use reverse engineering to obtain highly non-trivial, surprising theorems about unimodality.

			\begin{example}
				All of the following functions are not only polynomial for $n\ge0,0,0,4,4$, respectively, but they are also unimodal.
				\begin{align*}
					P_1(n)	&=	{\frac {1-{q}^{n+1}}{1-q}},\\
					P_2(n)	&=	{\frac {1- q^{2n+1}}{1-q}},\\
					P_3(n)	&=	{\frac {2q^{3n+2}-2q^{3n+1}+q^{3n}- q^{2n+1}-q^{n+1}+{q}^{2}-2\,q+2}{ \left( 1-q \right) ^{2}}},\\
					P_4(n)	&=	\bigg(5q^{4n+2}-5q^{4n+1}+3{q}^{4n-2}-3q^{4n-3}+4q^{4n-4}-4{q}^{2n+3}-4{q}^{2n-1}\\
							&\tab	+4{q}^{6}-3{q}^{5}+3{q}^{4}-5q+5\bigg)\bigg/\left(1-q\right) ^{2},\\
					P_5(n)	&=	\bigg(2q^{5n+3}-4q^{5n+2}+7{q}^{5n+1}-5{q}^{5n}+5{q}^{5n-3}-5{q}^{5n-4}+4{q}^{5n-5}-4{q}^{5n-6}+6{q}^{5n-9}\\
							&\tab	-6q^{5n-10}+3q^{5n-11}	-5{q}^{4n+2}+5{q}^{4n+1}-4{q}^{4n-2}+4{q}^{4n-3}-3{q}^{4n-4}-5{q}^{3n+4}\\
							&\tab	+5{q}^{3n+3}-3{q}^{3n+2}-6{q}^{3n-4}+6{q}^{3n-5}-3{q}^{3n-6}+3{q}^{2n+9}-6{q}^{2n+8}+6{q}^{2n+7}\\
							&\tab+3{q}^{2n+1}-5{q}^{2n}+5{q}^{2n-1}+3{q}^{n+7}-4{q}^{n+6}+4{q}^{n+5}-5{q}^{n+2}+5{q}^{n+1}\\
							&\tab	-3{q}^{14}+6{q}^{13}-6{q}^{12}+4{q}^{9}-4{q}^{8}+5{q}^{7}-5{q}^{6}+5{q}^{3}-7q^{2}+4q-2\bigg)\bigg/\left( -1+q \right) ^{3}.
				\end{align*}
				\label{exa:ImpressivePolynomials1}
			\end{example}
			
			\begin{example}
				The following functions are all unimodal for $n\ge0$.
				\begin{align*}
					Q_1(n)	&=	{\frac {4\left(1-{q}^{n+1}\right)}{1-q}},\\
					Q_2(n)	&=	5{\frac {8{q}^{2n+4}-2{q}^{n+3}-2{q}^{n+2}+8q-3[(q^5+1)(q^n-(-q)^n)+(q^4+q)(q^n+(-q)^n)]}{2q\left(1-q\right)\left(1-q^2\right)}},\\
					Q_3(n)	&=	\bigg(16\big[3-3q+16{q}^{2}-16{q}^{n+1}-16{q}^{n+3}-16{q}^{n+5}+16{q}^{2n+3}\\
							&\tab\tab	+16{q}^{2n+5}+16{q}^{2n+7}-16{q}^{3n+6}+3{q}^{3n+7}-3{q}^{3n+8}\big]\\
							&\tab+(1-(-1)^n)(-q)^{(3n-9)/2}\big[8{q}^{17}-8{q}^{16}+64{q}^{15}-9{q}^{14}-55{q}^{13}-60{q}^{11}+57{q}^{10}\\
							&\tab\hspace{3cm}	-125{q}^{9}+125{q}^{8}-57{q}^{7}+60{q}^{6}+55{q}^{4}+9{q}^{3}-64{q}^{2}+8{q}-8\big]\\
							&\tab+(1-(-1)^n)q^{(3n-9)/2}\big[8{q}^{17}-8{q}^{16}+64{q}^{15}+9{q}^{14}-73{q}^{13}-60{q}^{11}-65{q}^{10}\\
							&\tab\hspace{3.7cm}	-3{q}^{9}+3{q}^{8}+65{q}^{7}+60{q}^{6}+73{q}^{4}-9{q}^{3}-64{q}^{2}+8q-8\big]\\
							&\tab+(1+(-1)^n)(-q)^{(3n-6)/2}\big[12{q}^{14}-12{q}^{13}+64{q}^{12}-75{q}^{11}+74{q}^{10}-127{q}^{9}\\
							&\tab\hspace{4.4cm}	+127{q}^{5}-74{q}^{4}+75{q}^{3}-64{q}^{2}+12q-12\big]\\
							&\tab+(1+(-1)^n)q^{(3n-6)/2}\big[12{q}^{14}-12{q}^{13}+64{q}^{12}-53{q}^{11}-74{q}^{10}-{q}^{9}\\
							&\tab\hspace{3.2cm}+{q}^{5}+74{q}^{4}+53{q}^{3}-64{q}^{2}+12q-12\big]\bigg)
							\bigg/4(1-q)^2(1-q^6).
				\end{align*}
				\label{exa:ImpressivePolynomials2}
				$Q_3$ is truly amazing. It appears quite unwieldy and at first glance one might doubt that it has real coefficients. But for any nonnegative integer $n$, $Q_3(n)$ is guaranteed to be a unimodal polynomial in $\Z[q]$. Try simplifying $Q_3$ assuming $n$ is even/odd. One can obtain further simplifications assuming $n=0,1,2,3\modu{4}$.
			\end{example}

		\end{subsection}%

		\begin{subsection}{Symmetric and Unimodal}
			We recall several definitions and propositions from Zeilberger \cite{Zeilberger1989} for the sake of completeness.
			\bed[Unimodal]
				A sequence $\mathcal{A}=\{a_0,\ldots,a_n\}$ is {\it unimodal} if it is weakly increasing up to a point and then weakly decreasing, i.e., there exists an index $i$ such that $a_0\le a_1\le\cdots\le a_i\ge\cdots\ge a_n$.
			\eed
			\bed[Symmetric]
				A sequence $\mathcal{A}=\{a_0,\ldots,a_n\}$ is {\it symmetric} if $a_i=a_{n-i}$ for every $0\le i\le n$.
			\eed
			A polynomial is said to have either of the above properties if its sequence of coefficients has the respective property.
			
			\bed[Darga]
				The {\it darga} of a polynomial $p(q)=a_iq^i+\cdots+a_jq^j$, with $a_i\ne0\ne a_j$, is defined to be $i+j$, i.e., the sum of its lowest and highest powers.
			\eed
			darga is equal to $2\cdot C(p)$ used by Brent \cite{BRENTI1990}.
			\begin{example}
				$darga(q^2+3q^3)=5$ and $darga(q^2)=4$.
			\end{example}

			\begin{proposition}
				The sum of two symmetric and unimodal polynomials of darga $m$ is also symmetric and unimodal of darga $m$.
				\label{pro:SumSymUni}
			\end{proposition}
			\begin{proposition}
				The product of two symmetric and unimodal nonnegative\footnote{Nonnegative is necessary: $(-1+x-x^2)^2=1-2x+3x^2-2x^3+x^4$.} polynomials of darga $m$ and $m'$ is a symmetric and unimodal polynomial of darga $m+m'$.
				\bep
					A polynomial is symmetric and unimodal of darga $m$ if and only if it can be expressed as a sum of \tbl atomic\tbr entities of the form $c(q^{m-r}+q^{m-r-1}+\cdots+q^{r})$, for some {\bf positive} constant $c$ and integer $0\le r\le\frac{m}{2}$. By Proposition \ref{pro:SumSymUni}, it is enough to prove that the product of two such atoms of dargas $m$ and $m'$ is symmetric and unimodal of darga $m+m'$:
					\beq
						(q^{m-r}+\cdots+q^{r})(q^{m'-r'}+\cdots+q^{r'})=q^{m+m'-r-r'}+2q^{m+m'-r-r'-1}+\cdots+2q^{r+r'+1}+q^{r+r'}.
					\eeq
				\eep
				\label{pro:ProdSymUni}
			\end{proposition}
			\begin{proposition}
				If $p$ is symmetric and unimodal of darga $m$, then $q^\alpha p$ is symmetric and unimodal of darga $m+2\alpha$.
				\label{pro:PowerSymUni}
			\end{proposition}
			One example is the binomial polynomial: $(1+x)^m$. It is symmetric and unimodal of darga $m$.
			
			\bed[$\gamma$-nonnegative \cite{Branden2015}]
				If $h(x)$ is a symmetric function (in its coefficients), then we can write $h(x)=\sum_{k=0}^{\floor{m/2}}\gamma_kx^k(1+x)^{m-2k}$. We call $\{\gamma_k\}_{k=0}^{\floor{m/2}}$ the {$\gamma$-vector} of $h$. If the $\gamma$-vector is nonnegative, $h$ is said to be {\it $\gamma$-nonnegative}.
			\eed
			One can use Propositions \ref{pro:SumSymUni}, \ref{pro:ProdSymUni}, and \ref{pro:PowerSymUni} to prove that $\gamma$-nonnegative implies symmetric and unimodal of darga $m$.
			
			\begin{example}
				To prove the polynomials in Example \ref{exa:ImpressivePolynomials1} are actually unimodal, one can show
				\begin{align*}
					P_1(n)	&=	\sum_{i=0}^{n}{q}^{i},\\
					P_2(n)	&=	P_1(2n),\\
					P_3(n)	&=	2P_1\left(3n\right) +{q}^{2}P_1\left(2n-2\right)P_1\left(n-2\right),\\
					P_4(n)	&=	5P_1\left(4n\right)+3{q}^{4}P_2\left(2n-4\right)+4{q}^{6}P_1\left(2n-4\right)P_2\left(n-4\right),\\
					P_5(n)	&=	2P_1\left(5n\right)+5{q}^{2}P_1\left(4n-2\right)P_1\left(n-2\right)+4{q}^{8}P_2\left(2n-6\right)P_1\left(n-4\right)\\
							&\tab	+5{q}^{6}P_1\left(3n-4\right)P_2\left(n-4\right)+3{q}^{12}P_1\left(2n-6\right)P_3\left(n-6\right),
				\end{align*}
				and then use Propositions \ref{pro:SumSymUni}, \ref{pro:ProdSymUni}, and \ref{pro:PowerSymUni} along with the property that $darga(P_k(n))=nk$. One must initially assume $n\ge 0,0,2,4,6$, respectively, but the $k=3,5$ bounds can be lowered by checking smaller values of $n$ manually.
			\end{example}
			
			\begin{example}
				To prove the polynomials in Example \ref{exa:ImpressivePolynomials2} are in fact unimodal, one can show
				\begin{align*}
					Q_1(n)	&=	4\sum_{i=0}^{n}{q}^{i},\\
					Q_2(n)	&=	{q}^{2}Q_2\left(n-2\right)+5Q_1\left(2n\right),\\
					Q_3(n)	&=	{q}^{6}Q_3\left(n-4\right)+3Q_1\left(3n\right)+4{q}^{2}Q_1\left(2n-2\right)Q_1\left(n-2\right),
				\end{align*}
				and then use induction on $n$ for $Q_2,Q_3$, with base cases $n=\{0,1\},\{0,1,2,3\}$, respectively, and Propositions \ref{pro:SumSymUni}, \ref{pro:ProdSymUni}, and \ref{pro:PowerSymUni} along with the property that $darga(Q_k(n))=nk$.
			\end{example}
			The wonderful aspect about these proofs is that they can be easily verified by computer!%
			
			\smallbreak
			
			Being able to write a function in \tbl decomposed\tbr\ form allows us to quickly verify the unimodal nature of each part and hence the sum. However, the combined form generally does not lead to an obvious decomposition. This is the motivation behind using reverse engineering: to guarantee that the decomposed form exists.

			\begin{subsubsection}{Real-rootedness and Log-concavity}
				The following properties are greatly related to unimodality but turn out to not be applicable in our situation. They are included for completeness of discussion.
				\bed[Real-rootedness]
					The generating polynomial, $p_\mathcal{A}(x) := a_0 + a_1q +\cdots+ a_nq^n$, is called {\it real-rooted} if all its zeros are real. By convention, constant polynomials are considered to be real-rooted.
				\eed
				\bed[Log-concavity]
					A sequence $\mathcal{A}=\{a_0,\ldots,a_n\}$ is {\it log-concave (convex)} if $a_i^2\ge a_{i-1}a_{i+1}$ ($a_i^2\le a_{i-1}a_{i+1}$) for all $1\le i<n$.
				\eed
				There is also a notion of {\it $k$-fold log-concave} and {\it infinitely log-concave} with open problems that may be of interest to the reader \cite{Branden2015}.
				\begin{proposition}
					The Hadamard (term-wise) product of log-concave (convex) sequences is also log-concave (convex).
				\end{proposition}
				\bel[Br\"{a}nd\'en \cite{Branden2015}]
					Let $\mathcal{A}=\{a_k\}_{k=0}^n$ be a finite sequence of nonnegative numbers.
					\begin{itemize}
						\item	If $p_\mathcal{A}(x)$ is real-rooted, then the sequence $\mathcal{A}':=\{a_k/\binom{n}{k}\}_{k=0}^n$ is log-concave.
						\item	If $\mathcal{A}'$ is log-concave, then so is $\mathcal{A}$.
						\item	If $\mathcal{A}$ is log-concave and positive,\footnote{Nonnegative is not sufficient: $\{1,0,0,1\}$ is log-concave (and log-convex) but not unimodal.} then $\mathcal{A}$ is unimodal.
					\end{itemize}
					\label{lem:Branden}
				\eel
				For a self-contained proof using less general (and possibly easier to understand) results, see Lecture 1 from Vatter's Algebraic Combinatorics class \cite{Vatter2009}, which uses the book {\it A Walk Through Combinatorics} \cite{Bona2006}. The converse statements of Lemma \ref{lem:Branden} are false. There are log-concave polynomials that are not real-rooted and there are unimodal polynomials which are not log-concave.
				
				Br\"and\'en also references a proof by Stanley \cite{Stanley1989} that:
				\begin{lemma}
					If $A(x),B(x)$ are log-concave, then $A(x)B(x)$ is log-concave. And if $A(x)$ is log-concave and $B(x)$ is unimodal, then $A(x)B(x)$ is unimodal.
					\label{lem:Stanley}
				\end{lemma}
				It is not sufficient for $A,B$ to be unimodal: $(3+q+q^2)^2=9+6q+7q^2+2q^3+q^4$.\footnote{Stanley \cite{Stanley1989} incorrectly writes the coefficient of $q^3$ as $1$.}
			\end{subsubsection}%
			
		\end{subsection}%

		\begin{subsection}{Partitions}
			\bed[Partition]
				A {\it partition} of $k$ is a non-increasing sequence of positive integers $\lambda=[a_1,a_2,\ldots,a_s]$ s.t. $\sum_{i=1}^s a_i=k$.\\
				We will use {\it frequency representation} $\lambda=[a_1^{b_1},a_2^{b_2},\ldots,a_r^{b_r}]$ s.t. $a_i>a_{i+1}$ and $b_i>0$ to abbreviate repeated terms \cite{Andrews1998}.\footnote{If $b_i=1$, it is omitted.} Note now that $\sum_{i=1}^r a_i b_i=k$. There is possible ambiguity as to whether $x^y$ is $x$ repeated $y$ times, or $x^y$ counted once. We always reserve exponentials in partitions to indicate repetition.\\
				The {\it size} of a partition, denoted $|\lambda|$, will indicate the {\bf number} of parts.\footnote{Many papers use $|\lambda|=k$ to denote what $\lambda$ partitions.} In standard notation $|\lambda|=s$; in frequency notation $|\lambda|=\sum_{i=1}^r b_i$.\\
				$d_i$ is the number of parts of size $i$ in the partition. In frequency notation, $d_{a_i}=b_i$.\\
				$\lambda\vdash k$ indicates $\lambda$ is a partition of $k$.
			\eed
			We will use $\lambda$ to denote a partition. Unless otherwise stated, $\lambda$ is a partition of $k$.
		\end{subsection}%

		\begin{subsection}{Chapter Organization}
			This chapter is organized in the following sections:
			\ben[label=\bfseries Section \arabic*:,leftmargin=21mm]
				\setcounter{enumi}{1}
				\item Maple Program: Briefly describes the accompanying Maple package.
				\item	$q$-binomial Polynomials: Describes the original interesting polynomials.
				\item	Original Recurrence: Introduces the recurrence that generates the $q$-binomial polynomials and simultaneously proves their unimodality. We briefly discuss properties of the recurrence itself.
				\item	Altered Recurrence: Modifications are injected into the recurrence that still maintain unimodality for the resulting polynomials. We examine the effects various changes have. This section contains most of the reference to new unimodal polynomials.
				\item	OEIS: Uses some of the created unimodal polynomials to enumerate sequences for the OEIS.\footnote{Online Encyclopedia of Integer Sequences \cite{OEIS}.}
				\item	Conclusion and Future Work: Provides avenues for future research.
			\een
		\end{subsection}%
		
	\end{section}%

	\begin{section}{Maple Program}
		The backbone of this paper is based on experimental work with the Maple package {\bf Gnk} available at
		\begin{center}{\bf math.rutgers.edu/$\sim$bte14/Code/Gnk/Gnk.txt}.\end{center}
		I will mention how functions are implemented throughout this paper using {\sc Implemented as {\bf...}}. For help with any function, type {\bf Help(function)}.
		
		The most important function is {\bf KOHgeneral}. It modifies the original recurrence in Eqn.\ \eqref{eqn:KOH} in several different ways. Using this function, one can change the recurrence call, multiply the summand by any manner of constant, restrict the types of partitions, or even hand-pick which partitions should be weighted most highly. The key part of these changes is that the solution to the recurrence remains symmetric and unimodal of darga $nk$.
		
		{\bf math.rutgers.edu/$\sim$bte14/Code/Gnk/Polynomials/} contains examples of many non-trivial one-variable rational functions that are guaranteed to be unimodal polynomials (for $n\ge n_0$). They were created using the {\bf KOHrecurse} function and restricting the partitions summed over in \eqref{eqn:KOH} to have smallest part $\ge1,2,3,4$ and distance $\ge1,2,3,4$. In decomposed form, they are clearly unimodal; but in combined form one is hard-pressed to state unimodality with certainty.
		
		The \tbl impressive\tbr\ polynomials in Examples \ref{exa:ImpressivePolynomials1} and \ref{exa:ImpressivePolynomials2} were generated by {\bf RandomTheoremAndProof}. It uses the recurrence in Eqn.\ \eqref{eqn:KOH} with random multiplicative constants in each term. The parameters were $(n,q,5,5,0,false)$ and $(n,q,3,5,0,true)$, respectively. Initially, $Q_3(n)$ was an even larger behemoth. This will typically happen when Maple solves a recurrence equation so going beyond $K=3$ and $complicated=true$ is not recommended unless one wants to spend substantial time parsing the polynomial into something much more readable. Or better yet, one could come up with a way for Maple to do this parsing automatically.
		
		There are functions to test whether a polynomial is symmetric, unimodal, or both ({\bf isSymmetric(P,q)}, {\bf isUnimodal(P,q)}, {\bf isSymUni(P,q)}, respectively). As an extra take-away, {\bf generalPartitions} outputs all partitions restricted to a minimum/maximum integer/size as well as difference\footnote{Consecutive parts in the partition differ by $\ge d$.} and congruence requirements. One can also specify a finite set of integers which are allowed in {\bf genPartitions}.
	\end{section}%

	\begin{section}{$q$-binomial Polynomials}
		Of particular interest among unimodal polynomials are the $q$-binomial polynomials ${{n}\brack{k}}_q$ {\sc Implemented as {\bf qbin(n,k,q)}}. The $q$-binomials are analogs of the binomial coefficients. In the limit $q\to1^-$, we obtain the usual binomial coefficients. We parametrize $q$-binomial coefficients in order to eliminate the redundancy: ${{n}\brack{k}}_q={{n}\brack{n-k}}_q$:
		\begin{equation}
			G(n,k)	=	G(k,n)	=	{{n+k}\brack{k}}_q	=	\frac{[n+k]_q!}{[n]_q![k]_q!}	=	\frac{(1-q^{n+1})\cdots(1-q^{n+k})}{(1-q)\cdots(1-q^k)},
			\label{eqn:Qbinomial}
		\end{equation}
		{\sc Implemented as {\bf Gnk(n,k,q)}}. $[n]_q=\frac{1-q^n}{1-q}$ is the q-bracket. $[n]_q!=\prod_{i=1}^n[i]_q$ denotes the q-factorial\footnote{The q-factorial is the generating polynomial for the number of inversions over the symmetric group $S_n$ \cite{Branden2015}.} {\sc Implemented as {\bf qfac(n,q)}}. It is non-trivial that $G(n,k)$ is even a polynomial for integer $n,k$. Polynomial is proven since $G(n,k)$ satisfies the simple recurrence relation
		\begin{equation}
			G(n+1,k+1)	=	q^{k+1}G(n,k+1)+G(n+1,k).
			\label{eqn:Recurrence}
		\end{equation}
		This recurrence also shows that $G(n,k)$ is the generating function for the number of partitions in an $n\times k$ box. A partition either has a part of largest possible size ($q^{k+1}G(n,k+1)$) or it does not ($G(n+1,k)$). To prove this, one can simply check that $G(n,k)$ does satisfy this recurrence as well as the initial conditions $G(n,0)=1=G(0,k)$.\footnote{The coefficients are symmetric since each partition of $\ell$ inside an $n\times k$ box corresponds to a partition of $nk-\ell$.} This argument also confirms $G(n,k)=G(k,n)$.\footnote{$G(n-k,k)$ also counts the number of subspaces of dimension $k$ in a vector space of dimension $n$ over a finite field with $q$ elements.} %

		The question of unimodality for the $q$-binomial coefficients was first stated in the 1850s by Cayley and then proven by Sylvester in 1878 \cite{Branden2015}. The first \tbl elementary\tbr\ proof was given by Proctor in 1982. He essentially described the coefficients as partitions in a box (grid-shading problem) and then used linear algebra to finish the proof \cite{Proctor1982}.
		
		Lemma \ref{lem:Branden} does not apply here; q-brackets and q-factorials are in the class of log-concave (use Lemma \ref{lem:Stanley}) but not real-rooted (for $n>1$) polynomials while these $G(n,k)$ polynomials are in the class of unimodal, but not necessarily log-concave polynomials. For example, $G(2,2)=1+q+2q^2+q^3+q^4$ is unimodal but not log-concave. In fact, most $q$-binomial polynomials are not log-concave since the first (and last) 2 coefficients are 1. Explicitly: $G(n,k)$ is log-concave $\Leftrightarrow$ $k=1$.
		
	\end{section}%

	\begin{section}{Original Recurrence}
		The symmetric and structured nature of $G(n,k)$ and its coefficients is visually muddled by the following recurrence that Zeilberger \cite{Zeilberger1989} created to translate O'Hara's \cite{Ohara} combinatorial argument\footnote{Her argument is rewritten/explained by Zeilberger \cite{zeilberger1989kathy}.} to mostly algebra.
		\begin{align}
			G(n,k)	&=	\sum_{(d_1,\ldots,d_k);\sum_{i=1}^k id_i=k} q^{k\left(\sum_{i=1}^kd_i\right)-k-\sum_{1\le j<i\le k}(i-j)d_id_j}\nonumber\\
					&\tab\tab	\prod_{i=0}^{k-1}G\left((k-i)n-2i+2\sum_{j=0}^{i-1}(i-j)d_{k-j},d_{k-i}\right)
			\label{eqn:KOH}
		\end{align}
		is {\sc Implemented as {\bf KOH(n,k,q)}}. The product is part of the summand and the outside sum is over all partitions of $k$. $d_i$ is the number of parts of size $i$ in the partition. The initial conditions\footnote{$G(n,1)$ is explicitly given because Eqn.\ \eqref{eqn:KOH} only yields the tautology $G(n,1)=G(n,1)$.} are
		\begin{equation}
			G(n<0,k)	=	0,	\hspace{1.35em}	G(n,k<0)	=	0,	\hspace{1.35em}	G(0,k)	=	1,	\hspace{1.35em}	G(n,0)	=	1,	\hspace{1.35em}	G(n,1)	=	\frac{1-q^{n+1}}{1-q}.
			\label{eqn:InitialConditions}
		\end{equation}
		By using Propositions \ref{pro:SumSymUni}, \ref{pro:ProdSymUni}, and \ref{pro:PowerSymUni}, and induction on the symmetric unimodality of $G(a,b)$ for $a<n$ or $b<k$, we can see (after a straightforward calculation) that the right hand side is symmetric and unimodal of darga $nk$. For each partition, the darga will be:
		\beq
			2\left[k\left(\sum_{i=1}^kd_i\right)-k-\sum_{1\le j<i\le k}\hspace{-0.8em}(i-j)d_id_j\right]+\sum_{i=0}^{k-1}\left((k-i)n-2i+2\sum_{j=0}^{i-1}(i-j)d_{k-j}\right)d_{k-i}	=	nk.
		\eeq
		Recall that $\sum_{i=1}^k id_i=k$.
		
		\bigskip
		
		One is led to ask how long this more complicated recurrence will take to compute? What is the largest depth of recursive calls that will be made in Eqn.\ \eqref{eqn:KOH}? This is answered with a brute-force method in {\bf KOHdepth}. {\bf KOHcalls} returns the total number of recursive calls made. But is there an explicit answer?
		
		We will use $n',k'$ to denote the new $n$, $k$ in the recursive call for a given partition. We sometimes also treat them as functions of the index: $n'(i)=(k-i)n-2i+2\sum_{j=0}^{i-1}(i-j)d_{k-j}$. Once a partition only has distinct parts, the recurrence will cease after a final recursive call, since each $k'<=1$. Therefore to maximize depth, we need to maximize using partitions with repeated parts and maximize the size of $n'$, $k'$.
		
		One method is to begin by using $\left[\floor{\frac{k}{2}}^2\right]$ with a $0,1$ as needed and the recursive call with $i=\ceil{\frac{k}{2}}$. Then (since $k'=d_{k-i}=2$) repeatedly use the partition $[1^2]$ and $i=1$. This method calls
		\begin{align*}
			G&\left(\left(k-\ceil{\frac{k}{2}}\right)n-2\ceil{\frac{k}{2}}+2\sum_{j=0}^{\ceil{k/2}-1}\left(\ceil{\frac{k}{2}}-j\right)d_{k-j},d_{k-\ceil{k/2}}\right)\\
				&\tab=	G\left(\floor{\frac{k}{2}}n-2\ceil{\frac{k}{2}}+2\sum_{j=0}^{\ceil{k/2}-1}0,d_{\floor{k/2}}\right)\\
				&\tab=	G\left(\floor{\frac{k}{2}}n-2\ceil{\frac{k}{2}},2\right).
		\end{align*}
		Then using the partition $[1^2]$, and $i=1$, calls $G(n'-2,2)$. Thus, the total depth of calls is
		\beq
			\ceil{\frac{\floor{k/2}n-2\ceil{k/2}}{2}}+1	=	\ceil{\frac{\floor{k/2}n}{2}}-\ceil{\frac{k}{2}}+1.
		\eeq
		In fact
		\bel
			If $k=1,2$ or if $k\ge4$ is even (odd) and $n\ge2$ ($n\ge4$), then the maximum depth of recursive calls is
			\beq
				\ceil{\frac{\floor{k/2}n}{2}}-\ceil{\frac{k}{2}}+1.
			\eeq
			For $k=3$, the maximum depth of recursive calls is $\ceil{\frac{n}{4}}$.
			\bep
				Base case: $k=1$. There are no recursive calls:
				\beq
					\ceil{\frac{\floor{1/2}n}{2}}-\ceil{\frac{1}{2}}+1	=	\ceil{\frac{0\cdot n}{2}}-1+1	=	0.
				\eeq
				For $k=2$, the only option that yields recursive calls is $[1^2]$ and $i=1$: $G(n-2,2)$. Thus, we will have
				\beq
					\ceil{\frac{n}{2}}	=	\ceil{\frac{\floor{2/2}n}{2}}-\ceil{\frac{2}{2}}+1
				\eeq
				recursive calls, which matches with expected. For $k=3$, the only recursive call that does not immediately terminate is $\lambda=[1^3]$ and $i=2$. The call is
				\beq
					G\left((3-2)n-2\cdot2+2\sum_{j=0}^{2-1}(2-j)d_{3-j},d_{3-2}\right)	=	G(n-4,3).
				\eeq
				Thus, we will have $\ceil{\frac{n}{4}}$ recursive calls. The difference from other $k$ arises because $\floor{\frac{3}{2}}=1$.\\
				
				First consider $k=$ even and $n=2$. By Proposition \ref{pro:MaxDepth1}, the depth is
				\beq
					1	=	\frac{kn}{4}-\frac{k}{2}+1	=	\ceil{\frac{\floor{k/2}n}{2}}-\ceil{\frac{k}{2}}+1.
				\eeq
				
				Now let $k\ge4$ and $n\ge3$ ($n\ge4$ if $k=$ odd). Assume true for $(n',k')$ s.t. $k'<k$ or $k'=k$ and $n'<n$.\\
				What if we repeat part $m$, $\ell\ge2$ times? The recursive call uses $i=k-m$.
				\begin{align*}
					G&\left(\left(k-k+m\right)n-2\left(k-m\right)+2\sum_{j=0}^{k-m-1}\left(k-m-j\right)d_{k-j},d_{m}\right)\\
					&\tab=	G\left(mn-2k+2m+2\sum_{j=m+1}^{k}\left(j-m\right)d_{j},\ell\right).
				\end{align*}
				By the induction hypothesis, the depth of this call is (note $ml\le k$)
				\begin{align*}
					1-\ceil{\frac{\ell}{2}}+&\ceil{\frac{\floor{\ell/2}\left(mn-2k+2m+2\sum_{j=m+1}^{k}\left(j-m\right)d_{j}\right)}{2}}\\
							&\le	1-\ceil{\frac{l}{2}}+\ceil{\frac{\floor{l/2}\left(mn-2k+2m+2\sum_{j=m+1}^{k}jd_{j}\right)}{2}}\\
							&\le	1-\ceil{\frac{\ell}{2}}+\ceil{\frac{\floor{\ell/2}\left(mn-2k+2m+2(k-\ell m)\right)}{2}}\\
							&=	1-\ceil{\frac{\ell}{2}}+\ceil{\frac{\floor{\ell/2}m\left(n+2-2\ell\right)}{2}}\\
							&\le	1-\ceil{\frac{\ell}{2}}+\ceil{\frac{m\ell\left(n+2-2\ell\right)}{4}}\\
							&\le	1-\ceil{\frac{\ell}{2}}+\ceil{\frac{k\left(n+2-2\ell\right)}{4}}.
				\end{align*}
				The upper bound is maximized when $\ell$ is minimized: $\ell=2$, upper bound $=\ceil{\frac{kn}{4}-\frac{k}{2}}$. We now need to confirm that $\lambda=\left[\floor{\frac{k}{2}}^2\right]$ (with $1$ if needed), $i=\ceil{\frac{k}{2}}$ achieves this upper bound. The depth is
				\begin{align*}
					1-\ceil{\frac{2}{2}}+\ceil{\frac{\floor{2/2}\left(\floor{\frac{k}{2}}n-2k+2\floor{\frac{k}{2}}\right)}{2}}	=	\ceil{\frac{\left(\floor{\frac{k}{2}}n-2\ceil{\frac{k}{2}}\right)}{2}}	=	\ceil{\frac{\floor{k/2}n}{2}-\ceil{\frac{k}{2}}},
				\end{align*}
				which for $k=$ even directly matches the upper bound. For $k=$ odd, we adjust the upper bound by recognizing that if $\ell=$ even, then $m\ell\le k-1$ and splitting into 2 cases:
				\ben
					\item	$\sum_{j=m+1}^{k}jd_{j}=0$. Then the depth of the call is
						\begin{align*}
							1-\ceil{\frac{\ell}{2}}+&\ceil{\frac{\floor{\ell/2}\left(mn-2k+2m+2\sum_{j=m+1}^{k}\left(j-m\right)d_{j}\right)}{2}}\\
								&=	1-\frac{\ell}{2}+\ceil{\frac{m\ell(n+2)-2k\ell}{4}}\\
								&\le	1-\frac{\ell}{2}+\ceil{\frac{(k-1)(n+2)-2k\ell}{4}}:
						\end{align*}
						maximized for $\ell=2$ and gives an upper bound $\ceil{\frac{(k-1)(n+2)-4k}{4}}=\ceil{\frac{(k-1)n}{4}-\frac{k+1}{2}}$.
					\item	$\sum_{j=m+1}^{k}jd_{j}>0$. Then $m\ell\le k-(m+1)\le k-2$ and the depth of the call is
						\begin{align*}
							1-\ceil{\frac{\ell}{2}}+&\ceil{\frac{\floor{\ell/2}\left(mn-2k+2m+2\sum_{j=m+1}^{k}\left(j-m\right)d_{j}\right)}{2}}\\
								&\le	1-\frac{\ell}{2}+\ceil{\frac{\ell\left(mn-2k+2m+2(k-\ell m)-2m\right)}{4}}\\
								&=	1-\frac{\ell}{2}+\ceil{\frac{m\ell\left(n-2\ell\right)}{4}}\\
								&\le	1-\frac{\ell}{2}+\ceil{\frac{(k-2)\left(n-2\ell\right)}{4}},
						\end{align*}
						which again is maximized for $\ell=2$ and gives an upper bound $\ceil{\frac{(k-2)(n-4)}{4}}=\ceil{\frac{(k-1)n}{4}-\frac{k+1}{2}-\frac{2k+n-10}{4}}\le\ceil{\frac{(k-1)n}{4}-\frac{k+1}{2}}$ since $k\ge4$ and $n\ge2$.\footnote{Actually, we have $k=$ odd so $k\ge5$ and $n\ge4$; the bound improves further.}
				\een
				
				\medskip
				
				The upper bound from either case matches the achievable bound for $k=$ odd. If $\ell=$ odd, then the earlier $\floor{\frac{\ell}{2}}\le\frac{\ell}{2}$ can be replaced by $=\frac{\ell-1}{2}$. Thus, the upper bound is
				\beq
					\le	1-\ceil{\frac{\ell}{2}}+\ceil{\frac{m(\ell-1)\left(n+2-2\ell\right)}{4}}	\le	1-\ceil{\frac{\ell}{2}}+\ceil{\frac{(k-1)\left(n+2-2\ell\right)}{4}},
				\eeq
				which is now maximized by $\ell=3$ (since we assume $\ell=$ odd)\footnote{Actually, $\ell=5$ since we ignore $\ell=3$.} giving the upper bound: $\ceil{\frac{(k-1)\left(n-4\right)}{4}}-1$, which is worse\footnote{We only need $k\ge1$.} than the bound when $\ell=2$ so we can ignore $\ell=$ odd when it comes to maximizing depth. Therefore our chosen $\lambda$, and $i$, achieve the greatest depth: $1+\ceil{\frac{\floor{k/2}n}{2}}-\ceil{\frac{k}{2}}$.
				
				To confirm that this recursive depth actually happens, we must also confirm that\linebreak$G(n',k')\ne0$ for any other $i$. I.e., we must show that $n'\ge0$. For $i<\ceil{\frac{k}{2}}$
				\beq
					n'	\hspace{-0.15em}=\hspace{-0.15em}	(k-i)n-2i+2\sum_{j=0}^{i-1}(i-j)d_{k-j}	\hspace{-0.15em}=	\hspace{-0.15em}	(k-i)n-2i	\hspace{-0.1em}>\hspace{-0.3em}	\left(k-\ceil{\frac{k}{2}}\right)n-2\ceil{\frac{k}{2}}	\hspace{-0.15em}=\hspace{-0.15em}	\floor{\frac{k}{2}}n-2\ceil{\frac{k}{2}}	\ge	0
				\eeq
				since $n\ge3$ when $k=$ even and since $n\ge4$ when $k=$ odd. For $i>\ceil{\frac{k}{2}}$,
				\begin{align*}
					n'	&=	(k-i)n-2i+2\sum_{j=0}^{i-1}(i-j)d_{k-j}	=	(k-i)n-2i+2\left(i-\ceil{\frac{k}{2}}\right)2\\
						&=	(k-i)n+2i-4\ceil{\frac{k}{2}}	=	kn-i(n-2)-4\ceil{\frac{k}{2}}	\ge	kn-(k-1)(n-2)-4\ceil{\frac{k}{2}}\\
						&=	kn-kn+n+2k-2-4\ceil{\frac{k}{2}}	=	n-2+2k-4\ceil{\frac{k}{2}},
				\end{align*}
				which for $k=$ even, is $\ge0$ since $n\ge3$. And for $k=$ odd $n'\ge0$, since $n\ge4$.\\
				By induction, the claim holds for all $n,k$.\\
				
				We have glossed over a few details. We assumed that $k'=\ell\ne3$. If $\ell=3$, then this would lead to a different (shallower) recursive call than $\ell=2$:
				\begin{align*}
					\ceil{\frac{n'}{4}}	\le	\ceil{\frac{n'}{2}}	=	\ceil{\frac{\floor{2/2}n'}{2}}-\ceil{\frac{2}{2}}+1.
				\end{align*}
				We also need to ensure that if $k'\ge4$, then $n'\ge3$ or $n'\ge4$ as needed. Actually, if $k'\ge4$, and $n'<3$ ($n'<4$) then by Proposition \ref{pro:MaxDepth1} their depth is 1. And since the extra depth of using $\lambda=\left[\floor{\frac{k}{2}}^2\right]$ (with $0,1$ as needed), is\footnote{Recall that $k\ge4$ and $n\ge3$ ($n\ge4$). Actually, this is the only spot that needs $n\ge3$ when $k=$ even. Otherwise $n\ge2$ is sufficient.}
				\beq
					\ceil{\frac{\floor{k/2}n}{2}}-\ceil{\frac{k}{2}}	\ge	1,
				\eeq
				$\ell=2$ produces a recursive call that is at least as deep.
			\eep
			\label{lem:MaxDepth}
		\eel
		
		\begin{proposition}
			If $k=$ even, then for $n=1,2$, the maximum depth of recursive calls is $1$. If $k\ge3$ is odd, then for $n=1,2,3$, the maximum depth is $1$.
			\bep
				Consider a partition that repeats part $m$, $\ell\ge2$ times. First consider $n\le2$. The recursive call with $i=k-m$ is
				\begin{align*}
					n'	&=	mn-2k+2m+2\sum_{j=m+1}^{k}\left(j-m\right)d_{j}	\le	4m-2k+2\sum_{j=m+1}^{k}\left(j-m\right)d_{j}\\
						&\le	4m-2k+2(k-m\ell)	=	m(4-2\ell)	\le	m(4-4)	=	0.
				\end{align*}
				Thus, the chain terminates after at most 1 step. For all other partitions (distinct), we already know the chain terminates after 1 step.\\
				Now consider $k=$ odd and $n=3$. Note that $d_k=0$ otherwise we have a distinct partition. The recursive call with $i=k-1$ is
				\begin{align*}
					n'	&=	(k-(k-1))3-2(k-1)+2\sum_{j=0}^{(k-1)-1}\left((k-1)-j\right)d_{k-j}\\
						&=	5-2k+2\sum_{j=0}^{k-1}\left(k-j\right)d_{k-j}-2\sum_{j=0}^{k-1}d_{k-j}	=	5-2k+2k-2\sum_{j=0}^{k-1}d_{k-j}	=	5-2\sum_{j=1}^kd_j	<	0
				\end{align*}
				if $|\lambda|\ge3$. Thus, the only partitions that actually create recursive calls are those with $|\lambda|\le2$. And since $k=$ odd, $\lambda$ will be distinct and therefore terminate after at most 1 step.
			\eep
			\label{pro:MaxDepth1}
		\end{proposition}
		
		The explicit depth formula is {\sc Implemented as {\bf KOHdepthFAST(n,k)}}. It is somewhat odd that the depth is not symmetric in $n,k$. $G(n,k)$ is symmetric so we can choose to calculate $G(k,n)$ if $n>k$ and have a shallower depth of recursive calls. The total number of calls may be the same; this was not analyzed.
		
	\end{section}%

	\begin{section}{Altered Recurrence}

		\begin{subsection}{Restricted Partitions}
			We can alter the recurrence of Eqn.\ \eqref{eqn:KOH} in several ways to create $G'(n,k)$ and maintain the property that $G'(n,k)$ is symmetric and unimodal of darga $nk$.\footnote{Dependent on restrictions, one may achieve $G'(n,k)=0$.} One simple way to do this is to restrict the partitions over which we sum. This simply reduces the use of Proposition \ref{pro:SumSymUni}. We can restrict the minimum/maximum size of a partition, min/max integer in a partition, distinct parts, modulo classes, etc.
		
			Suppose we restrict to partitions with size $\le s$ and denote this new function as $G_s$.
			\begin{align*}
				G_s(n,k)	&=	\sum_{(d_1,\ldots,d_k);\sum_{i=1}^k id_i=k;\sum_{i=1}^k d_i\le s} q^{k\left(\sum_{i=1}^kd_i\right)-k-\sum_{1\le j<i\le k}(i-j)d_id_j}\\
						&\hspace{9em}	\prod_{i=0}^{k-1}G_s\left((k-i)n-2i+2\sum_{j=0}^{i-1}(i-j)d_{k-j},d_{k-i}\right)\\
						&=	\sum_{(d_1,\ldots,d_k);\sum_{i=1}^k id_i=k;\sum_{i=1}^k d_i\le s} q^{k\left(\sum_{i=1}^kd_i\right)-k-\sum_{1\le j<i\le k}(i-j)d_id_j}\\
						&\hspace{9em}	\prod_{i=0}^{k-1}G\left((k-i)n-2i+2\sum_{j=0}^{i-1}(i-j)d_{k-j},d_{k-i}\right).
			\end{align*}
			For $k\le s$, $G_s(n,k)=G(n,k)$ since ALL partitions of $k$ will have size $\le s$. And since $d_{k-i}\le s$ by definition of the sum, we can replace $G_s$ by $G$ (the actual $q$-binomial) in the product. I found conjectured recurrence relations of order 1 and degree 1 in both $n,k$ for $s=1,2,3,4$. I then conjectured for general $s$:
			\begin{conjecture}
				\begin{equation}
					q^n(q^{k+s}-1)G_s(n+1,k)-q^{k+1}(q^n-q^{s-2})G_s(n,k+1)+(q^n-q^{k+s-1})G_s(n+1,k+1)	=	0.
					\label{eqn:RecurrenceSizeS}
				\end{equation}
				\label{con:sizeS}
			\end{conjecture}
			Eqn.\ \eqref{eqn:RecurrenceSizeS} was verified for $n,k\le20$ and $s\le10$. The bounds were chosen to make the verification run in a short time ($\approx2$ minutes). In the algebra of formal power series, taking the limit as $s\to\infty$ (thus obtaining the original $G(n,k)$ with no restriction), recovers
			\begin{align*}
				q^n(-1)G_\infty(n+1,k)-q^{k+1}(q^n)G_\infty(n,k+1)+(q^n)G_\infty(n+1,k+1)	&=	0,\\
				G(n+1,k)+q^{k+1}G(n,k+1)	&=	G(n+1,k+1).
			\end{align*}
			which matches the recurrence relation in Eqn.\ \eqref{eqn:Recurrence}. Recall that $G_s(n,k)=G(n,k)$ for $k\le s$. It is somewhat surprising that $G(n,k)$ follows many recurrence relations for bounded $k$.
			
			We can try to enumerate all $G_s(n,k)$ using a translated Eqn.\ \eqref{eqn:RecurrenceSizeS}:
			\beq%
				G_s(n,k)	=	\frac{1}{q^{n}-q^{k+s-1}}\left(q^{n}(q^{k+s-1}-1)G_s(n,k-1)-q^{k}(q^{n}-q^{s-1})G_s(n-1,k)\right).
			\eeq%
			{\sc Implemented as \bf{Gs(s,n,k,q)}}. However, this recurrence by itself cannot enumerate all $G_s(n,k)$ as we will eventually hit the singularity-causing $n=k+s-1$ if $n\ge s$. However, we do obtain an interesting relation along that line:
			\beq
				G_s(k+s,k)	=	\frac{1-q^{k+1}}{1-q^{k+s}}G_s(k+s-1,k+1) \hspace{0.8em}\mbox{or}\hspace{0.8em}	G_{n-k+1}(n+1,k)	=	\frac{1-q^{k+1}}{1-q^{n+1}}G_{n-k+1}(n,k+1).
			\eeq
			So far, Eqn.\ \eqref{eqn:RecurrenceSizeS} has been verified for $s=1,2,3$ using the explicit formulas in Eqns. \eqref{eqn:ExplicitSize1}, \eqref{eqn:ExplicitSize2}, and \eqref{eqn:ExplicitSize3}, respectively.
			
			It is quite remarkable that the generalized recurrence in Eqn.\ \eqref{eqn:RecurrenceSizeS} yields polynomials using the initial conditions of Eqn.\ \eqref{eqn:InitialConditions},\footnote{$G_s(n,0)=G_s(0,k)=1$.} let alone unimodal polynomials. Are there other sets of initial conditions (potentially for chosen values of $s$) that will also yield unimodal polynomials when iterated in Eqn.\ \eqref{eqn:RecurrenceSizeS}?\footnote{At least until the singularity.}\\
			
			There are a couple of avenues for attack of Conjecture \ref{con:sizeS}. A good starting place may be to tackle the (easier?) question of whether the recurrence even produces symmetric polynomials.
			
			Since we have base case verification, we could try to use induction to prove the conjecture. This would likely involve writing $G_{s+1}=G_s+g_{s+1}$, where $g_{s+1}$ is the contribution from all of the partitions of size $s+1$:
			\begin{align*}
				g_{s+1}(n,k)	&=	\sum_{(d_1,\ldots,d_k);\sum_{i=1}^k id_i=k;\sum_{i=1}^k d_i=s+1} q^{ks-\sum_{1\le j<i\le k}(i-j)d_id_j}\\
							&\hspace{9em}	\prod_{i=0}^{k-1}G\left((k-i)n-2i+2\sum_{j=0}^{i-1}(i-j)d_{k-j},d_{k-i}\right).
			\end{align*}
			Another possibility, since we have an explicit expression for $G_s$ (and for $G$ if needed), may be to simply plug in to Eqn.\ \eqref{eqn:RecurrenceSizeS} and simplify cleverly.\\
			
			A more general question would be to analyze recurrences of type $$a(n,k)F(n+1,k+1)=b(n,k)F(n+1,k)+c(n,k)F(n,k+1)$$ that produce symmetric and unimodal polynomials. What are the requirements on $a,b,c$? If $a_1,b_1,c_1$ and $a_2,b_2,c_2$ both produce unimodal polynomials, what, if anything, can be said about $a_1+a_2$, $b_1+b_2$, $c_1+c_2$? About $a_1a_2$, $b_1b_2$, $c_1c_2$? If one can build from previously known valid recurrences, one could potentially build up to Eqn.\ \eqref{eqn:RecurrenceSizeS} from Eqn.\ \eqref{eqn:Recurrence}.			

			\begin{subsubsection}{\tbl Natural Partitions\tbr}
			
				\bel
					Suppose we restrict partitions to only use integers $\le p$. If $p<\frac{2k}{n+2}$ then $G'(n,k)=0$.\footnote{The bound $n+2<2\ceil{\frac{k}{p}}$ is sufficient to make $G'(n,k)=0$.}
					\bep
						We have $d_i=0$ for $i>p$. In the product of Eqn.\ \eqref{eqn:KOH}, consider $i=k-1$ for any restricted partition $\lambda$. Then for the new $n$ in the recursive call we obtain
						\begin{align*}
							n'(k-1)	&=	(k-(k-1))n-2(k-1)+2\sum_{j=0}^{(k-1)-1}((k-1)-j)d_{k-j}\\
									&=	n-2(k-1)+2\sum_{j=k-p}^{k-2}(k-1-j)d_{k-j}\\
																		&=	n-2(k-1)+2\sum_{j=2}^{p}(k-1-(k-j))d_{k-(k-j)}\\
																		&=	n-2(k-1)+2\sum_{j=1}^{p}(j-1)d_{j}\\
																		&=	n-2(k-1)+2k-2|\lambda|\\
																		&\le	n+2-2\ceil{\frac{k}{p}}\\
																		&\le	n+2-2\frac{k}{p}\\
																		&<	n+2-2\frac{k}{2k/(n+2)}	=	0.
						\end{align*}
						Then $G'(n',d_{1})=0$ for any partition and therefore $G'(n,k)=\sum_{\lambda}0=0$.
					\eep
					\label{lem:BoundedPartitions}
				\eel
				What do other natural partitions yield for a recursive call? Let us first consider $\lambda=[1^k]$. Then $d_1=k$ and $d_{i\ne1}=0$.
				\begin{align*}
					&q^{k\left(\sum_{i=1}^kd_i\right)-k-\sum_{1\le j<i\le k}(i-j)d_id_j}		\prod_{i=0}^{k-1}G\left((k-i)n-2i+2\sum_{j=0}^{i-1}(i-j)d_{k-j},d_{k-i}\right)\\
					&\tab=	q^{k\cdot k-k}G\left(n-2(k-1),k\right)\prod_{i=0}^{k-2}G\left((k-i)n-2i+2\sum_{j=0}^{i-1}(i-j)\cdot0,0\right)\\
					&\tab=	q^{k\cdot(k-1)}G\left(n-2(k-1),k\right)\prod_{i=0}^{k-2}G\left((k-i)n-2i,0\right)	=	q^{k(k-1)}G\left(n-2(k-1),k\right),
				\end{align*}
				since $(k-i)n-2i\ge2n-2(k-2)$ and if $2n-2(k-2)<0$, then $n-2(k-1)<0$. If we define $G'(n,k)$ as only being over partitions $=[1^k]$, then $G'(n,k)=q^{nk/2}$ if $2(k-1)|n$ and $=0$ otherwise. One can confirm the monomial by either looking at how many times $q^{k(k-1)}$ is multiplied or by recognizing that $G'(n,k)$ must be a monomial of darga $nk$.
				
				Now consider $\lambda=[k]$. Then $d_k=1$ and $d_{i\ne k}=0$.
				\begin{align*}
					&q^{k\left(\sum_{i=1}^kd_i\right)-k-\sum_{1\le j<i\le k}(i-j)d_id_j}	\prod_{i=0}^{k-1}G\left((k-i)n-2i+2\sum_{j=0}^{i-1}(i-j)d_{k-j},d_{k-i}\right)\\
					&\tab=	q^{k-k}G\left(kn,1\right)\prod_{i=1}^{k-1}G\left((k-i)n-2i+2\sum_{j=1}^{i-1}(i-j)d_{k-j}+2(i-0)d_k,0\right)\\
					&\tab=	G\left(kn,1\right)\prod_{i=1}^{k-1}G\left((k-i)n,0\right)	=	G(nk,1)	=	\frac{1-q^{nk+1}}{1-q},
				\end{align*}
				since $(k-i)n\ge n\ge0$. So $[k]$ provides the \tbl base\tbr\ of $G(n,k)$: $1+q+\cdots+q^{nk}$.
				
				If we restrict to partitions of size 1, then
				\begin{equation}
					G_1(n,k)	=	\frac{1-q^{nk+1}}{1-q}.
					\label{eqn:ExplicitSize1}
				\end{equation}
				We can use this to verify the recurrence relation in Eqn.\ \eqref{eqn:RecurrenceSizeS} for $s=1$.\\
				
				\bel
					Let us restrict partitions to those with distinct parts. Fix $k\in\N$. Let $\ell=\floor{\frac{1}{2}(\sqrt{8k+1}-1)}$.\footnote{I.e., the index of the largest triangular number $\le k$ or the maximum size of a distinct partition of $k$.} Then $G'(n,k)$ is a polynomial of degree $k$ in $q^n$ for $n\ge2\ell-2$.\footnote{It appears that $n=2\ell-3$ also fits the polynomial, though this was left unproven.}%
					\bep
						Since $\lambda=[k]$ is distinct, $G'(n,k)$ is of degree $k$ in $q^{n}$. Let $\lambda$ be a partition with distinct parts. Then $d_{i}\le1$ and $|\lambda|=\sum_{j=1}^kd_k\le l$. And for each recursive call, $G'(n',d_i)=0,1,$ or $\frac{1-q^{n'+1}}{1-q}$. Since $k$ is fixed, it has a fixed number of partitions into distinct parts; thus, it remains to show $G'(n',d_i)\ne0$ by showing $n'(i)\ge0$ for all $i$.
						\begin{align*}
							n'(i)			&=	(k-i)n-2i+2\sum_{j=0}^{i-1}(i-j)d_{k-j},\\
							n'(i)-n'(i-1)	&=	2\sum_{j=0}^{i-1}(i-j)d_{k-j}-2\sum_{j=0}^{i-1}(i-j)d_{k-j}+2\sum_{j=0}^{i-2}d_{k-j}+2d_{k-i+1}-n-2\\
											&=	2\sum_{j=0}^{i-1} d_{k-j}-n-2	\le	2\ell-n-2	\le	0.
						\end{align*}
						Thus, the minimum occurs at $i=k-1$ and has value
						\begin{align*}
							n'(k-1)	&=	(k-(k-1))n-2(k-1)+2\sum_{j=0}^{k-1-1}(k-1-j)d_{k-j}\\
									&=	n-2(k-1)+2\sum_{j=0}^{k-1}(k-1-j)d_{k-j}\\
									&=	n-2(k-1)+2\sum_{j=0}^{k-1}(k-j)d_{k-j}-2\sum_{j=0}^{k-1}d_{k-j}\\
									&=	n-2(k-1)+2k-2\sum_{j=0}^{k-1}d_{k-j}	=	n+2-2\ell	\ge	0.
						\end{align*}
					\eep
				\eel
				
				Let us return to the minimum of $n'$ for general partitions. It occurs either at $i=k-1$, or if $n\le2|\lambda|-2-2(\#$ of 1s in $\lambda)$, then at the final $i$ s.t. $n+2\ge2\sum_{j=0}^{i-1} d_{k-j}$: see the simplification of $n'(i)-n'(i-1)$ above. The second difference is $n'(i)-2n'(i-1)+n'(i-2)=2d_{k-i+1}\ge0$. Thus, the extrema of $n'$ found is in fact a minimum. And if we fix the $i$ such that there exists a minimum (or equivalent to the minimum), we obtain a value of
				\begin{align*}
					(k-i)n-2i+2\sum_{j=0}^{i-1}(i-j)d_{k-j}	&\approx	(k-i)\left(2\sum_{j=0}^{i-1} d_{k-j}-2\right)-2i+2\sum_{j=0}^{i-1}(i-j)d_{k-j}\\
												&=	2\sum_{j=0}^{i-1}(k-j)d_{k-j}-2k	\le	2k-2k	=	0.
				\end{align*}
				To avoid the approximation step, consider the next $i$:
				\begin{align*}
					(k-(i+1))n-&2(i+1)+2\sum_{j=0}^{i}(i+1-j)d_{k-j}\\
									&<	(k-(i+1))\left(2\sum_{j=0}^{i} d_{k-j}-2\right)-2(i+1)+2\sum_{j=0}^{i}(i+1-j)d_{k-j}\\
														&=	2\sum_{j=0}^{i}(k-j) d_{k-j}-2k	\le	2k-2k	=	0.
				\end{align*}
				The recursive call is then seen to be 0 for any partition that has $|\lambda|-|\{1s\mbox{ in }\lambda\}|\ge \frac{n+2}{2}$.
				
				Now we consider $\lambda=[k-\ell,\ell]$ for some $1\le \ell<\frac{k}{2}$.
				\begin{align*}
					&q^{k\left(\sum_{i=1}^kd_i\right)-k-\sum_{1\le j<i\le k}(i-j)d_id_j}	\prod_{i=0}^{k-1}G\left((k-i)n-2i+2\sum_{j=0}^{i-1}(i-j)d_{k-j},d_{k-i}\right)\\
					&\tab=	q^{2k-k-((k-\ell)-\ell)}\prod_{i=0}^{\ell-1}G\left((k-i)n-2i+2\sum_{j=0}^{i-1}(i-j)d_{k-j},0\right)\\
					&\tab\tab	G\left((k-\ell)n-2\ell+2\sum_{j=0}^{\ell-1}(\ell-j)d_{k-j},d_{k-\ell}\right)\\
					&\tab\tab	\prod_{i=\ell+1}^{k-\ell-1}G\left((k-i)n-2i+2\sum_{j=0}^{i-1}(i-j)d_{k-j},0\right)\\
					&\tab\tab	G\left((k-(k-\ell))n-2(k-\ell)+2\sum_{j=0}^{k-\ell-1}(k-\ell-j)d_{k-j},d_{\ell}\right)\\
					&\tab\tab	\prod_{i=k-\ell+1}^{k-1}G\left((k-i)n-2i+2\sum_{j=0}^{i-1}(i-j)d_{k-j},0\right)\\
					&\tab=	q^{2\ell}\prod_{i=0}^{\ell-1}G\left((k-i)n-2i,0\right)\\
					&\tab\tab	G\left((k-\ell)n-2\ell,1\right)\\
					&\tab\tab	\prod_{i=\ell+1}^{k-\ell-1}G\left((k-i)n-2\ell,0\right)\\
					&\tab\tab	G\left(\ell n-2\ell,1\right)\\
					&\tab\tab	\prod_{i=k-\ell+1}^{k-1}G\left((k-i)n+2i-2k,0\right),
				\end{align*}
				which to be nonzero requires that (assuming $n\ge2$ otherwise we get $0$ in the recursive call)
				\begin{align*}
					0	&\le	(k-i)n-2i	\ge	(k-\ell+1)n-2(\ell-1)	\ge	(2\ell-\ell+1)2-2(\ell-1)	\ge	4,	\tab	\mbox{and}\\
					0	&\le	(k-i)n-2\ell	\ge	(k-(k-\ell-1))n-2\ell	=	(\ell+1)n-2\ell	\ge	2,	\tab	\mbox{and}\\
					0	&\le	(k-i)n+2i-2k	\ge	(k-(k-1))n+2(k-1)-2k	=	n-2	\ge	0.
				\end{align*}
				The leftmost inequality is what we need and the other inequalities deduce that to be true. Thus, the recursive call is
				\beq
					q^{2\ell}\frac{1-q^{(k-\ell)n-2\ell+1}}{1-q}\cdot\frac{1-q^{\ell n-2\ell+1}}{1-q}.
				\eeq
				
				For a visual of how the different partitions can contribute to the final polynomial, see Figure \ref{fig:PartitionDivision}.
				\begin{figure}[H]
					\caption{Contributions from Different Partitions to Produce $G(5,5)$.}
					\includegraphics[width=\textwidth]{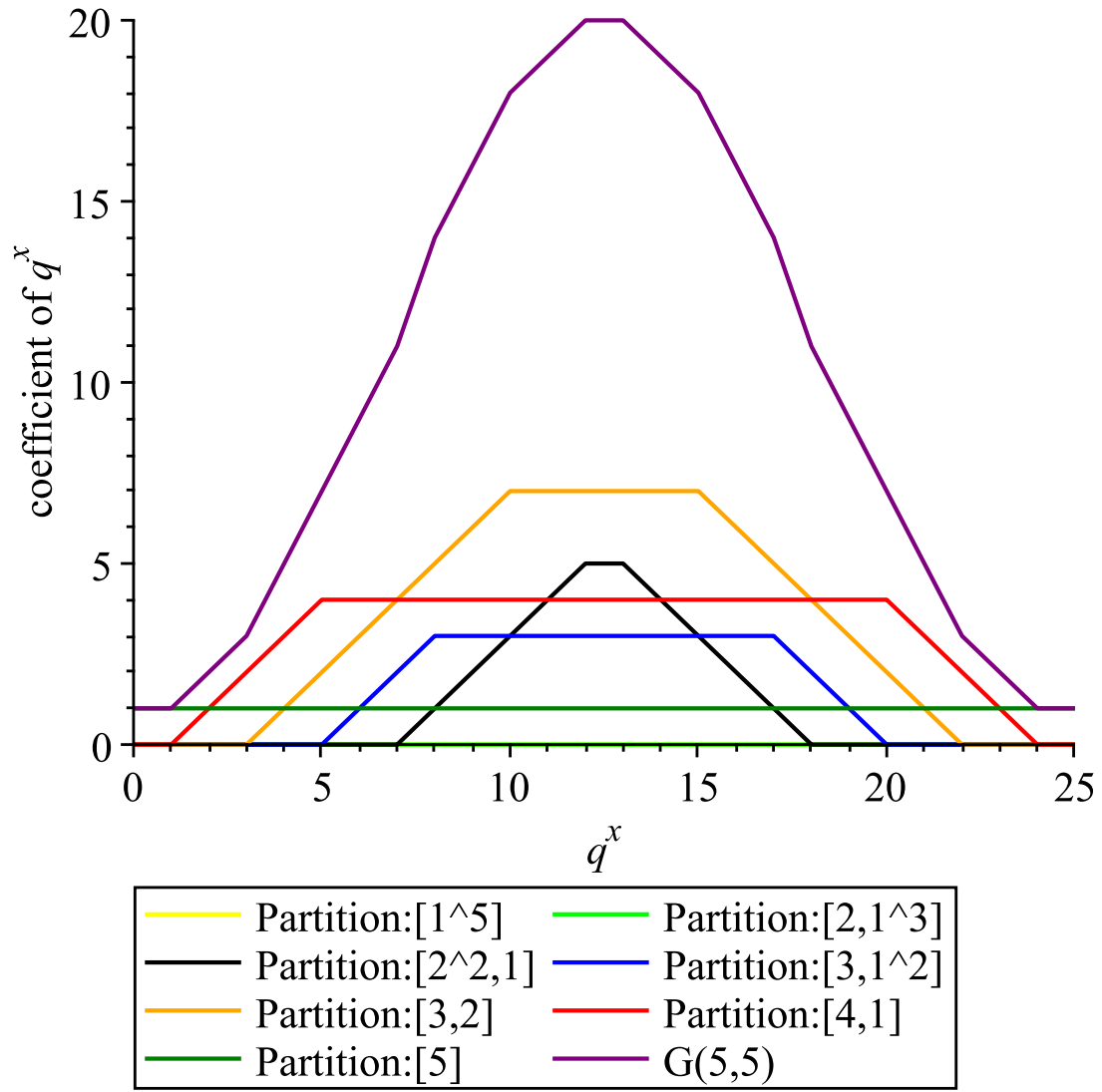}
					 The partitions $[1^5]$ and $[2,1^3]$ both contribute the 0 polynomial.
					\label{fig:PartitionDivision}
				\end{figure}
	
				\bel
					If we restrict to partitions with consecutive difference $\ge k-d$, and $d<\frac{2}{3}k+1$, then the coefficient of $q^{nk/2}$ in $G'(n,k)$ is $n\frac{D(D+1)}{2}-(D^2-1)$ for $D=\floor{\frac{d}{2}}$.
					\bep
						The bound on $d$ implies that the maximum size of a partition is $2$ by Proposition \ref{pro:RestrictMinDiffSize}. $D$ can be seen to be the relevant number because it is the number of possible size 2 partitions of $k$. $[k],[k-1,1],\ldots,[k-D,D]$ are valid since $k-D-D=k-2\floor{\frac{d}{2}}\ge k-d$. Then
						\begin{equation}
							G'(n,k)	=	\frac{1-q^{nk+1}}{1-q}+\sum_{i=1}^D q^{2i}\frac{1-q^{(k-i)n-2i+1}}{1-q}\cdot\frac{1-q^{in-2i+1}}{1-q}.
							\label{eqn:PartitionsDifferenceKD}
						\end{equation}
						Note that the middle coefficient of $\left(\frac{1-q^{a}}{1-q}\right)\cdot\left(\frac{1-q^{b}}{1-q}\right)$ is going to be $\min(a,b)$. Also note that the darga of the summand in Eqn.\ \eqref{eqn:PartitionsDifferenceKD} is
						\beq
							2\cdot(2i)+[(k-i)n-2i+1-1]+[in-2i+1-1]	=	nk,
						\eeq
						so $q^{nk/2}$ is the middle monomial.\footnote{As one should expect since the summand of Eqn.\ \eqref{eqn:PartitionsDifferenceKD} is a specific case of the summand in Eqn.\ \eqref{eqn:KOH}.} Thus, the coefficient is
						\beq
							1+\sum_{i=1}^D \min\left[(k-i)n-2i+1,in-2i+1\right]	=	1+\sum_{i=1}^D (in-2i+1)	=	\frac{1}{2}\left(  \left( n-2 \right)D +2 \right)  \left( D+1 \right).
						\eeq
					\eep
					\label{lem:CoefficientSize2Partition}
				\eel
				
				\begin{proposition}
					If we restrict to partitions with difference $\ge k-d$ and minimum part $\ell$, then to have a partition of size $s$, it is required that
					\begin{align*}
						d	&\ge	k\frac{(s+1)(s-2)}{s(s-1)}+\frac{2\ell}{s-1},\\
						\ell	&\le	\frac{ds^2-ks^2-ds+ks+2k}{2s},\\
						k	&\le	\frac{(ds-d-2\ell)s}{(s+1)(s-2)}.
					\end{align*}
					\bep
						A partition $\lambda$ of size $s$, minimum part $\ell$, and difference $k-d$ will have total
						\beq
							k	=	\sum_{i=1}^s\lambda_i	\ge	\sum_{i=1}^s \left[\ell+(i-1)(k-d)\right]	=	\frac{s}{2}(2\ell-k+d+(k-d)s).
						\eeq
					\eep
					\label{pro:RestrictMinDiffSize}
				\end{proposition}
	
			\end{subsubsection}%

			\begin{subsubsection}{Size 2 Partitions}
			
				Another special partition to consider is $\lambda=\left[\frac{k}{2}^2\right]$ when $k$ is even. Then $d_{k/2}=2$ and $d_i=0$ otherwise.
				\begin{align*}
					&q^{k\left(\sum_{i=1}^kd_i\right)-k-\sum_{1\le j<i\le k}(i-j)d_id_j}	\prod_{i=0}^{k-1}G\left((k-i)n-2i+2\sum_{j=0}^{i-1}(i-j)d_{k-j},d_{k-i}\right)\\
						&\tab=	q^{2k-k-0}\prod_{i=0}^{k/2-1}G\left((k-i)n-2i+2\sum_{j=0}^{i-1}(i-j)d_{k-j},d_{k-i}\right)\\
						&\tab\tab	G\left(\left(k-\frac{k}{2}\right)n-2\frac{k}{2}+2\sum_{j=0}^{k/2-1}\left(\frac{k}{2}-j\right)d_{k-j},d_{k/2}\right)\\
						&\tab\tab	\prod_{i=k/2+1}^{k-1}G\left((k-i)n-2i+2\sum_{j=0}^{i-1}(i-j)d_{k-j},d_{k-i}\right)\\
						&\tab=	q^{k}G\left(\frac{k}{2}n-k,2\right)	\prod_{i=0}^{k/2-1}G\left((k-i)n-2i,0\right)	\prod_{i=k/2+1}^{k-1}G\left((k-i)n-2i+2(i-k/2)2,0\right)\\
						&\tab=	q^{k}G\left(\frac{k}{2}n-k,2\right)	\prod_{i=0}^{k/2-1}G\left((k-i)n-2i,0\right)	\prod_{i=k/2+1}^{k-1}G\left((k-i)n+2i-2k,0\right),
				\end{align*}
				which is nonzero only if (assuming $n\ge2$ otherwise we get $0$ in the recursive call)
				\begin{align*}
					0	&\le	(k-i)n-2i	\ge	(k-(k/2-1))2-2(k/2-1)	=	(k/2+1)2-k+2	=	4,\\
					0	&\le	(k-i)n+2i-2k	\ge	(k-(k-1))n+2(k-1)-2k	=	n-2	\ge	0.
				\end{align*}
				Thus, the recursive call is $q^{k}G\left(\frac{k}{2}n-k,2\right)$.\\
				
				While proving Lemma \ref{lem:CoefficientSize2Partition}, Maple produced this useful gem:
				\begin{align}
					\sum_{i=1}^D q^{2i}\frac{1-q^{(k-i)n-2i+1}}{1-q}\ \cdot\ \frac{1-q^{in-2i+1}}{1-q}	&=	q^2\frac{(1+q^{nk-2D})(1-q^{2D})}{(1-q)^2(1-q^2)}\nonumber\\
																				&\tab	-q^{n+1}\frac{(1+q^{nk-nD-n})(1-q^{nD})}{(1-q)^2(1-q^n)}.
					\label{eqn:UsefulGem}
				\end{align}
				We can then look at using only partitions restricted to maximum size $2$. The possible partitions are $[k]$ and $[k-\ell,\ell]$ for $1\le\ell\le\frac{k}{2}$. We can find the exact form of $G_2(n,k)$ by including all partitions with size $\le2$.
				\bel
					For $n\ge1$,
					\begin{align}
						G_2(n,k)	&=	\bigg(q^{nk+n+4}-q^{nk+n+3}+q^{nk+n+1}-q^{nk+4}-q^{(n-1)k+n+3}+q^{(n-1)k+3}+q^{k+n+1}\nonumber\\
								&\tab	-q^{k+1}-q^n+q^3-q+1\bigg)\bigg/(1-q)^2(1-q^2)(1-q^n).
						\label{eqn:ExplicitSize2}
					\end{align}
					If $n=0$, then $G_2(0,k)=1$.
					\bep
						If $k=$ odd, then, utilizing Eqn.\ \eqref{eqn:UsefulGem},
						\begin{align*}
							G_2(n,k)	&=	\frac{1-q^{nk+1}}{1-q}+\sum_{\ell=1}^{(k-1)/2} q^{2\ell}\frac{1-q^{(k-\ell)n-2\ell+1}}{1-q}\cdot\frac{1-q^{\ell n-2\ell+1}}{1-q}\\
									&=	\frac{1-q^{nk+1}}{1-q}+q^2\frac{(1+q^{nk-2(k-1)/2})(1-q^{2(k-1)/2})}{(1-q)^2(1-q^2)}\\
									&\tab	-q^{n+1}\frac{(1+q^{nk-n(k-1)/2-n})(1-q^{n(k-1)/2})}{(1-q)^2(1-q^n)}\\
									&=	\frac{1-q^{nk+1}}{1-q}+q^2\frac{(1+q^{nk-k+1})(1-q^{k-1})}{(1-q)^2(1-q^2)}-q^{n+1}\frac{1-q^{n(k-1)}}{(1-q)^2(1-q^n)}.
						\end{align*}
						If $k=$ even, then
						\begin{align*}
							G_2(n,k)	&=	\frac{1-q^{nk+1}}{1-q}+\sum_{\ell=1}^{k/2-1} q^{2\ell}\frac{1-q^{(k-\ell)n-2\ell+1}}{1-q}\cdot\frac{1-q^{\ell n-2\ell+1}}{1-q}+q^{k}G_2\left(\frac{k}{2}n-k,2\right)\\
									&=	\frac{1-q^{nk+1}}{1-q}+q^2\frac{(1+q^{nk-k+2})(1-q^{k-2})}{(1-q)^2(1-q^2)}-q^{n+1}\frac{(1+q^{nk/2})(1-q^{n(k-2)/2})}{(1-q)^2(1-q^n)}\\
									&\tab	+q^k\frac{(1-q^{nk/2-k+1})(1-q^{nk/2-k+2})}{(1-q)(1-q^2)}.
						\end{align*}
						We again utilized Eqn.\ \eqref{eqn:UsefulGem} and noted that $G_2(n',2)=G(n',2)$ and used Eqn.\ \eqref{eqn:Qbinomial}. One can use Maple or any other mathematics software to verify that both cases \tbl simplify\tbr\ to the expression given above.
						
						Notice also that the bound of $n\ge2$ required from characterizing partition calls has been lowered to $n\ge1$. If $n=1$, then the size 2 partition calls will be 0, which is exactly what the extra summands reduce to.
					\eep
					\label{lem:ExplicitSize2}
				\eel
				We can use this to verify the recurrence relation in Eqn.\ \eqref{eqn:RecurrenceSizeS} for $s=2$. To compute $G_2(n,k)$ using the explicit expression, use {\bf G2(n,k,q)}.\\
				
				The expression given in Eqn.\ \eqref{eqn:ExplicitSize2} obfuscates that $G_2(n,k)$ is unimodal of darga $nk$. The expanded expressions in the proof allow for a heuristic argument by reasoning the darga of each summand. This would be a great way to prove unimodality of general functions. However, the useful decomposition is often difficult to ascertain.

			\end{subsubsection}%

			\begin{subsubsection}{Size 3 Partitions}
			
				By looking at $\lambda=[\ell_1,\ell_2,\ell_3]\vdash k$ for the separate cases $\ell_1=\ell_2=\ell_3$, $\ell_1=\ell_2\ne\ell_3$, $\ell_1\ne\ell_2=\ell_3$, and $\ell_1\ne\ell_2\ne\ell_3$, one can add this to Lemma \ref{lem:ExplicitSize2}, and find an explicit expression for $G_3(n,k)$.\\
				
				We start with $[\ell_1\ne\ell_2\ne\ell_3]$. Then $d_{\ell_i}=1$. We will need $n\ge4$ for nonzero calls.
				\begin{align*}
					&q^{k\left(\sum_{i=1}^kd_i\right)-k-\sum_{1\le j<i\le k}(i-j)d_id_j}	\prod_{i=0}^{k-1}G\left((k-i)n-2i+2\sum_{j=0}^{i-1}(i-j)d_{k-j},d_{k-i}\right)\\
						&\tab=	q^{2(\ell_2+2\ell_3)}\prod_{i=0}^{k-\ell_1-1}G\left((k-i)n-2i,0\right)\\
						&\tab\tab	G\left(\ell_1n-2(\ell_2+\ell_3),1\right)\\
						&\tab\tab	\prod_{i=k-\ell_1+1}^{k-\ell_2-1}G\left((k-i)n-2(k-\ell_1),0\right)\\
						&\tab\tab	G\left(\ell_2n-2(\ell_2+\ell_3),1\right)\\
						&\tab\tab	\prod_{i=k-\ell_2+1}^{k-\ell_3-1}G\left((k-i)(n-2)-2\ell_3,0\right)\\
						&\tab\tab	G\left(\ell_3(n-4),1\right)\\
						&\tab\tab	\prod_{i=k-\ell_3+1}^{k-1}G\left((k-i)(n-4),0\right)\\
						&\tab=	q^{2(\ell_2+2\ell_3)}\cdot\frac{1-q^{\ell_1n-2(\ell_2+\ell_3)+1}}{1-q}\cdot\frac{1-q^{\ell_2n-2(\ell_2+\ell_3)+1}}{1-q}\cdot\frac{1-q^{\ell_3(n-4)+1}}{1-q}.
				\end{align*}

				Now $[k-2\ell,\ell^2]$, for $\ell<\frac{k}{3}$. Then $d_{k-2\ell}=1$, $d_{\ell}=2$. Again, we need $n\ge4$.
				\begin{align*}
					&q^{k\left(\sum_{i=1}^kd_i\right)-k-\sum_{1\le j<i\le k}(i-j)d_id_j}	\prod_{i=0}^{k-1}G\left((k-i)n-2i+2\sum_{j=0}^{i-1}(i-j)d_{k-j},d_{k-i}\right)\\
						&\tab=	q^{6\ell}\prod_{i=0}^{2\ell-1}G\left((k-i)n-2i,0\right)\\
						&\tab\tab	G\left((k-2\ell)n-4\ell,1\right)\\
						&\tab\tab	\prod_{i=2\ell+1}^{k-\ell-1}G\left((k-i)n-4\ell,0\right)\\
						&\tab\tab	G\left(\ell n-4\ell,2\right)\\
						&\tab\tab	\prod_{i=k-\ell+1}^{k-1}G\left((k-i)(n-4),0\right)\\
						&\tab=	q^{6\ell}\cdot\frac{1-q^{(k-2\ell)n-4\ell+1}}{1-q}\cdot\frac{(1-q^{\ell n-4\ell+1})(1-q^{\ell n-4\ell+2})}{(1-q)(1-q^2)}.
				\end{align*}

				We continue by analyzing $[\ell^2,k-2\ell]$, for $\frac{k}{2}>\ell>\frac{k}{3}$. Then $d_{\ell}=2$, $d_{k-2\ell}=1$. $n\ge4$ is needed as usual.
				\begin{align*}
					&q^{k\left(\sum_{i=1}^kd_i\right)-k-\sum_{1\le j<i\le k}(i-j)d_id_j}	\prod_{i=0}^{k-1}G\left((k-i)n-2i+2\sum_{j=0}^{i-1}(i-j)d_{k-j},d_{k-i}\right)\\
						&\tab=	q^{4k-6\ell}\prod_{i=0}^{k-\ell-1}G\left((k-i)n-2i,0\right)\\
						&\tab\tab	G\left(\ell n-2k+2\ell,2\right)\\
						&\tab\tab	\prod_{i=k-\ell+1}^{2\ell-1}G\left((k-i)n+2i-4(k-\ell),0\right)\\
						&\tab\tab	G\left((k-2\ell)n+8\ell-4k,1\right)\\
						&\tab\tab	\prod_{i=2\ell+1}^{k-1}G\left((k-i)(n-4),0\right)\\
						&\tab=	q^{4k-6\ell}\cdot\frac{(1-q^{\ell n-2k+2\ell+1})(1-q^{\ell n-2k+2\ell+2})}{(1-q)(1-q^2)}\cdot\frac{1-q^{(k-2\ell)n+8\ell-4k+1}}{1-q}.
				\end{align*}

				Finally, consider $[\frac{k}{3}^3]$ for $k\cong0\modu{3}$. Then $d_{k/3}=3$. Again, we need $n\ge4$.
				\begin{align*}
					&q^{k\left(\sum_{i=1}^kd_i\right)-k-\sum_{1\le j<i\le k}(i-j)d_id_j}	\prod_{i=0}^{k-1}G\left((k-i)n-2i+2\sum_{j=0}^{i-1}(i-j)d_{k-j},d_{k-i}\right)\\
						&\tab=	q^{2k}\prod_{i=0}^{2k/3-1}G\left((k-i)n-2i,0\right)\\
						&\tab\tab	G\left((n-4)k/3,3\right)\\
						&\tab\tab	\prod_{i=2k/3+1}^{k-1}G\left((k-i)(n-4),0\right)\\
						&\tab=	q^{2k}\frac{(1-q^{(n-4)k/3+1})(1-q^{(n-4)k/3+2})(1-q^{(n-4)k/3+3})}{(1-q)(1-q^2)(1-q^3)}.
				\end{align*}

				\bel
					For $n\ge2$, 
					\begin{align}
						G_3(n,k)	&=	\bigg[q^{4k}\cdot{q}^{3} \left( 1-{q}^{n} \right)  \left( q-{q}^{n} \right)%
									\tab	-	\tab	q^{3k}\cdot q\left( 1+q \right)  \left( 1-{q}^{n} \right)  \left( q-{q}^{2}+{q}^{5}-{q}^{n} \right)\nonumber\\
									&\hspace{1.2em}-q^{2k}\bigg(q^{nk}\left({q}^{2n}({q}^{9}- {q}^{8}-{q}^{7}+{q}^{6}+ {q}^{5}-{q}^{3}+q)-{q}^{n}({q}^{10}-{q}^{8}+{q}^{6}+{q}^{5})+{q}^{10}\right)\nonumber\\
									&\tab\tab\tab-q^{2n}+{q}^{n}({q}^{5}+{q}^{4}-{q}^{2}+1)-{q}^{9}+{q}^{7}-{q}^{5}-{q}^{4}+{q}^{3}+{q}^{2}-q\bigg)\nonumber\\
									&\hspace{1.2em}+q^k\cdot {q}^{n k}{q}^{3} \left( 1+q \right) \left( 1-{q}^{n} \right)  \left( q^5-q^n+q^{n+3}-{q}^{n+4}\right)%
										-	{q}^{n k}{q}^{6} \left( 1-{q}^{n} \right) \left( q-{q}^{n} \right)\bigg]\nonumber\\
									&\tab\bigg/(1-q)^2(1-q^2)^2(1-q^3)(1-q^{n-1})(1-q^{n})q^{2k+1}.
						\label{eqn:ExplicitSize3}
					\end{align}
					For $n<2$, $G_3(n,k)=G_2(n,k)$.
					\bep
						For $n\ge4$,
						\begin{align*}
							G_3(n,k)	&=	G_2(n,k)+\sum_{\ell_3=1}^{\floor{k/3}-1}\sum_{\ell_2=\ell_3+1}^{\floor{(k-\ell_3-1)/2}}q^{2(\ell_2+2\ell_3)}\cdot\frac{1-q^{(k-\ell_3-\ell_2)n-2(\ell_2+\ell_3)+1}}{1-q}\\
									&\hspace{6.95cm}	\cdot\frac{1-q^{\ell_2n-2(\ell_2+\ell_3)+1}}{1-q}\cdot\frac{1-q^{\ell_3(n-4)+1}}{1-q}\\
									&\tab	+	\sum_{\ell=1}^{\floor{(k-1)/3}}q^{6\ell}\cdot\frac{1-q^{(k-2\ell)n-4\ell+1}}{1-q}\cdot\frac{(1-q^{\ell n-4\ell+1})(1-q^{\ell n-4\ell+2})}{(1-q)(1-q^2)}\\
									&\tab	+	\sum_{\ell=\ceil{(k+1)/3}}^{\floor{(k-1)/2}}q^{4k-6\ell}\hspace{-0.2em}\cdot\hspace{-0.2em}\frac{(1-q^{\ell n-2k+2\ell+1})(1-q^{\ell n-2k+2\ell+2})}{(1-q)(1-q^2)}\hspace{-0.2em}\cdot\hspace{-0.2em}\frac{1-q^{(k-2\ell)n+8\ell-4k+1}}{1-q}.
						\end{align*}
						The first nested summation is for the distinct size 3 partitions. The other 2 sums are for partitions of type $[k-2\ell,\ell,\ell]$ and $[\ell,\ell,k-2\ell]$, respectively. If $k\cong0\modu{3}$ then, for the $\left[\frac{k}{3}^3\right]$ partition, we need to add the extra term:
						\beq
							q^{2k}\frac{(1-q^{(n-4)k/3+1})(1-q^{(n-4)k/3+2})(1-q^{(n-4)k/3+3})}{(1-q)(1-q^2)(1-q^3)}.
						\eeq
						Amazingly, Maple is able to simplify the above expression into closed form once one makes assumptions about $k\modu{6}$. It is also necessary to divide the first sum into 2 parts based on the parity of $l_3$. In the divided sum, we are replacing $\ell_3$ by $2\ell_3-1$ and $2\ell_3$, respectively.%
						\begin{align*}
							\sum_{\ell_3=1}^{\floor{(\floor{k/3}-1)/2}}\Bigg[&\sum_{\ell_2=2\ell_3}^{\floor{k/2}-\ell_3}q^{2(\ell_2+2(2\ell_3-1))}\cdot\frac{1-q^{(k-(2\ell_3-1)-\ell_2)n-2(\ell_2+2\ell_3)+3}}{1-q}\\
							&\hspace{3.62cm}	\cdot\frac{1-q^{\ell_2 n-2(\ell_2+2\ell_3)+3}}{1-q}\cdot\frac{1-q^{(2\ell_3-1)(n-4)+1}}{1-q}\\
										&\tab+\sum_{\ell_2=2\ell_3+1}^{\floor{(k-1)/2}-\ell_3}q^{2(\ell_2+4\ell_3)}\cdot\frac{1-q^{(k-2\ell_3-\ell_2)n-2(\ell_2+2\ell_3)+1}}{1-q}\\
							&\hspace{4.62cm}	\cdot\frac{1-q^{\ell_2n-2(\ell_2+2\ell_3)+1}}{1-q}\cdot\frac{1-q^{2\ell_3(n-4)+1}}{1-q}\Bigg].
						\end{align*}
						If $k\cong0,1,2\modu{6}$ then $\floor{\frac{k}{3}}-1\cong1\modu{2}$ and so we need to add the final inner sum when $\ell_3=\floor{\frac{k}{3}}-1$:
						\begin{align*}
							\sum_{\ell_2=\floor{k/3}}^{\floor{(k-\floor{k/3})/2}}q^{2(\ell_2+2(\floor{k/3}-1))}&\cdot\frac{1-q^{(k-\ell_2-\floor{k/3}+1)n-2(\ell_2+\floor{k/3})+3}}{1-q}\\
							&\cdot\frac{1-q^{\ell_2n-2(\ell_2+\floor{k/3})+3}}{1-q}\cdot\frac{1-q^{(\floor{k/3}-1)(n-4)+1}}{1-q}.
						\end{align*}
						For each $k\cong0,\ldots,5\modu{6}$, the decomposition simplifies to the same result: Eqn.\ \eqref{eqn:ExplicitSize3}.\footnote{This was only fully done for the $k\cong5\modu{6}$ case. However, the other cases were partly simplified and found to be in agreement empirically; after removing the assumption on $k\modu{6}$, all cases still gave the same result. And the explicit formula in Eqn.\ \eqref{eqn:ExplicitSize3} matches the polynomials found by recursion for all $n=2,\ldots,30$ and $k=0,\ldots30$. If the reader does not wish to take our word for it, they are encouraged to show the other cases for themselves. We highly recommend using Maple or another computer assistant if you do not have time to waste.} 
						
						If $n<4$, then the additions from size 3 partition calls will all be 0 (see the characterization of size 3 partition calls above) and so $G_3(n,k)=G_2(n,k)$. It is a simple exercise to check that the difference evaluates to 0 when $n=2,3$. The difference between the two polynomials is given by
						\begin{align*}
							G_3(n,k)-G_2(n,k)	&=	\bigg[q^{4k}\cdot q^2\left( 1-{q}^{n} \right)  \left( q-{q}^{n} \right)\\
											&\tab	-q^{3k}\cdot q\left( 1+q \right)  \left( 1-{q}^{n} \right)  \left( {q}^{n+3}-{q}^{n+2}-q^n+{q}^{3}\right)\\
											&\tab	+q^{2k}\hspace{-0.13em}\cdot\hspace{-0.13em} q\left( {q}^{3}-{q}^{n} \right)  \left({q}^{nk}\left({q}^{n+1}-{q}^{3}-{q}^{2}+1\right)+q^n(q^3-q-1)+{q}^{2} \right)\\
											&\tab	-q^{k}\cdot{q}^{nk}{q}^{3} \left( 1+q \right)  \left( 1-{q}^{n} \right)  \left(1-q -{q}^{3}+{q}^{n}\right)\\
											&\tab	-{q}^{nk}{q}^{5} \left( 1-{q}^{n} \right)  \left( q-{q}^{n} \right)\bigg]\\
											&\tab\bigg/(1-q)^2(1-q^2)^2(1-q^3)(1-q^{n-1})(1-q^{n})q^{2k}.
						\end{align*}
					\eep
					\label{lem:ExplicitSize3}
				\eel
				We can use this result to verify the recurrence relation in Eqn.\ \eqref{eqn:RecurrenceSizeS} for $s=3$. To compute $G_3(n,k)$ for numeric $n,k$ using the explicit expression, use {\bf G3(n,k,q)}.
			
			\end{subsubsection}%
			
			Because we will always have $k'\le s$, it is possible to automate the process of finding $G_s(n,k)$ for any $s$ and have it terminate in finite time. Characterize all of the new size $s$ partitions and add their contributions to $G_{s-1}(n,k)$. However, the number of new types of partitions to consider grows exponentially as $2^{s-1}$; each new part is either the same as the previous part, or smaller. Also, the partitions become more and more complex to describe. One general partition that we can tackle is $\left[\frac{k}{\ell}^\ell\right]$. Then $d_{k/\ell}=\ell$ and the recursive call is
			\begin{align*}
				q^{k\left(\sum_{i=1}^kd_i\right)-k-\sum_{1\le j<i\le k}(i-j)d_id_j}	&	\prod_{i=0}^{k-1}G\left((k-i)n-2i+2\sum_{j=0}^{i-1}(i-j)d_{k-j},d_{k-i}\right)\\
							&=	q^{k(\ell-1)}\prod_{i=0}^{k(\ell-1)/\ell-1}G\left((k-i)n-2i,0\right)\\
							&\tab	G\left(\frac{k}{\ell}n-2k\frac{\ell-1}{\ell},\ell\right)\\
							&\tab	\prod_{i=k(\ell-1)/\ell+1}^{k-1}G\left((k-i)n-2i+2(i-k(\ell-1)/\ell)\ell,0\right)\\
							&	=	q^{k(\ell-1)}G\left(\frac{k}{\ell}(n-2(\ell-1)),\ell\right),
			\end{align*}
			which to be nonzero requires that $n\ge2(\ell-1)$.\\
			
			The final functions found in Lemmas \ref{lem:ExplicitSize2} and \ref{lem:ExplicitSize3} are by no means obviously unimodal. If one instead writes them in decomposed form, they can be reasoned to be unimodal by looking at the darga of each term. Knowing the proper decomposition beforehand allows for easy proof of unimodality, but finding any decomposition is typically extremely difficult. This is the motivation behind reverse engineering.

			\label{subsec:RestrictedPartitions}
		\end{subsection}%

		\begin{subsection}{Adjusted Initial Conditions}
			We can also multiply each recursive call by constant factors or multiply initial conditions while still maintaining base case darga $=nk$.
			\begin{equation}
				G(n<0,k)	=	0,	\hspace{1em}	G(n,k<0)	=	0,	\hspace{1em}	G(0,k)	=	\nu,	\hspace{1em}	G(n,0)	=	\nu,	\hspace{1em}	G(n,1)	=	\nu\frac{1-q^{n+1}}{1-q}.
				\label{eqn:AdjustedInitialConditions}
			\end{equation}
			We can also adjust the recursive call on each $\lambda\vdash k$ by multiplying it by $\rho$. $\rho$ need not be constant; it can be a function of $|\lambda|$, the largest part of $\lambda$, or even more generally, each $\lambda$ gets its own specific value. 
			We choose to \tbl normalize\tbr\ the final answer by dividing by $\nu^k$ (and $\rho$ when it is constant) so that we do not have arbitrary differences. Even after normalizing, there will still be differences from $G(n,k)$.\\
			
			Though we have not yet found meaningful results using the techniques of this section, we are hopeful that in the future we can find a pattern for the difference between the largest and second-largest coefficients.
			\label{subsec:AdjustedInitialConditions}
		\end{subsection}%

		\begin{subsection}{Adjusted Recursive Call}
			A more complex method is to adjust the recursive call. We return to using the initial conditions in Eqn.\ \eqref{eqn:InitialConditions}.%
			\begin{align*}
				G(n,k;a,b,c)	&=	\sum_{(d_1,\ldots,d_k);\sum_{i=1}^k id_i=k} q^{k\left(\sum_{i=1}^kd_i\right)-k-\sum_{1\le j<i\le k}(i-j)d_id_j}\\
						&\hspace{3.5cm}	\cdot q^{k a+\left(\sum_{i=1}^kd_i\right)b+c[k(k+1)(n/2-a)-k(k+2b)+\sum_{i=1}^kd_i i^2]}\nonumber\\
						&	\prod_{i=0}^{k-1}G\left((k-i)n-2i+2\sum_{j=0}^{i-1}(i-j)d_{k-j}-2a(k-i)-2b,d_{k-i}-2c;a,b,c\right).
			\end{align*}
			The reduction in darga in the recursive call is exactly balanced by the extra factors of $q$.\\
			However, if $c>0$, since there is some $i$ s.t. $d_{k-i}=0$ for any partition, the product will be $0$ by initial conditions, and thus $G(n,k;a,b,c>0)=0$. And if $c<0$, we obtain infinite recursion for nontrivial choices of $n,k$. Therefore, we should ignore this parameter.
			
			\begin{align}
				G(n,k;a,b)	=&	\sum_{(d_1,\ldots,d_k);\sum_{i=1}^k id_i=k} q^{(k+b)\left(\sum_{i=1}^kd_i\right)+k(a-1)-\sum_{1\le j<i\le k}(i-j)d_id_j}\nonumber\\
						&\prod_{i=0}^{k-1}G\left((k-i)(n-2a)-2(i+b)+2\sum_{j=0}^{i-1}(i-j)d_{k-j},d_{k-i};a,b\right).
				\label{eqn:KOHgeneral}
			\end{align}
			{\sc Implemented as {\bf KOHgeneral}}.\\
			
			Unless otherwise stated, we are now assuming the type of partition is unrestricted.
			\begin{proposition}
				\ben
					\item	$G(n,k;a,b)$ will have smallest degree $ka+b$.
					\item	If $a+b>\frac{n}{2}$, then $G(n,k;a,b)=0$.
					\item	If $a+b=\frac{n}{2}$, then $G(n,k;a,b)=q^{ka+b}\frac{1-q^{nk+1-2(ka+b)}}{1-q}$.
				\een
				\bep
					\ben
						\item	This follows from inspecting factors of $q$. Note that $k\left(\sum_{i=1}^kd_i\right)-k-\sum_{1\le j<i\le k}(i-j)d_id_j\ge0$ since $G(n,k)$ has smallest degree $1$. Recall it $=0$ for $\lambda=[k]$.
						\item	Let $\lambda\vdash k$. Consider $n'(k-1)$ in the recursive call:
							\begin{align*}
								n'(k-1)	&=	(k-(k-1))(n-2a)-2((k-1)+b)+2\sum_{j=0}^{(k-1)-1}((k-1)-j)d_{k-j}\\
											&=	(1)(n-2a)-2(k-1+b)+2\sum_{j=0}^{k-1}(k-1-j)d_{k-j}\\
											&=	n-2a-2k+2-2b+2\sum_{j=1}^{k}(k-1-(k-j))d_{j}\\
											&=	n-2a-2b+2+\left(-2k+2\sum_{j=1}^{k}jd_{j}\right)-2\sum_{j=1}^kd_{j}\\
											&=	n-2(a+b)+2-2\sum_{j=1}^kd_{j}\\
											&\le	n-2(a+b)+2-2(1)	<	0.
							\end{align*}
							Thus, $G(n'(k-1),d_{1};a,b)=0$ and since $\lambda$ was arbitrary, $G(n,k;a,b)=0$.
						\item	The above has equality only if $|\lambda|=1$, i.e., $\lambda=[k]$. Then
							\begin{align*}
								G(n,k;a,b)	&=	q^{(k+b)\left(\sum_{i=1}^kd_i\right)+k(a-1)-\sum_{1\le j<i\le k}(i-j)d_id_j}\\
										&\tab	\prod_{i=0}^{k-1}G\left((k-i)(n-2a)-2(i+b)+2\sum_{j=0}^{i-1}(i-j)d_{k-j},d_{k-i};a,b\right)\\
										&=	q^{ka+b}G\left((k)(n-2a)-2(b),1;a,b\right)\\
										&\tab	\prod_{i=1}^{k-1}G\left((k-i)(n-2a)-2(i+b)+2\sum_{j=0}^{i-1}(i-j)d_{k-j},0;a,b\right)\\
										&=	q^{ka+b}G\left(k(n-2a)-2b,1;a,b\right)\\
										&\tab	\prod_{i=1}^{k-1}G\left((k-i)(n-2a)-2(i+b)+2(i-0),0;a,b\right).
							\end{align*}
							It only remains to confirm $(k-i)(n-2a)-2b\ge0$, which follows from $a=\frac{n}{2}-b$.
					\een
				\eep
			\end{proposition}
		
		\end{subsection}%
		
	\end{section}%

	\begin{section}{OEIS}
		When we take $q\to1^-$ for the normal $q$-binomials, we obtain the common binomial coefficients. What happens if we do that to the modified $G_s$ polynomials? There is no direct combinatorial interpretation (yet), but we can derive some cases to start the search. Recall the explicit forms of $G_1(n,k)$, $G_2(n,k)$, and $G_3(n,k)$ in Eqns. \eqref{eqn:ExplicitSize1}, \eqref{eqn:ExplicitSize2}, and \eqref{eqn:ExplicitSize3}, respectively.
		\begin{align*}
			\lim_{q\to1^-}G_1(n,k)	&=	nk+1,\\
			\lim_{q\to1^-}G_2(n,k)	&=	\frac{1}{12}(k+1)(k^2n^2-3k^2n-kn^2+2k^2+9kn-8k+12),\\
			\lim_{q\to1^-}G_3(n,k)	&=	\frac{1}{720}(k+2)(k+1)\bigg(k^3n^3-9k^3n^2-3k^2n^3+26k^3n+42k^2n^2+2kn^3\\
								&\hspace{2.7cm}	-24k^3-153k^2n-33kn^2+162k^2+247kn-378k+360\bigg).
		\end{align*}
		We can also conjecture the form of $G_4(n,n)$ and $G_5(n,n)$ from experimental data:
		\begin{align*}
			\lim_{q\to1^-}G_4(n,n)	&=	\frac{1}{120960}n(n+2)(n+1)\bigg(n^8-32n^7+462n^6-3836n^5+20013n^4\\
								&\hspace{4.75cm}	-66836n^3+140804n^2-171216n+100800\bigg),\\
			\lim_{q\to1^-}G_5(n,n)	&=	\frac{n(n+3)(n+2)(n+1)}{43545600}\bigg(n^{10}-50n^9+1140n^8-15420n^7+136533n^6\\
								&\hspace{4.75cm}	-824370n^5+3436190n^4-9762880n^3\\
								&\hspace{4.75cm}	+18198936n^2-20242080n+10886400\bigg).
		\end{align*}
		The sequences for $G_1(n,n),\ldots,G_5(n,n)$ are now in the OEIS \cite{A002522,A302612,A302644,A302645,A302646}. Only $s=1$ was already in the OEIS.
		\begin{table}[H]
			\caption{Unimodal Sequences from the OEIS}
			\vspace{-1cm}
			\begin{align*}
				\seqnum{A002522}:\hspace{1em}	s=1	\hspace{1em}	&	\tab	1, 2, 5, 10, 17, 26, 37, 50, 65, 82, 101, 122, 145, 170, 197, 226, 257, 290,\\
				\seqnum{A302612}:\hspace{1em}	s=2	\hspace{1em}	&	\tab	1, 2, 6, 20, 65, 186, 462, 1016, 2025, 3730, 6446, 10572, 16601, 25130,\\
				\seqnum{A302644}:\hspace{1em}	s=3	\hspace{1em}	&	\tab	1, 2, 6, 20, 70, 252, 896, 2976, 8955, 24310, 60038, 136500, 289016,\\
				\seqnum{A302645}:\hspace{1em}	s=4	\hspace{1em}	&	\tab	1, 2, 6, 20, 70, 252, 924, 3432, 12705, 45430, 152438, 472836, 1352078,\\
				\seqnum{A302646}:\hspace{1em}	s=5	\hspace{1em}	&	\tab	1, 2, 6, 20, 70, 252, 924, 3432, 12870, 48620, 183755, 683046, 2443168.%
			\end{align*}
		\end{table}[H]
		The sequences will approach the central binomial coefficients \seqnum{A000984} \cite{A000984} since $G_s\to G_\infty=G$.
	\end{section}%

	\begin{section}{Conclusion and Future Work}
		In this work we found a new way to generate unimodal polynomials for which unimodality is, {\it a priori}, very difficult to prove. Many of the recurrences found produce surprisingly unimodal polynomials given the correct initial conditions. The proof is \tbl trivial\tbr\ because of our method of reverse engineering. There are many promising avenues for future work.
		\ben
			\item	%
				We were not able to make interesting use of the techniques in Section 5.2. With much more care, it may be possible to glean some interesting results by adjusting initial conditions.
			\item	An ultimate goal would be to find a restriction on the partitions and other parameters that yields non-trivial polynomials that can be expressed as a function in 2 variables. We have a few for the restricted partition size, but more is always better!
			\item	Find other initial conditions that still yield unimodal polynomials when iterated in the recurrence of Eqn.\ \eqref{eqn:RecurrenceSizeS}.
			\item	Try to fully automate the process of finding $G_s(n,k)$ in Section \ref{subsec:RestrictedPartitions}. Use them to verify Eqn.\ \eqref{eqn:RecurrenceSizeS} for larger values of $s$. The ideal result would be to prove Eqn.\ \eqref{eqn:RecurrenceSizeS} for ALL $s$.
			\item	I also tried looking at how restricting to odd partitions compares to restricting to distinct partitions. As these are conjugate sets, I thought there may be a relation in the resulting polynomials. Alas, I could see no connection. One could reexamine these pairs or look at other partition sets that are known to be equinumerous, e.g., odd and distinct partitions compared to self-congruent partitions.
			\item	An appealing result would be to find combinatorial interpretations for any of the restricted polynomials. Particularly the $G_s$ polynomials and the OEIS sequences \seqnum{A302612} \cite{A302612}, \seqnum{A302644} \cite{A302644}, \seqnum{A302645} \cite{A302645}, and \seqnum{A302646} \cite{A302646}.
			\item	Instead of only using an adjusted Eqn.\ \eqref{eqn:KOHgeneral}, one might be able to \tbl arbitrarily\tbr\ combine polynomials of known darga to create another polynomial with known darga. Do this in a similar recursive manner. Instead of using partitions, are there other combinatorial objects that can be applied? What other unimodal and symmetric polynomials are out there to use as starting blocks? What recurrences do they satisfy that can be tweaked to reverse-engineer unimodal polynomials?
			\item Unimodal probability distributions have special properties. One could use these polynomials for that purpose after normalizing by $G'(n,k;q=1)$. Various properties can be discovered using Gauss's inequality \cite{GaussInequality} or others.
			\item	%
				Unimodal polynomials can be constructed from multiplication very easily. Is the inverse process of factoring them easier than for general polynomials? If so, there could be applications in cryptography and coding theory where factoring is a common theme.
		\een
	\end{section}%

	\label{chap:UnimodalPolynomials}
\end{chapter}%

\begin{chapter}{Lattice Walk Enumeration}

	\begin{section}{Introduction}
		This chapter is available as a stand-alone paper on arXiv.org, number 1803.10920. It is being submitted to the {\it Journal of Integer Sequences}.\\

		Trying to enumerate all of the walks in a 2D lattice is a fun combinatorial problem and there are numerous applications, from polymers to sports. Computers provide a wonderful tool for analyzing these walks; we provide a Maple package for automatically describing generating functions of walks restricted to any step set in a 2D lattice. We always obtain a closed system of relations for generating functions of walks that are bounded, semi-bounded, or unbounded. For bounded walks, this leads to explicit rational solutions! For semi-bounded or unbounded walks, we may get lucky and obtain algebraic solutions; if not, we still have a short self-referential description of the generating function.\\

		We consider walks in the two-dimensional square lattice with an \tbl arbitrary\tbr\ set of integral steps $(x,y)$ subject to $x\ge0$. In addition to unbounded walks, we also separately constrain the walks to lie in regions bounded above and below as well as bounded only below.
		
		One would then like to count all possible walks of a certain length possibly with a specific total change in $y$ value. Rather than a brute-force search of the entire space, looking for 1 value, one could use generating function relations. As a bonus, one would obtain not only the initial generating function of desire, but also many related ones that may be of interest.
		
		Studies of this kind have been done in the literature with simple step sets such as Dyck paths ($\{[1,1],[1,-1]\}$) \cite{DR,BORW} and old-time basketball games \cite{oldtimebasketball}. Philippe Duchon analyzed the case of nonnegative bridges with step set $\{[1,-2],[1,3]\}$; see OEIS\footnote{Online Encyclopedia of Integer Sequences \cite{OEIS}.} sequence \seqnum{A060941} \cite{A060941}. For further developments on the subject, see \cite{R} and the references therein.
		
		The accompanying Maple package is able to extend and inform on old sequences and create many new sequences. Much of the analysis, thus far, has been on steps with $x$-value exactly equal to $1$. One of the aspects of this paper that sets it apart is the ease with which it can analyze more generic cases.

		\begin{subsection}{Motivation}
			We want to enumerate walking paths constrained to specific allowable steps. Most of the time we are looking for paths that begin and end on the $x$-axis.
			
			An earlier motivation for bounded walks came from Physics: analyzing polymers constrained between plates \cite{BORW}. Ayyer and Zeilberger gave one solution in an earlier paper \cite{oldtimebasketball} that provided the main motivation for this research.
			
			The kernel method has received attention lately for analyzing specific cases of walks \cite{basketballKernelMethod}. There are several advantages to describing walks using the method of this paper. The main idea is the same: writing functional equations to describe possible steps in a walk. The difference is that this method then describes \tbl new\tbr\ components of the functional equation in an iterative manner. The kernel method uses analytical number theory on the roots and can be reliant on very case-specific techniques. Compared to the kernel method, we believe our method is a lot easier to understand combinatorially, is more insightful, faster, and easier to produce.\\
			
			Trying to picture an entire walk at once can be difficult. This is where the awesome powers of dynamical programming come into play. Instead of trying to think about the entirety of a walk, think about a part: either the beginning step, end step, or the middle step across the $x$-axis. Break the walk down into different parts (irreducible versus reducible). This is what the generating function equations accomplish.
			
			Originally, we wrote out a single equation to describe the initial walk of interest, then the next, then the next, until it eventually became a closed system. (And the wonderful part is that our descriptions ALWAYS lead to closed systems.) Finally, we solved the system we created. Since this is a very algorithmic approach to answering the question, why not have a computer work for us?
		\end{subsection}%

		\begin{subsection}{Definitions}
			I will generically use the term {\it walk} to indicate any sequence of points $\{(x_0,y_0),\ldots,(x_s,y_s)\}$ in the $xy$-plane. Though the walk will not necessarily start at the origin, if the starting point is not given, it is assumed to begin at $(0,0)$. The steps of a walk are $\{(x_1-x_0,y_1-y_0),\ldots,(x_s-x_{s-1},y_s-y_{s-1})\}$ and are built from some step set $\S$, a set of ordered pairs. I exclusively consider walks that are monotonically weakly increasing: $x_{i+1}-x_i\ge0$. A walk is {\it nonnegative} ({\it nonpositive}) if the walk never crosses below (above) the $x$-axis. I will sometimes refer to the $y$-value as the {\it altitude} of the walk.
			\bed[Walks]
				A {\it bridge} is an unbounded walk that begins at the origin and ends on the $x$-axis. I say {\it bounded bridge} for a bounded walk that begins at the origin and ends on the $x$-axis.\\
				An {\it excursion} is a semi-bounded (not necessarily nonnegative) walk that begins at the origin and ends on the $x$-axis.\\
				A {\it free} walk can end anywhere (not ending at any specific altitude). A {\it meander} is a semi-bounded free walk.\\
				The {\it length} of a walk is $n=\sum_{i=0}^s x_i$. The {\it size} of a walk is $s$.
			\eed
			See Figure \ref{fig:WalkExamples} for examples of the different types of walks.
			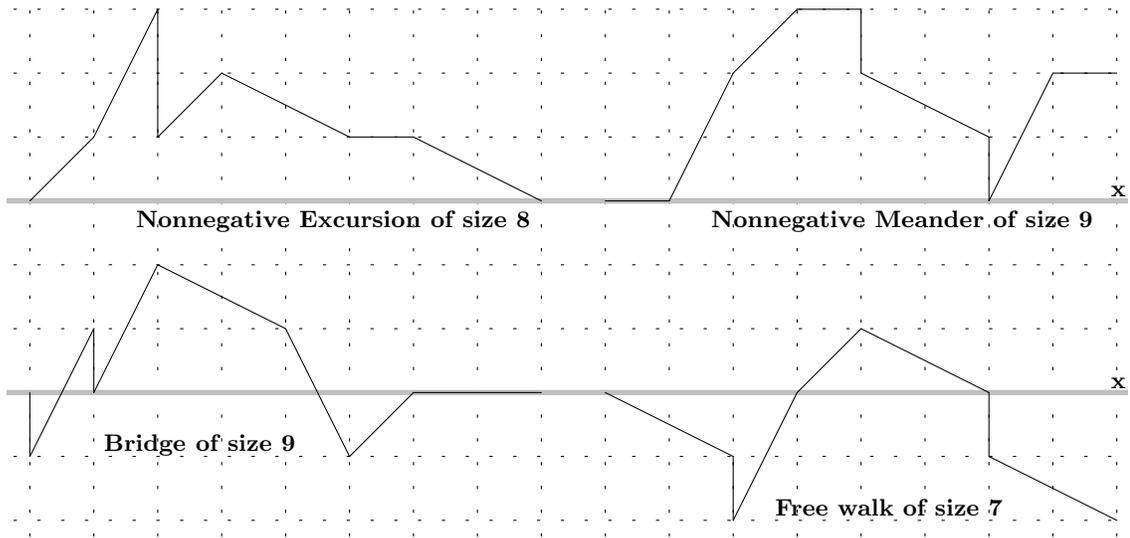
\begin{figure}[H]
				\caption{Walk Examples}
				\begin{tikzpicture}[line cap=round,line join=round,>=triangle 45,x=0.85cm,y=0.85cm]
					\draw [color=black,dash pattern=on 1pt off 8pt, xstep=0.85cm,ystep=0.85cm] (-0.25,-5.25) grid (17.25,3.25);
					\clip(-0.35,-5.25) rectangle (17.17,3.36);
					\draw [line width=2pt,color=lightgray] (-1,0)-- (18,0);
					\draw (0,0)-- (1,1);
					\draw (1,1)-- (2,3);
					\draw (2,3)-- (2,2);
					\draw (2,2)-- (2,1);
					\draw (2,1)-- (3,2);
					\draw (5,1)-- (6,1);
					\draw (6,1)-- (8,0);
					\draw (3,2)-- (5,1);
					\draw [line width=2pt,color=lightgray] (-1,-3)-- (18,-3);
					\draw (0,-3)-- (0,-4);
					\draw (0,-4)-- (1,-2);
					\draw (1,-2)-- (1,-3);
					\draw (1,-3)-- (2,-1);
					\draw (2,-1)-- (4,-2);
					\draw (4,-2)-- (5,-4);
					\draw (5,-4)-- (6,-3);
					\draw (6,-3)-- (7,-3);
					\draw (7,-3)-- (8,-3);
					\draw (9,0)-- (10,0);
					\draw (10,0)-- (11,2);
					\draw (11,2)-- (12,3);
					\draw (12,3)-- (13,3);
					\draw (13,3)-- (13,2);
					\draw (13,2)-- (15,1);
					\draw (15,1)-- (15,0);
					\draw (15,0)-- (16,2);
					\draw (16,2)-- (17,2);
					\draw (9,-3)-- (11,-4);
					\draw (11,-4)-- (11,-5);
					\draw (11,-5)-- (12,-3);
					\draw (12,-3)-- (13,-2);
					\draw (13,-2)-- (15,-3);
					\draw (15,-3)-- (15,-4);
					\draw (15,-4)-- (17,-5);
					\draw (16.75,0.4) node[anchor=north west] {\footnotesize \bf x};
					\draw (16.75,-2.6) node[anchor=north west] {\footnotesize \bf x};
					\draw (1.5,0) node[anchor=north west] {\footnotesize \bf Nonnegative Excursion of size 8};
					\draw (10.5,0) node[anchor=north west] {\footnotesize \bf Nonnegative Meander of size 9};
					\draw (1,-3.5) node[anchor=north west] {\footnotesize \bf Bridge of size 9};
					\draw (11.5,-4.5) node[anchor=north west] {\footnotesize \bf Free walk of size 7};
				\end{tikzpicture}
				All of the example walks are of length 8 and are considered to have begun at the origin. They are all built from the same step set $\S=\{[0,-1],[1,0],[1,1],[1,2],[2,-1]\}$. Note the vertical line in the excursion is actually 2 steps. The excursion happens to be irreducible while the rest are not.
				\label{fig:WalkExamples}
			\end{figure}
			Banderier, {\it et al.} \cite{basketballKernelMethod} uses {\it walk/path} to reference any sequence of steps that start at the origin. The only difference in definitions is that they always count excursions and meanders as nonnegative. I will mostly consider them as nonnegative and explicitly say when I am not.
			\bed
				The {\it interior} of a walk consists of every point other than the endpoints: $\{(x_1,y_1),\ldots,(x_{s-1},y_{s-1})\}$.\\
				An {\it irreducible} walk is one in which the interior has a strictly higher altitude than the lower endpoint: $\min\{y_1,\ldots,y_{s-1}\}>\min\{y_0,y_s\}$. For purposes of this paper and Maple package, the stationary walk\footnote{A single point.} and walks that are direct steps to the right are NOT considered to be irreducible walks.\footnote{As the interior (the edge between points) is not strictly higher than the endpoints.}
			\eed
			Irreducible is also used to refer to walks that do not exactly hit the final altitude until the final step.
		\end{subsection}%

		\begin{subsection}{Chapter Organization}
			I will exclusively use $t$ as the variable in generating functions. I also abbreviate generating function(s) as g.f.(s). This chapter is organized in the following sections:
			\ben[label=\bfseries Section \arabic*:,leftmargin=21mm]
				\setcounter{enumi}{1}
				\item Bounded: This is the section in which a computer does the best; it can give an exact g.f. solution since we have a \tbl linear\tbr\ system of equations. Computers are very good at solving these quickly and efficiently.
				\item	Semi-bounded: An algebraic expression for the g.f. is no longer guaranteed. However, we can always find a polynomial in $\Z[t]$ for which the g.f. is a root.\footnote{In the formal power series sense.}
				\item	Guess-and-check: Finding a minimal polynomial is guaranteed so why not just guess? This section gives a time and memory comparison showing why that would be a poor decision.
				\item	Algebraic-to-Recursive: Using the minimal polynomial of the g.f. is not always the fastest for enumerating. This section introduces another Maple package that converts the polynomial into a 1D recurrence.
				\item	Unbounded: We find minimal polynomials for g.f.s of unbounded walks and show an alternative method of producing recurrences in specific cases.
				\item	Asymptotics: We discover asymptotic results for several step sets to show relationships between the number of excursions, bridges, and, to a lesser degree, meanders.
				\item	Applications: We use the Maple package to find some extended results and talk about probabilistic behavior. We also express minimal polynomials for excursions of small step sets.
				\item	Conclusions and Future Work: A brief recap of what we can do with the {\bf ScoringPaths} Maple package and how we can extend the work.
			\een
			This paper is produced in conjunction with the Maple package {\bf ScoringPaths}. It is downloadable from
			\begin{center}{\bf math.rutgers.edu/$\sim$bte14/Code/ScoringPaths/ScoringPaths.txt}.\end{center}
			I will mention various functions of the package in {\bf bold}. All functions mentioned in this paper are included in the accompanying Maple package. Some functions have been borrowed from other packages with credit and included in {\bf ScoringPaths} for completeness. All comparisons of time and memory are done with Maple 17's {\bf CodeTools[Usage]} on Linux version 3.10.0-514.el7.x86\_64 with 8GB of RAM. All values are averages of at least 30 trials unless otherwise specified. Pay careful attention to the units in some examples.
			
			We are considering discrete walks. As such, we can consider steps with only integral values. If the $x$-steps are fractional, there will likely not be any issues. If the $y$-steps are fractional, many functions will not work as intended. If all of the steps share a common factor in the $y$-value, it should be factored out leaving an equivalent problem. The bounds will need to be factored and truncated as well. Leaving the common factor may cause problems with some functions.
			
			The g.f. may produce non-zero coefficients only every $m^{th}$ value. It may be desirable to make the substitution $t\to t^{1/m}$ for ease of reading. I commonly use $B(n)$, and occasionally $C(n)$ as the number of walks of length $n$. The step set and walk restrictions will be obvious from context or given explicitly.
			
			A few functions to create sets of random steps have been included for quick demonstrations; {\bf RandomStepSet} produces a generic set, {\bf RandomZeroStepSet} parses for walks that begin and end on the $x$-axis, {\bf RandomSemiBoundedStepSet} parses for walks bounded below, and {\bf RandomUnBoundedStepSet} produces a step set without any $[0,y]$ steps.
			
			There are several \tbl paper\tbr\ functions that automatically produce an article with information about a given step set. See the \tbl Paper-Producing Functions\tbr\ section of the {\bf Help} function.\\
		
			The package {\bf gfun} by Salvy and Zimmermann \cite{salvy1994gfun}, a staple included with recent Maple versions, is also very useful for manipulating g.f.s. It contains many functions for translating between algebraic expressions, recurrences, and differential equations satisfied by g.f.s.
		\end{subsection}%
		
	\end{section}%

	\begin{section}{Bounded}
		The first case is walks that are bounded above and below. Consider an arbitrary set of steps $\S=\{(x_1,y_1),\ldots,(x_n,y_n)\}$. The goal is to find the g.f., denoted $f_{a,b}$, for walks with step set $\S$, starting at $(0,0)$, and bounded above and below by $a\ge0$ and $b\le0$, respectively. All walks/paths in this section are constrained to a step set $\S$ and bounded above and below.

		\begin{subsection}{Walking \tbl Anywhere\tbr}
			First assume that the walk can end anywhere between the lines $y=a$ and $y=b$. A walk either never goes anywhere, $+1$, or it takes a step ($t^x$) and continues as if it is a walk starting at a new point: $f_{a-y,b-y}$. The following relation accomplishes that.
			\beq
				f_{a,b}	=	1+\sum_{(x,y)\in\S}t^{x}f_{a-y,b-y}.
			\eeq
			This by itself gets us nowhere. But if $a-y<0$ or $b-y>0$, then $f_{a-y,b-y}=0$ since it is already starting in a prohibited region. Now write out all $(a-b)$ equations:
			\begin{align*}
				f_{0,b-a}	&=	1+\sum_{(x,y)\in\S}t^{x}f_{0-y,b-a-y},\\
				f_{1,b-a+1}	&=	1+\sum_{(x,y)\in\S}t^{x}f_{1-y,b-a+1-y},\\
						&\hspace{2mm}\vdots\\
				f_{a-b,0}	&=	1+\sum_{(x,y)\in\S}t^{x}f_{a-b-y,0-y}.
			\end{align*}
			This is a (linear!) system of $a-b$ equations with $a-b$ variables, after discarding the constant $0$ g.f.s, for which Maple's solve function can easily find the solution. As a bonus, we not only have $f_{a,b}$, but also every $f_{m,n}$ such that $m-n=a-b$ (and $m\ge0$, $n\le0$).
			
			The function {\bf BoundedScoringPathsEqsVars} will produce all of the equations and variables for this system extremely quickly. {\bf BoundedScoringPaths} will output only the g.f. $f_{a,b}$. Again, this is very fast since it is solving a linear system of only $(a-b)$ equations. To verify the values, one can use {\bf BoundedScoringPathsNumber}, which computes the number of walks by recursion.
			
			Allowable steps are any as long as there is not a $(0,-)$ AND a $(0,+)$ such that they can sum to $0$ within the bounds. Specifically, if $(0,m)$ and $(0,-n)$ are steps, then it is allowable (will not produce infinite values) as long as the width is small enough: $a-b<n+m-gcd(m,n)$.
			\bep
				Let $\S=\{[0,m],[0,-n]\}$ and assume $m>n$. Also assume $gcd(m,n)=1$. If $a-b\ge n+m-1$, then we can always take at least one of the $\S$ steps. Since there is a finite number of altitudes, then must be acollision at some point: a loop.
				
				If $gcd(m,n)=g>1$, then reconsider the same problem after factoring $g$ out of everything: $a'=\floor{\frac{a}{g}}$, $b'=\ceil{\frac{b}{g}}$, $\S'=\{[0,m/g],[0,-n/g]\}$.
			\eep
			This bound only works for two 0-steps. For more than two 0-steps, the width must be even shorter. For some combinations (e.g., $\{[0,u-v],[0,-u],[0,u+v]\}$), there is no allowable width that includes all steps.\footnote{Assuming that one can get to the edges of the boundary. $(u+v)+(-u)+(u-v)+(-u)=0$: a loop.}

			\begin{subsubsection}{Examples}
				\begin{example}[Close (American) Football Games]
					Consider trying to enumerate the number of American football games in which the teams are never separated by more than 1 score. On a given play,\footnote{Counting untimed downs as part of the previous play.} it is possible for one team to score $2,3,6,7$, or $8$ points. It is also possible for one team to score $6$ points and the opposing team to score 1 or 2 points, though these are pretty rare occurrences. What if we want to enumerate the number of games with $n$ scoring plays that end separated by no more than 1 score? We can use the step set
					\begin{align*}
						\S=\{&[1,2],[1,3],[1,6],[1,7],[1,8],[1,5],[1,4],\\
							&[1,-2],[1,-3],[1,-6],[1,-7],[1,-8],[1,-5],[1,-4]\}.
					\end{align*}
					Each step represents 1 scoring play and how many points the home team gained relative to the away team. Since the maximum scoring play is 8 points, we make the bounds $y=-8,8$. The g.f. is then found with {\bf BoundedScoringPaths($\S$,8,-8,t)}:
					\beq
						\frac{1+10t+13t^2-37t^3-40t^4+28t^5+26t^6-2t^7}{1-4t-59t^2-77t^3+170t^4+234t^5-92t^6-142t^7-4t^8+6t^9}.
					\eeq
					This sequence is new in the OEIS: \seqnum{A301379} \cite{A301379}.
					\label{exa:CloseAmericanFootballGames}
				\end{example}
				
				\begin{example}[Speed Enumeration - Bounded 1]
					Consider finding the first 1000 coefficients of the g.f. found in Example \ref{exa:CloseAmericanFootballGames} above. We could use brute-force recursion in {\bf BoundedScoringPathsNumber} or take a taylor series expansion since we have an explicit form for the g.f. from {\bf BoundedScoringPaths}.
					\begin{table}[H]\begin{center}
						\caption{Bounded Walk Enumeration}
						\begin{tabular}{|c|c|c|c|c|}
							\hline
							Method	&	Memory Used	&	Memory Allocation	&	CPU Time	&	Real Time\\
							\hline
							Brute-Force Recursion	&	60.23MiB	&	24MiB	&	309.53ms	&	308.90ms\\
							\hline
							G.F. Construction	&	3.27MiB	&	0 bytes	&	30.17ms	&	33.50ms\\
							Taylor Enumeration	&	4.86MiB	&	0 bytes	&	4.93ms	&	5.07ms\\
							\hdashline
							Total	&	8.13MiB	&	0 bytes	&	35.10ms	&	38.57ms\\
							\hline
						\end{tabular}
						\label{tab:BoundedEnumeration}
					\end{center}\end{table}
					Using the g.f. method of enumeration is about 8 times as fast and uses a much smaller amount of memory: ($1/8$th).
				\end{example}
			\end{subsubsection}%
		\end{subsection}%

		\begin{subsection}{Returning to the $x$-axis}
			\label{subsec:EqualBounded}
			The above equations are for g.f.s of walks that are able to end ANYWHERE. What if we want a g.f. for walks that must end back at the $x$-axis? Ayyer and Zeilberger solved the bounded bridge case (but not the bounded free walk case) using different relations but the same general method \cite{BoundedCase}; they also use the equivalent problem of walks between parallel lines of positive slope instead of between horizontal lines.
			
			Let $f_{a,b}$ now denote the g.f. for walks that begin at $(0,0)$ and end on the $x$-axis. Again, we either never get moving, $+1$, or we take a step and then must take a path back to the $x$-axis. Let $e_{a,b,c}$ denote the g.f. for paths that start at $(0,c)$, end on the $x$-axis and never touch the $x$-axis beforehand. Now that $e_{a,b,c}$ is introduced, we can write the relation for $f_{a,b}$:
			\begin{equation}
				f_{a,b}	=	1+\left(\sum_{(x,0)\in\S}t^x\right)f_{a,b}+\left(\sum_{(x,y)\in\S;y\ne0}t^x e_{a,b,y}\right)f_{a,b}.
				\label{eqn:EqualBounded}
			\end{equation}
			There is also the possibility of just moving along the $x$-axis to start. After returning to the $x$-axis, we can now take any walk as before. Hence the multiplication by $f_{a,b}$. We now need the equations for $e_{a,b,c}$ for $a\ge c\ge b$. If any step would return to the $x$-axis, then we are done. Otherwise, it must continue as a new $e_{a,b,c'}$ path.
			\beq
				e_{a,b,c}	=	\sum_{(x,-c)\in\S}t^x	+	\sum_{(x,y)\in\S;y\ne-c}t^x e_{a,b,c+y}.
			\eeq
			We now have a set of linear equations to solve for the $e_{a,b,c}$. To produce all of the equations and variables, use {\bf EqualBoundedScoringPathsEqsVars}. Once the $e_{a,b,c}$ system is solved, Eqn.\ \eqref{eqn:EqualBounded} is linear in $f_{a,b}$ so we (or rather a computer) will find a rational g.f. solution! Use {\bf EqualBoundedScoringPaths} for the single solution. To verify the result, one can check using the enumeration in {\bf EqualBoundedScoringPathsNumber}. The irreducible g.f.s can be checked with {\bf SpecificIrreducibleBoundedScoringPathsNumber}.

			\begin{subsubsection}{Old-Time Basketball}
				\label{subsubsec:Basketball}
				The methods of relating generating functions were inspired by Ayyer and Zeilberger's work \cite{oldtimebasketball}. They found the following relation.
				\bet
					Let $F_w$ denote the generating function for the number of walks subject to step set $\S=\{[\frac{1}{2},1],[\frac{1}{2},-1],[1,2],[1,-2]\}$, that start at $(0,0)$, end at $(n,0)$, and never go below the $x$-axis or above the line $y=w$. Then $F_w$ satisfies the following recurrence relation:
					\begin{align*}
						F_w	&=	1	-	t F_w	+	2t F_w F_{w-1}	+	2 t^2 F_w F_{w-1} F_{w-2}\\
							&\tab	- (t^3 + t^4)F_w F_{w-1} F_{w-2} F_{w-3}	+	t^5 F_w F_{w-1} F_{w-2} F_{w-3} F_{w-4}.
					\end{align*}
					\label{thm:OldTimeBasketball}
				\eet
				\begin{proposition}
					The initial conditions are given as follows:
					\begin{align*}
						F_0	=	1,	\tab\tab	
						F_1	=&	\frac{1}{1-t},	\tab\tab
						F_2	=	\frac{1-t}{1-2t+3t^2},\\
						F_3	=	\frac{1-2t+3t^2}{1-3t-5t^2-2t^3+t^4},	&\tab\tab
						F_4	=	\frac{1-3t-5t^2-2t^3+t^4}{1-4t-6t^2+2t^3}.
					\end{align*}
				\end{proposition}
				If we compute {\bf EqualBoundedScoringPaths(\{[1/2,1],[1/2,-1],[1,2],[1,-2]\},w,0,t)} for $w=0,\ldots,4$, then the initial conditions match. And we can verify Theorem \ref{thm:OldTimeBasketball} empirically for any fixed $w$.
			\end{subsubsection}%

			\begin{subsubsection}{Examples}
				\begin{example}[Tied (American) Football Games]
					Consider the similar Example \ref{exa:CloseAmericanFootballGames}, but this time we want the teams to be tied at the end of the game. Our step set is the same. The only difference is in which method we use: {\bf EqualBoundedScoringPaths($\S$,8,-8,t)}:
					\beq
						\frac{1-4t-45t^2-43t^3+98t^4+108t^5-24t^6-30t^7}{1-4t-59t^2-77t^3+170t^4+234t^5-92t^6-142t^7-4t^8+6t^9}.
					\eeq
					This sequence is new in the OEIS: \seqnum{A301380} \cite{A301380}.
					\label{exa:TiedAmericanFootballGames}
				\end{example}
				
				\begin{example}[Speed Enumeration - Bounded 2]
					Once again let us find the first 1000 coefficients of a g.f.: this time the one found in Example \ref{exa:TiedAmericanFootballGames} above.
					\begin{table}[H]\begin{center}
						\caption{Bounded Bridge Enumeration 1000 terms}
						\begin{tabular}{|c|c|c|c|c|}
							\hline
							Method	&	Memory Used	&	Memory Allocation	&	CPU Time	&	Real Time\\
							\hline
							Brute-Force Recursion	&	55.98MiB	&	24MiB	&	278.93ms	&	279.03ms\\
							\hline
							G.F. Construction	&	10.13MiB	&	8MiB	&	102.43ms	&	103.20ms\\
							Taylor Enumeration	&	4.85MiB	&	0 bytes	&	5.20ms	&	5.10ms\\
							\hdashline
							Total	&	14.98MiB	&	8MiB	&	107.63ms	&	108.30ms\\
							\hline
						\end{tabular}
						\label{tab:EqualBoundedEnumeration}
					\end{center}\end{table}
					Again, the g.f. method is faster: this time, about 2.5 times as fast and $1/4$th the memory.
				\end{example}
			\end{subsubsection}%
	
		\end{subsection}%
	
		\label{sec:Bounded}
	\end{section}%

	\begin{section}{Semi-bounded}
		We now remove the restriction of bounding walks from above.\footnote{To look at walks solely bounded above, simply change the sign of the $y$-value of every step. Or note that every nonnegative excursion is in bijection with a nonpositive excursion by reversing the order of steps.} Ayyer and Zeilberger also provided relations for describing excursions that are similar to those included here \cite{BoundedCase}. Duchon had previously tackled the case of excursions using much different language, but the same overall method \cite{Duchon}.
		
		Let $\S$ denote the set of steps as in the previous section. Let $f$ now denote the g.f. for nonnegative excursions with step set $\S$ that begin at $(0,0)$ and end on the $x$-axis. How will we describe $f$ in an equation? We could try to use the same method as Section \ref{subsec:EqualBounded} and write
		\beq
			f	=	1+\left(\sum_{(x,0)\in\S}t^x\right)f+\left(\sum_{(x,y)\in\S;y\ne0}t^x e_{y}\right)f.
		\eeq
		where $e_{y}$ is the g.f. for walks that start at $(0,y)$, end at $(n,0)$ and never hit the $x$-axis beforehand. However, to describe $e_y$ as in Section \ref{subsec:EqualBounded} (by looking at the first step) would require writing equations for all $e_{i>0}$ since a walk could get arbitrarily far from the $x$-axis before returning.\footnote{We could instead look at the final step, but that is essentially the method described in the next paragraph.} At some point, all of the $e_i$ equations would look essentially the same and it may be possible to take limits of their form to find a solution for $e_y$. However, there is an easier method.\\
		
		Ayyer and Zeilberger \cite{oldtimebasketball} used a standard idea in combinatorics of \tbl irreducible\tbr\ walks to describe $F_w$ in Section \ref{subsubsec:Basketball}. We will use them as well. Let $f_{a,b}$ denote the g.f. for nonnegative walks that start at $(0,a)$ and end at $(n,b)$. Note that $f_{0,0}=f$ is what we are typically looking for. Let $g_{a,b}$ denote the g.f. for walks that start at $(0,a)$, end at $(n,b)$, and stay above the line $y=\min(a,b)$ except at the respective endpoint. Note then that $g_{a,b}=g_{a-b,0}$ (or $=g_{0,b-a}$ if $b>a$).\footnote{In this definition, $g_{0,0}$ does not include stationary walks or walks that are solely a step directly to the right. If this is what one desires, then use $g=g_{0,0}+1+\left(\sum_{(x,0)\in\S}t^x\right)$.} Then
		\begin{equation}
			f_{0,0}	=	1+\left(g_{0,0}+\sum_{(x,0)\in\S}t^x\right)f_{0,0}.
		\end{equation}
		Either the walk is stationary (+1) or it returns to the $x$-axis ($g_{0,0}+\sum_{(x,0)\in\S}t^x$) and then continues as if it were new (multiplication by $f_{0,0}$).
		\begin{equation}
			g_{0,0}	=	\sum_{(x_1,y_1)\in\S;y_1>0}\sum_{(x_2,y_2)\in\S;y_2<0}t^{x_1}f_{y_1-1,-y_2-1}t^{x_2}.
		\end{equation}
		An irreducible $(0,0)$-walk can be characterized by how it departs from the $x$-axis ($t^{x_1}$) and how it returns ($t^{x_2}$). There must be an intermediate walk between these two steps ($f_{y_1-1,-y_2-1}$) that does not touch the $x$-axis: hence the shift. Now we need to describe each new $f_{a,b}$.
		\begin{align}
			a>b	&\tab	f_{a,b}	=	\sum_{i=0}^b g_{a-i,0} f_{0,b-i},\\
			\label{eqn:SemiBoundedFaa}
			a=b	&\tab	f_{a,a}	=	\sum_{i=0}^{b-1} g_{a-i,0} f_{0,b-i}+f_{0,0},\\
			a<b	&\tab	f_{a,b}	=	\sum_{i=0}^{a-1} g_{a-i,0} f_{0,b-i}+f_{0,0}g_{0,b-a}.
		\end{align}
		We can characterize a walk by how close it gets to the $x$-axis. An $f_{a,b}$ walk can go $i$ levels below the level of $y=\min(a,b)$. In this case, we have an irreducible walk down to the lowest point ($g_{a-i,0}$), and then an arbitrary walk to the final level that does not go any lower ($f_{0,b-i}$). If $a=b$, then a walk that does not go below level $y=a$ is equivalent to $f_{0,0}$. If $a<b$, then a walk that does not go below level $y=a$ consists of an arbitrary walk back to the same level without going lower ($f_{0,0}$), followed by an irreducible walk to the final level of $y=b$ ($g_{0,b-a}$). Now we must describe the irreducible walks that have not been covered.
		\begin{align}
			\label{eqn:IrreducibleDown}
			g_{a,0}	&=	\sum_{(x,y)\in\S;y<0}f_{a-1,-y-1}t^x,\\
			\label{eqn:IrreducibleUp}
			g_{0,a}	&=	\sum_{(x,y)\in\S;y>0}t^x f_{y-1,a-1}.
		\end{align}
		An irreducible walk can be characterized by how it reaches ($t^x$) the lowest point. The rest of the walk is arbitrary as long as it does not hit the $x$-axis: again, hence the shift.
		
		We need to show that this iterated process does eventually terminate. This can easily be justified because the largest index of either $g$ or $f$ will be $<\max_{(x,y)\in\S}|y|$: the maximum step size. Thus, there is a finite number of $g$, $f$ that we need to describe in order to have a closed system. Use {\bf EqualSemiBoundedScoringPathsEqsVars} to generate these variables and their equations.\\
		
		The equations are no longer linear and so Maple's built-in solve function does not work as nicely. Instead, we use the {\bf Basis} procedure in the {\bf Groebner} package to find
		\bed[Minimal Polynomial]
			the (minimal) algebraic equation satisfied by the generating function: $p$ such that $p(f)=0$ in terms of formal power series. We refer to $p$ as the {\it minimal polynomial}. When discussing the {\it degree} of the minimal polynomial, we are considering the degree of $f$ in $p$ unless explicitly stated that we are considering $t$.
		\eed
		To produce this polynomial, use {\bf EqualSemiBoundedScoringPaths}. To produce an \tbl ideal\tbr\ polynomial for each $g,f$, use {\bf AllEqualSemiBoundedScoringPaths}.\footnote{This can be slow as it requires finding a Groebner basis for {\it every} variable.} If you want the minimal polynomial for only a single variable, use {\bf SpecificEqualSemiBoundedScoringPaths}. The minimal polynomial is typically found when the Groebner basis receives the desired variable as the lowest order lexicographically. Occasionally the Groebner basis will not produce the minimal polynomial, but instead a product of it and another polynomial. One can recover the minimal polynomial by factoring the output and testing which factor satisfies $p(f)=0$ in terms of formal power series.\footnote{By checking a truncated version of $f$. In most cases $f=1$ is sufficient to see which factor works. In general, enumerate $m$ terms first using the proper enumerating function and see if the corresponding polynomial is $O(t^m)$, signifying a root.} Every function that outputs an algebraic expression has had {\bf FindProperRoot} appended to the end to properly parse the minimal polynomial.\\
		
		\begin{example}[Semi-bounded Example]
			Let $\S=\{[0,-3],[1,-2],[2,0],[3,1]\}$. The g.f., $F$, for the number of excursions that do not go below $y=-1$ and have step set $\S$ is found to satisfy, using {\bf EqualSemiBoundedScoringPaths($\S$,-1,t,F)},
			\beq
				{t}^{18}{F}^{4}+{t}^{14} \left( t^2-1 \right) {F}^{3}+2{t}^{7} \left( t^2-1 \right) ^{2}{F}^{2}+ \left( t^2-1 \right) ^{5}F+ \left( t^2-1 \right) ^{4}	=	0.
			\eeq
			This is actually a fairly simple answer. If we change $[0,-3]\to[0,3]$, then the minimal polynomial is degree $10$ in $F$ and takes 4 lines to write.
			\label{exa:EqualSemiBounded}
		\end{example}
		This polynomial can then be used to discover the enumeration hidden in the coefficients of the taylor series expansion by setting $f_0=1$ and then iterating $f_0\to p(f_0)+f_0$ to find a fixed point solution.\footnote{This glosses over why the convergence works. It becomes an issue in the unbounded case.} Finding the coefficients in this manner requires finding $p$: use {\bf EqualSemiBoundedScoringPathsCoefficients}.
		
		Finding the minimal polynomial can be a time-consuming process. An alternative method is to iterate a fixed-point solution of a vector of all $g,f$ and then pick out $f_{0,0}$. Initialize $\{f_{a,a}=1\}$ and every other g.f. to 0.\footnote{Staying stationary is only valid for general walks that begin and end at the same level.} Now iterate all of the $g,f$ into their respective equations, truncating to the desired coefficient. Eventually we will reach a fixed-point for the vector of solutions. To do this, use {\bf EqualSemiBoundedScoringPathsSeries}. However, for most reasonable calculations of the coefficients ($<1000$ terms), finding the values via brute-force recursion (with {\bf SpecificEqualBoundedScoringPathsNumber}) is faster than either iterating technique.%
		\begin{example}[Speed Enumeration - Excursions]
			We want to obtain the number of nonnegative excursions with step set $\S=\{[1,-2],[1,-1],[1,0],[1,1],[1,2]\}$. Let $F$ denote the corresponding g.f.; then\footnote{Using {\bf EqualSemiBoundedScoringPaths($\S$,0,t,F)}.}
			\begin{equation}
				t^4 F^4	-	t^2(t+1)F^3	+	t(t+2) F^2	-	(t+1) F	+	1	=	0.
				\label{eqn:SpeedSemiBounded}
			\end{equation}
			A truncated solution in formal power series, and the one that makes sense in terms of our problem, is
			\beq
				F	=	1+t+3t^2+9t^3+32t^4+120t^5+473t^6+1925t^7+8034t^8+34188t^9+147787t^{10}.%
			\eeq
			We compare and contrast the amount of time and memory to enumerate the first 500 and 1000 coefficients of $F$ in the following tables.
			\ben
				\item	Set up the g.f. equations and then iterate a vector of solutions using {\bf EqualSemiBoundedScoringPathsSeries}.
				\item	Solve for the minimal polynomial (Eqn.\ \eqref{eqn:SpeedSemiBounded}) and iterate a single solution with {\bf EqualSemiBoundedScoringPaths}.
				\item	Use Maple's {\bf taylor} function on the minimal polynomial.\footnote{Requires replacing $F$ by $\_Z$ and using {\bf RootOf(p)}.}
				\item	Use brute-force recursion and Maple's option remember: {\bf SpecificEqualSemiBoundedScoringPathsNumber}.
			\een
			\begin{table}[H]\begin{center}
				\caption{500 term Excursion Enumeration}
				\begin{tabular}{|c|cccc|}
					\hline
					Method	&	Memory Used	&	Memory Allocation	&	CPU Time	&	Real Time\\
					\hline
					Vector Set-Up	&	16.58KiB	&	0 bytes	&	300$\mu$s	&	766$\mu$s\\
					Iterating	&	100.04GiB	&	55.64MiB	&	4.45m	&	4.09m\\
					\hdashline
					Total	&	100.04GiB	&	55.64MiB	&	4.45m	&	4.09m\\
					\hline
					Polynomial	&	2.63MiB	&	0 bytes	&	19.23ms	&	18.63ms\\
					Iterating	&	23.61GiB	&	20.84MiB	&	64.52s	&	58.56s\\
					\hdashline
					Total	&	23.61GiB	&	20.84MiB	&	64.54s	&	58.58s\\
					\hline
					Polynomial	&	2.63MiB	&	0 bytes	&	19.23ms	&	18.63ms\\
					taylor	&	328.97MiB	&	53.85MiB	&	3.20s	&	3.18s\\
					\hdashline
					Total	&	331.60MiB	&	53.85MiB	&	3.22s	&	3.20s\\
					\hline
					Brute-Force Recursion	&	186.60MiB	&	328MiB	&	1.569s	&	1.472s\\
					\hline
				\end{tabular}
				\label{tab:EqualSemiBoundedEnumeration500}
			\end{center}\end{table}
			Finding the polynomials takes negligible time and memory in comparison to actually enumerating the coefficients beyond the first few terms.
			\begin{table}[H]\begin{center}
				\caption{1000 term Excursion Enumeration}
				\begin{tabular}{|c|cccc|}
					\hline
					Method	&	Memory Used	&	Memory Allocation	&	CPU Time	&	Real Time\\
					\hline
					Vector Iterating	&	1.31TiB	&	159.26MiB	&	63.84m	&	51.86m\\
					\hline
					Single Iterating	&	300.84GiB	&	15.32MiB	&	14.50m	&	11.93m\\
					\hline
					taylor	&	1.82GiB	&	489.68MiB	&	12.20s	&	11.12s\\
					\hline
					Brute-Force Recursion	&	0.85GiB	&	508MiB	&	7.53s	&	7.13s\\
					\hline
				\end{tabular}
				\label{tab:EqualSemiBoundedEnumeration1000}
			\end{center}\end{table}
			So using recursion is the fastest, but also must reserve the most memory. The other methods use more memory in total but can recycle much of what they used. We can best all of the methods with another recursion that is faster and uses less memory; see Example \ref{exa:EqualSemiBoundedConverted}. 
			\label{exa:EqualSemiBoundedSpeed}
		\end{example}
		Before implementing {\bf FindProperRoot}, the Groebner basis output a degree 5 polynomial that took over 5 times as long to iterate as the degree 4 minimal polynomial. Whether iterating a single polynomial or the entire vector of solutions is faster depends on the degree of the minimal polynomial and the number of variables in the closed system. And the number of terms to enumerate. An interesting question would be to look at the time-complexity of this method of enumeration. The vector and single polynomial iteration could potentially be optimized to pick out coefficients from each monomial instead of expanding everything and then picking only relatively few terms. This would help the minimal polynomial iteration much more as the vector iteration only relies on degree 2 expressions.

		\begin{subsection}{Arbitrary Lower Limit}
			What if we want to consider walks that stay above an arbitrary lower bound $y=-c$ as in Example \ref{exa:EqualSemiBounded}? This is actually very easy. By shifting the walk to be nonnegative, we are now looking for $f_{c,c}$. Describe $f_{c,c}$ using Eqn.\ \eqref{eqn:SemiBoundedFaa}. Again, iterating on all new $g,f$ will eventually yield a closed system since the indices are bounded by $\max(\max_{(x,y)\in\S}|y|-1,c)$.
			
			This shifting technique is what the Maple package utilizes when it is given a lower limit other than 0 for semi-bounded walks. The Maple functions are always focused on obtaining g.f.s for walks that begin at the origin: other produced g.f. are purely bonus. To obtain a g.f. that starts at $(0,c)$, input a lower limit of $-c$, and take the g.f. that corresponds to starting at the origin.
		\end{subsection}%

		\begin{subsection}{Meanders}
			What if we do not care where the walk ends, as long as it stays above $y=-c$? Neither Ayyer and Zeilberger \cite{BoundedCase} nor Duchon \cite{Duchon} investigated semi-bounded free walks. In this case, we can actually utilize the irreducible walks we have just created. Recall the irreducible g.f.s in Eqns. \eqref{eqn:IrreducibleDown} and \eqref{eqn:IrreducibleUp}.
			
			Let $k_{a}$ denote the g.f. for nonnegative meanders that begin at $(0,a)$, restricted to step set $\S$. Then
			\beq
				k_0	=	1	+	\left(g_{0,0}+\sum_{(x,0)\in\S}t^x\right)k_0	+	\sum_{(x,y)\in\S;y>0}t^x k_{y-1}.
			\eeq
			The walk can be stationary ($+1$), it can return to the $x$-axis ($g_{0,0}+\sum_{(x,0)\in\S}t^x$) and continue as a new meander ($k_0$), or we take the first step ($t^x$) and continue as a new meander that never returns to the $x$-axis ($k_{y-1}$). $g_{0,0}$, the g.f. for irreducible walks that return to the $x$-axis, already has a description from our previous work. We need only describe the new $k_i$.\footnote{If they exist.}
			\beq
				k_a	=	\sum_{i=0}^{a-1}g_{a,i}k_0	+	k_0	=	\left(\sum_{i=1}^{a}g_{i,0}	+	1\right)k_0.
			\eeq
			The meander can drop down to any lower altitude $(g_{a,i}=g_{a-i,0})$ and then continue as a new meander ($k_0$), never dropping further. Or the meander will never go below $y=a$ so it is equivalent to a meander from the origin ($k_0$).
			
			And since we have already described the irreducible walks earlier, we now have a closed system that we can use to solve for $k_0$. To produce the entire system of equations, use {\bf SemiBoundedScoringPathsEqsVars}. Again, we use {\bf Groebner[Basis]} to find a minimal polynomial: {\bf SemiBoundedScoringPaths}. To find the minimal polynomials for ALL\footnote{Including all of the irreducible and specific altitude walks from the previous section.} of the variables, use {\bf AllSemiBoundedScoringPaths}, though this will likely take a while. For a specific $k_a$, use {\bf SpecificSemiBoundedScoringPaths}. For any $g_{a,b},f_{a,b}$, it is best\footnote{And necessary. Compatibility was removed for ease of use.} to use {\bf SpecificEqualSemiBoundedScoringPaths}.
			\begin{example}
				For comparison, we use the same step set as in Example \ref{exa:EqualSemiBounded}: \linebreak$\S=\{[0,-3],[1,-2],[2,0],[3,1]\}$. The g.f., $K$, for the number of meanders that do not go below $y=-1$ and have step set $\S$ is found to satisfy, using {\bf SemiBoundedScoringPaths($\S$,-1,t,K)},
				\begin{align*}
					{t}^{19} \left( {t}^{2}+t+1 \right)& {K}^{4}	+	{t}^{12} \left( 4{t}^{6}+{t}^{4}+3{t}^{3}+6{t}^{2}+t-3 \right) {K}^{3}\\
					&	+	3{t}^{6} \left( 2{t}^{9}-2{t}^{8}+{t}^{7}+4{t}^{6}+2{t}^{5}-2{t}^{4}-{t}^{3}+3{t}^{2}-1 \right) {K}^{2}\\
					&	+	\left( 4{t}^{12}-8{t}^{11}+7{t}^{10}+11{t}^{9}-7{t}^{8}-3{t}^{7}+7{t}^{6}+6{t}^{5}-6{t}^{4}-2{t}^{3}+4{t}^{2}-1 \right) K\\
					&	+	{t}^{9}-3{t}^{8}+4{t}^{7}+2{t}^{6}-6{t}^{5}+4{t}^{4}+3{t}^{3}-3{t}^{2}+1	=	0.
				\end{align*}
				Sadly, the coefficients do not factor as nicely as in Example \ref{exa:EqualSemiBounded}.
			\end{example}
			
			Now that we have a method of describing the g.f. for meanders, let us compare how fast it is for enumeration.
			\begin{example}[Speed Enumeration - Meanders]
				We use the simpler step set \linebreak$\S=\{[1,-2],[1,-1],[1,0],[1,1],[1,2]\}$ and let $K$ be the g.f. for the number of nonnegative meanders with step set $\S$. Then\footnote{Using {\bf SemiBoundedScoringPaths($\S$,0,t,K)}.}
				\begin{equation}
					t^2(5t-1)^2 K^4	+	t(5t-1)^2 K^3	+	3t(5t-1) K^2	+	(5t-1) K	+	1	=	0.
					\label{eqn:SpeedMeander}
				\end{equation}
				A truncated solution is 
				\beq
					K	=	1+3t+12t^2+51t^3+226t^4+1025t^5+4724t^6+22022t^7+103550t^8+490191t^9+2333057t^{10}.
				\eeq
				
				We repeat our analysis of differing methods of enumeration from Example \ref{exa:EqualSemiBoundedSpeed}.
				\ben
					\item	Iterating a vector of solutions using {\bf SemiBoundedScoringPathsSeries}.%
					\item	Iterating a fixed point solution after solving for the polynomial (Eqn.\ \eqref{eqn:SpeedMeander}) with {\bf SemiBoundedScoringPaths}.\footnote{All of this can be accomplished with the one function {\bf SemiBoundedScoringPathsCoefficients}.}%
					\item	Using Maple's {\bf taylor} function on the minimal polynomial.\footnote{Requires replacing $K$ by $\_Z$ and using {\bf RootOf(p)}.}%
					\item	Enumerating using brute-force recursion and Maple's option remember: \linebreak{\bf SpecificSemiBoundedScoringPathsNumber}.%
				\een
				\begin{table}[H]\begin{center}
					\caption{1000 term Meander Enumeration}
					\begin{tabular}{|c|cccc|}
						\hline
						Method	&	Memory Used	&	Memory Allocation	&	CPU Time	&	Real Time\\
						\hline
						Vector Set-Up	&	21.07KiB	&	0 bytes	&	233$\mu$s	&	233$\mu$s\\
						Vector Iterating	&	2.18TiB	&	15.29MiB	&	102.64m	&	84.48m\\
						\hdashline
						Total	&	2.18TiB	&	15.29MiB	&	102.64m	&	84.48m\\
						\hline
						Polynomial	&	6.37MiB	&	0 bytes	&	43.00ms	&	43.03ms\\
						Single Iterating	&	304.99GiB	&	30.58MiB	&	13.88m	&	11.64m\\
						\hdashline
						Total	&	305.00GiB	&	30.58MiB	&	13.88m	&	11.64m\\
						\hline
						Polynomial	&	6.37MiB	&	0 bytes	&	43.00ms	&	43.03ms\\
						taylor	&	1.85GiB	&	490MiB	&	11.33	&	10.62s\\
						\hdashline
						Total	&	1.86GiB	&	490MiB	&	11.37s	&	10.66s	\\
						\hline
						Brute-Force Recursion	&	1.08GiB	&	0.52GiB	&	7.93s	&	7.465s\\
						\hline
					\end{tabular}
					\label{tab:MeanderEnumeration1000}
				\end{center}\end{table}
				\label{exa:MeanderSpeed}
			\end{example}
			Again, there is an even faster recursive formula. See Example \ref{exa:MeanderConverted}.
		\end{subsection}%
	
		\label{sec:SemiBounded}
	\end{section}%

	\begin{section}{Guess-and-Check Method}
		Zeilberger provided a guess-and-check method and Maple package, {\bf W1D} \cite{GuessandCheck}, for finding algebraic expressions. In the semi-bounded case where all steps have $x$-step 1, Phillippe Duchon guaranteed that the results are algebraic. Thus, the guess-and-check method WILL work, eventually, if you set the search parameter high enough. For semi-bounded cases with differing $x$-steps and unbounded cases, there is no guarantee that guessing will eventually work \cite{nonconverging}. Though you may get lucky and produce an algebraic equation that has the minimal polynomial as a root.
		
		Another advantage of this paper's method over guess-and-check is that this method is typically much much faster. Setting up the equations takes a set amount of time that is linear in the maximum step size, $|\S|$, and the lower limit. The time sink comes in finding a Groebner basis, but this is still typically faster.
		
		{\bf Empir} will take a list of the first few coefficients of the g.f. $F$ and attempt to find an algebraic equation that $F$ satisfies by guessing the degree and trying to solve for the coefficients. {\bf EmpirF} will do this faster\footnote{Almost always.} by utilizing the {\bf gfun} package. Both {\bf Empir} and {\bf EmpirF} require enumerating the first few terms.\\
		
		\begin{example}[Excursion Time Comparison]
			Consider looking for an algebraic equation for the g.f. for nonnegative excursions with step set $\S=\{[1,-2],[1,-1],[1,0],[1,1],[1,2]\}$. It takes\footnote{The actual minimal polynomial is shown in Example \ref{exa:EqualSemiBoundedSpeed}.}
			\begin{table}[H]\begin{center}
				\caption{Finding Excursion Minimal Polynomial $\S=\{[1,-2],[1,-1],[1,0],[1,1],[1,2]\}$}
				\begin{tabular}{|c|cccc|}
					\hline
					Method	&	Memory Used	&	Memory Allocation	&	CPU Time	&	Real Time\\
					\hline
					{\bf ESBSP}	&	2.63MiB	&	0 bytes	&	19.23ms	&	18.63ms\\
					\hline
					{\bf Empir}	&	91.61MiB	&	5.60MiB	&	734ms	&	735ms\\
					\hline
					{\bf EmpirF}	&	4.33MiB	&	0 bytes	&	33.87ms	&	35.90ms\\
					\hline
				\end{tabular}\end{center}
				\label{tab:EqualSemiBoundedMinPolynomial}
			\end{table}
			{\bf EqualSemiBoundedScoringPaths} ({\bf ESBSP}) takes about $3\%$ of the time and memory that {\bf Empir} requires and $60\%$ of the time and memory as {\bf EmpirF}.%
			\label{exa:SemiBoundedSpeedGuess}
		\end{example}
		
		\begin{example}[Meander Time Comparison]
			Now try to find an algebraic equation satisfied by the g.f. for nonnegative meanders with step set $\S=\{[1,-2],[1,-1],[1,0],[1,1],[1,2]\}$.\footnote{The minimal polynomial is shown in Example \ref{exa:MeanderSpeed}.}
			\begin{table}[H]\begin{center}
				\caption{Finding Meander Minimal Polynomial $\S=\{[1,-2],[1,-1],[1,0],[1,1],[1,2]\}$}
				\begin{tabular}{|c|cccc|}
					\hline
					Method	&	Memory Used	&	Memory Allocation	&	CPU Time	&	Real Time\\
					\hline
					{\bf SBSP}	&	6.37MiB	&	0 bytes	&	43.00ms	&	43.03ms\\
					\hline
					{\bf Empir}	&	91.68MiB	&	5.45MiB	&	749ms	&	770ms\\
					\hline
					{\bf EmpirF}	&	4.74MiB	&	0 bytes	&	36.53ms	&	38.63ms\\
					\hline
				\end{tabular}\end{center}
				{\bf SBSP}={\bf SemiBoundedScoringPaths}.
				\label{tab:SemiBoundedMinPolynomial}
			\end{table}
			As expected, this package is still much faster than {\bf Empir}. Interestingly, {\bf EmpirF} appears to be as fast, if not a little faster. {\bf EmpirF} is still handicapped in its range of applications and so the slight speed-up is sacrificed for versatility.
			\label{exa:MeanderSpeedGuess}
		\end{example}
		In fact, adding $[1,3]$ to the step set already makes {\bf EmpirF} fail for the preset search bound. The minimal polynomial has degree and order 10 in that case: easily found by {\bf SemiBoundedScoringPaths} in $1/4$s.
		
		The other bonus is that this paper's method allows for any size $x$-step. {\bf W2D} could potentially be used to solve this problem once it has been restricted to look for the coefficient of $x^{n}\cdot y^{0}$. However, this would make the already slow guessing method even slower.
	\end{section}%

	\begin{section}{Algebraic to Recursive}

		\begin{subsection}{Conversion}
			There is a classical method for deducing from the algebraic function satisfying the g.f. a linear recurrence with polynomial coefficients satisfied by the coefficients of the g.f. in question. See Chapter 6 of
\tbl The Concrete Tetrahedron\tbr\ by Kauers and Paule \cite{ConcreteTetrahedron}.

			The method is implemented in the Maple package gfun \cite{salvy1994gfun} and also in Zeilberger's Maple package {\bf SCHUTZENBERGER}, that is used here. 
			The {\bf SCHUTZENBERGER} package also contains {\bf EmpirF} for obtaining the minimal polynomial, but we now have a much better method of producing the minimal polynomial. To convert an algebraic formula to a recurrence formula, use {\bf algtorec}. Let $B(n)$ denote the number of walks of length $n$.\\
			Interestingly, sometimes a larger (non-minimal) polynomial produces a better (lower-order) recurrence.\\
			
			\begin{example}[Better Excursion Recursion]
				Let $\S=\{[1,-2],[1,-1],[1,0],[1,1],[1,2]\}$. Let $F$ denote the g.f. for nonnegative excursions.\footnote{The minimal polynomial is given in Example \ref{exa:EqualSemiBoundedSpeed}.} Then its coefficients satisfy
				\begin{align}
					0	&=	3125(n+1)(n+2)(n+3)(n+4)B(n)\nonumber\\
						&\tab	-	250(n+4)(n+3)(n+2)(27n+122)B(n+1)\nonumber\\
						&\tab	+	25(n+4)(n+3)(107n^2+1457n+4316)B(n+2)\nonumber\\
						&\tab	+	10(n+4)(304n^3+3233n^2+9864n+6513)B(n+3)\nonumber\\
						&\tab	-	(2821n^4+56794n^3+425771n^2+1407974n-1731540)B(n+4)\nonumber\\
						&\tab	+	2(n+7)(413n^3+6986n^2+39356n+73830)B(n+5)\nonumber\\
						&\tab	-	(n+8)(n+7)(99n^2+1241n+3900)B(n+6)\nonumber\\
						&\tab	+	2(2n+15)(n+9)(n+8)(n+7)B(n+7).
					\label{eqn:ImprovedRecursion}
				\end{align}
				The following are the time and memory requirements for various stages of enumerating. We need to create the minimal polynomial with {\bf EqualSemiBoundedScoringPaths}, convert it to an improved recursive formula, Eqn.\ \eqref{eqn:ImprovedRecursion}, for the coefficients of $F$ with {\bf algtorec}, and then enumerate with {\bf SeqFromRec}.
				\begin{table}[H]\begin{center}
					\caption{Enumerating Excursions More Efficiently}
					\begin{tabular}{|c|cccc|}
						\hline
						Method	&	Memory Used	&	Memory Allocation	&	CPU Time	&	Real Time\\
						\hline
						{\bf ESBSP}	&	2.63MiB	&	0 bytes	&	19.23ms	&	18.63ms\\
						{\bf algtorec}	&	35.84MiB	&	3.12MiB	&	267.93ms	&	262.40ms\\
						{\bf SeqFromRec}	&	123.05MiB	&	24MiB	&	375.27ms	&	375.37ms\\
						\hdashline
						Total	&	161.52MiB	&	27.12MiB	&	662.43ms	&	656.40ms\\
						\hline
					\end{tabular}\end{center}
					{\bf ESBSP}={\bf EqualSemiBoundedScoringPaths}
					\label{tab:ExcursionEfficientEnumeration}
				\end{table}
				As expected, this streamlined recurrence is much faster and less memory-intensive than the basic recurrence in Example \ref{exa:EqualSemiBoundedSpeed}. There is an up-front cost for creating the improved recurrence, but if the goal is to enumerate enough terms, it can be worth it. In this case, enough is less than 500.%
				\label{exa:EqualSemiBoundedConverted}
			\end{example}
			Before {\bf FindProperRoot} was implemented to automatically parse the minimal polynomial, we converted a larger polynomial into a sixth-order recurrence in about 2.4 seconds.\\
			
			\begin{example}[Better Meander Recursion]
				Let $\S=\{[1,-2],[1,-1],[1,0],[1,1],[1,2]\}$ and $K$ denote the g.f. for nonnegative meanders.\footnote{The minimal polynomial is given in Example \ref{exa:MeanderSpeed}.} The coefficients now satisfy
				\begin{align*}
					0	&=	625(n+1)(n+2)(n+3)(n+4)B(n)\\
						&\tab	-	250(5n+21)(n+4)(n+3)(n+2)B(n+1)\\
						&\tab	+	50(n+4)(n+3)(7n^2+95n+270)B(n+2)\\
						&\tab	+	20(n+4)(32n^3+367n^2+1365n+1620)B(n+3)\\
						&\tab	-	(n+5)(463n^3+6691n^2+32442n+52704)B(n+4)\\
						&\tab	+	2(n+5)(n+6)(53n^2+593n+1674)B(n+5)\\
						&\tab	-	4(n+5)(2n+13)(n+7)(n+6)B(n+6).
				\end{align*}
				Oddly it is a lower order recurrence despite the slightly higher complexity of describing meanders. The following are the time and memory requirements for various stages of enumerating.
				\begin{table}[H]\begin{center}
					\caption{Enumerating Meanders More Efficiently}
					\begin{tabular}{|c|cccc|}
						\hline
						Method	&	Memory Used	&	Memory Allocation	&	CPU Time	&	Real Time\\
						\hline
						{\bf SBSP}	&	6.37MiB	&	0 bytes	&	43.00ms	&	43.03ms\\
						{\bf algtorec}	&	110.77MiB	&	29.39MiB	&	659.73ms	&	633.10ms\\
						{\bf SeqFromRec}	&	89.58MiB	&	4MiB	&	284.57ms	&	257.23ms\\
						\hdashline
						Total	&	206.72MiB	&	33.39MiB	&	987.30ms	&	933.36ms\\
						\hline
					\end{tabular}\end{center}
					{\bf SBSP}={\bf SemiBoundedScoringPaths}
					\label{tab:MeanderEfficientEnumeration}
				\end{table}
				This improved recurrence is about 8 times as fast as the basic recurrence in Example \ref{exa:MeanderSpeed}. It is not quite as much of a savings as the case of excursions, but it is still extremely good.
				\label{exa:MeanderConverted}
			\end{example}
	
		\end{subsection}%

		\begin{subsection}{Searching}
			There is an alternative to converting the minimal polynomial into a linear recurrence. Since we know that this will be possible, we could simply guess at the form of the recurrence and use a suitable number of starting values to determine the coefficients. {\bf Findrec} will accomplish this guessing method. The problem is that we do not know an upper bound for the order and degree. Setting the bound extremely high (or searching until a recurrence is found) would suffice. However, this is not ideal.
			\begin{example}
				Let $\S=\{[1,2],[1,-3]\}$. The minimal polynomial is given later in Section \ref{subsubsec:Duchon}. Let $B(n)$ denote the number of nonnegative excursions of length $5n$. $B(n)$ was found to satisfy a $4^{th}$ order, degree $11$ recurrence relation using {\bf Findrec} and a $7^{th}$ order, degree $9$ recurrence relation using {\bf algtorec}.\footnote{2 trials.} Both relations\footnote{The {\bf Findrec} recurrence matches with Andrew Lohr's calculation as expected since they are computing the same result.} are fairly large and so are not included here.\footnote{They are available at {\bf math.rutgers.edu/$\sim$bte14/Articles/ScoringPaths/23Recurrence.txt}.} We compare and contrast the methods in the table below.
				\begin{table}[H]\begin{center}
					\caption{Alternative Recurrence Step Set $\{[1,2],[1,-3]\}$}
					\begin{tabular}{|c|cccccc|}
						\hline
						Method	&	Degree	&	Order	&	Memory Used	&	Allocated	&	CPU Time	&	Real Time\\
						\hline
						{\bf algtorec}	&	9	&	7	&	165.84GiB	&	2.46GiB	&	46.99h	&	12.06h\\
						\hline
						{\bf Findrec}	&	11	&	4	&	0.50GiB	&	388.01MiB	&	4.47s	&	4.31s\\
						\hline
					\end{tabular}\end{center}
					\label{tab:AlternativeRecurrence23}
				\end{table}
				Simply guessing at the form appears to be the {\it MUCH} better choice here.
			\end{example}
			
			\begin{example}
				Let $\S=\{[1,-2],[1,-1],[1,0],[1,1],[1,2]\}$. See Example \ref{exa:EqualSemiBoundedSpeed} for the minimal polynomial. Let $B(n)$ denote the number of nonnegative excursions of length $n$. Then
				\begin{align*}
					0	&=	2 \left( n+4 \right)  \left( 2n+11 \right)  \left( n+7 \right)  \left( n+6 \right) B(n+5)\\
						&\tab	-	\left( n+6 \right)  \left( 43{n}^{3}+597{n}^{2}+2738n+4142 \right) B(n+4)\\
						&\tab	+\left( 124{n}^{4}+2110{n}^{3}+13305{n}^{2}+36815n+37686 \right) B(n+3)\\
						&\tab	-	5 \left( n+3 \right)  \left( 2{n}^{3}-22{n}^{2}-305n-726 \right) B(n+2)\\
						&\tab	-25 \left( n+3 \right)  \left( n+2 \right)  \left( 8{n}^{2}+72n+159 \right) B(n+1)\\
						&\tab	+	125 \left( n+5 \right)  \left( n+3 \right)  \left( n+2 \right)  \left( n+1 \right)B(n).
				\end{align*}
				This was found with {\bf Findrec} while {\bf algtorec} found
				\begin{align*}
					0	&=	2 \left( 2n+15 \right)  \left( n+9 \right)  \left( n+8 \right)  \left( n+7 \right) B(n+7)\\
						&\tab-	\left( n+8 \right)  \left( n+7 \right)  \left( 99{n}^{2}+1241n+3900 \right)B(n+6)\\
						&\tab	+2 \left( n+7 \right)  \left( 413{n}^{3}+6986{n}^{2}+39356n+73830 \right) B(n+5)\\
						&\tab	-	\left( 2821{n}^{4}+56794{n}^{3}+425771{n}^{2}+1407974n+1731540 \right)B(n+4)\\
						&\tab	+10 \left( n+4 \right)  \left( 304{n}^{3}+3233{n}^{2}+9864n+6513 \right) B(n+3)\\
						&\tab	+	25 \left( n+4 \right)  \left( n+3 \right)  \left( 107{n}^{2}+1457n+4316 \right)B(n+2)\\
						&\tab	-250 \left( n+4 \right)  \left( n+3 \right)  \left( n+2 \right)  \left( 27n+122 \right) B(n+1)\\
						&\tab	+	3125 \left( n+1 \right)  \left( n+2 \right)  \left( n+3 \right)  \left( n+4 \right)B(n).
				\end{align*}
				The comparison in methods is given below.
				\begin{table}[H]\begin{center}
					\caption{Alternative Recurrence Step Set $\{[1,-2],[1,-1],[1,0],[1,1],[1,2]\}$}
					\begin{tabular}{|c|cccccc|}
						\hline
						Method	&	Degree	&	Order	&	Memory Used	&	Allocated	&	CPU Time	&	Real Time\\
						\hline
						{\bf algtorec}	&	4	&	7	&	38.32MiB	&	42.61MiB	&	355.60ms	&	354.44ms\\
						\hline
						{\bf Findrec}	&	4	&	5	&	62.38MiB	&	28.00MiB	&	618.07ms	&	630.87ms\\
						\hline
					\end{tabular}\end{center}
					\label{tab:AlternativeRecurrence21012}
				\end{table}
				Conversion is actually more efficient, though not ideal, in this case. Though the difference is not as extreme as that shown in Table \ref{tab:AlternativeRecurrence23}.
			\end{example}
			
			With smaller step sizes, conversion appears to be \tbl better\tbr. This is likely due to the intensive task of converting higher degree polynomials, while guessing avoids that hurdle. For small enough examples, guessing is a larger search space than conversion.

			Because of the way {\bf Findrec} is programmed, it is guaranteed that a found recurrence has order less than or equal to that of {\bf algtorec}. Depending on the goal when enumerating, lower degree (save computation) or lower order (save memory) is optimal.
			
			One can use {\bf algotrec}'s order and degree as a proper bound for use in {\bf Findrec}. But instead of directly converting the minimal polynomial, it may be possible to find a bound on the order and degree of a recurrence based on the degrees of the g.f. and $t$ in the minimal polynomial. We may also be able to derive an upper bound based on the number of variables in the closed system we created.
			
			I was unable to prove any useful bounds, and these examples are evidence against bounding the order and degree directly by the degrees of the g.f. and $t$. All three examples have minimal polynomials with smaller degree (and total degree) than the order of the found recurrences. In fact, we know that the sum of the degree and order of any recurrence has to be greater than the total degree of the minimal polynomial in these cases since we searched everything lower.\\
			
			The first example showed that {\bf Findrec} can be an {\it EXTREME} improvement over {\bf algtorec} and the second example showed when {\bf algtorec} can just edge out {\bf Findrec}. The following example demonstrates that {\bf algtorec} can be the significantly better choice.
	
			\begin{example}
				Let $\S=\{[1,-1],[3,-1],[1,0],[3,0],[2,1],[1,2],[2,2]\}$. The minimal polynomial, letting $K$ denote the g.f. of meanders with step set $\S$, is 
				\begin{align*}
					1	+	(t^3+5t^2+4t-1)K	+	t(4t+3)(2t^3+2t^2+3t-1)K^2	+	t(t+1)(2t^3+2t^2+3t-1)K^3.
				\end{align*}
				The recurrences found by conversions and searching are much too large to include here. For their information, see \\{\bf math.rutgers.edu/$\sim$bte14/Articles/ScoringPaths/RecurrenceOutputFile}.
				\begin{table}[H]\begin{center}
					\caption{Alternative Recurrence Step Set $\{[1,-1],[3,-1],[1,0],[3,0],[2,1],[1,2],[2,2]\}$}
					\begin{tabular}{|c|cccccc|}
						\hline
						Method	&	Degree	&	Order	&	Memory Used	&	Allocated	&	CPU Time	&	Real Time\\
						\hline
						Polynomial	&		&		&	2.55MiB	&	24.00MiB	&	20.80ms	&	24.50ms\\
						{\bf algtorec}	&	3	&	31	&	73.92MiB	&	11.39MiB	&	491.03ms	&	475.13ms\\
						\hdashline
						Total	&		&		&	76.47MiB	&	35.39MiB	&	511.83ms	&	499.63ms\\
						\hline
						161 terms	&		&		&	24.11MiB	&	0 bytes	&	190.00ms	&	193.00ms\\
						{\bf Findrec}	&	5	&	15	&	2.46GiB	&	504MiB	&	23.30s	&	21.20s\\
						\hdashline
						Total	&		&		&	2.48GiB	&	504MiB	&	23.49s	&	21.39s\\
						\hline
					\end{tabular}\end{center}
					\label{tab:AlternativeRecurrenceS}
				\end{table}
				Conversion is surprisingly 43 times faster and has roughly $1/32^{nd}$ of the memory requirements!
				
				We chose to enumerate 161 terms of the sequence because that was sufficient, once we knew the degree and order of the converted recurrence, to guarantee {\bf Findrec} would encounter a solution. Typically, we do not use both {\bf Findrec} and {\bf algtorec}; we should simply find a set large bound of numbers. But enumerating was a small portion of {\bf Findrec}'s time so it does not matter too much to this example.
			\end{example}
			We have used meanders instead of excursions for our $3^{rd}$ example, but they are similar enough to compare. The important facets of the problem are the degree of the minimal polynomial and the number of steps.
			
			It would appear that increasing the number of steps does not affect the runtime of {\bf algtorec} as much as the runtime of {\bf Findrec}. The runtime of {\bf algtorec} is highly dependent on the degree of the minimal polynomial.\footnote{The degree of $t$ is much less relevant.} The degree of the minimal polynomial is generally correlated with the number of variables in the system, though not necessarily equal or bounded one way or the other. Example \ref{exa:EqualSemiBoundedSpeed} has 7 variables in the closed system and a minimal polynomial of degree only 4. Section \ref{subsec:2stepExamples} shows that a step set of $\{[1,3],[1,-5]\}$ for nonnegative excursions yields a closed system with 18 variables yet the minimal polynomial has degree 56!
			
			The most important aspects of the step set to {\bf algtorec} runtime are then the maximum and minimum\footnote{Most negative.} $y$-steps since these dictate how many other walks we must consider.
			
			{\bf Findrec} is much weaker with larger step sets. Inherently, the recurrences will become more complicated, which means {\bf Findrec} must consider many more terms in its guesses. The worst case (most number of terms to address) will occur when the order and degree of our search are the same.\\
			
			Since {\bf algtorec} works better with minimal polynomials of lower degree, and finding the minimal polynomial is generally very fast (see Example \ref{exa:MeanderSpeed}), I have chosen to have programs such as {\bf PaperSemiBounded} and {\bf BookEqualSemiBounded} use {\bf Findrec} when the minimal polynomial has degree $\ge8$ and {\bf algtorec} otherwise. The cutoff was arbitrarily chosen based on a few random examples. A rigorous cutoff would be useful for more widespread implementation.
	
		\end{subsection}%
	
		\label{sec:AlgebraictoRecursive}
	\end{section}%

	\begin{section}{Unbounded}
		We remove the lower bound in this section and transition from analyzing excursions to analyzing bridges. I was unable to find any current literature that analyzes the unbounded case with generating functions; thus, I am led to believe that this section is novel work. Again, all walks are assumed to be constructed from a step set $\S$.

		\begin{subsection}{Walking \tbl Anywhere\tbr}
			Enumerating walks that go anywhere {\it can} be a trivial task. First off, the $y$-values do not matter except for describing multiple steps with the same $x$-step. If all of the steps have $x$-value $m>0$, then the number of walks of length $m\cdot n$ is simply $|\S|^n$. If we have steps with varying $x$-values, then the problem becomes finding combinations of the $x$-values that sum to $n$. In general, the g.f. will be
			\beq
				\frac{1}{1-\sum_{(x,y)\in\S}t^x}.
			\eeq
			This g.f. is produced by {\bf UnBoundedScoringPaths}. The terms are generated by (using Taylor series) {\bf UnBoundedScoringPathsCoefficients} or (via recursion) {\bf UnBoundedScoringPathsNumber}. Using Maple's taylor series expansion is typically faster, but both are extremely quick at enumerating the first 1000 terms.\footnote{By utilizing Maple's {\bf option remember} in the recursive case.}
		\end{subsection}%

		\begin{subsection}{Returning to the $x$-axis}
			A lot of the work we have done in the semi-bounded case will prove useful here. We cannot use the exact same method as in semi-bounded Section \ref{sec:SemiBounded}. Suppose we tried describing a walk with a negative change in altitude. The first and last steps could both be positive. Then we would need to describe a walk that has a larger negative change in altitude. And so on.
			
			We try a slightly different method. Let $G$ denote the g.f. for walks with step set $\S$ that begin at $(0,0)$ and end on the $x$-axis: a bridge. We choose to introduce another g.f., which we denote $h_{i}$, for the number of walks from $(0,i)$ to $(n,0)$ that do not touch the $x$-axis beforehand, with the exception of $h_0$ starting on the $x$-axis.\footnote{$h_0$ does NOT include stationary walks nor walks that are steps directly to the right. For that, one will want to use $h=h_0+1+\sum_{(x,0)\in\S}t^x$.} Recall the definitions of $g_{a,b}$ and $f_{a,b}$. The first equation is similar to $f_{0,0}$: 
			\beq
				G	=	1+\left(h_0+\sum_{(x,0)\in\S}t^x\right)G.
			\eeq
			A walk can be stationary, or it returns to the $x$-axis after some number of steps, at which point it can take another $G$-walk.
			\beq
				h_0	=	2g_{0,0}+\sum_{(x,y)\in\S;y\le-2}\sum_{i=1}^{-y-1}g_{0,i} t^x h_{i+y}+\sum_{(x,y)\in\S;y\ge2}\sum_{i=1}^{y-1}g_{y-i,0} t^x h_i.
			\eeq
			An $h_0$-walk can be purely positive ($g_{0,0}$) or purely negative ($g_{0,0}$).\footnote{This is also partially why irreducible walks do not include steps to the right. If $g_{0,0}$ did, then I would have to write $h_0=2g_{0,0}-\sum_{(x,0)\in\S}t^x$, which is less elegant.} Note that walks below the $x$-axis are in bijection with walks above the $x$-axis by reversing the order of steps.\footnote{Equivalent to rotating the walk $180^\circ$.} We can then use $g_{0,-a}=g_{a,0}$ and $f_{-a,-b}=f_{b,a}$. It is also possible for the $h$-walk to cross the $x$-axis without touching it. This would involve an irreducible walk ($g_{0,i}$), a step across the $x$-axis ($t^x$), and another walk ($h_{i+y}$) to return to the $x$-axis. The previous sentence described crossing the $x$-axis from above to below. We could instead cross from below to above. This would consist of an irreducible walk ($g_{0,i-y}=g_{y-i,0}$), a step across the $x$-axis ($t^x$), and another walk ($h_i$) to return to the $x$-axis.
			
			We already have equations to describe the irreducible walks. We now need to describe the other $h$-walks.
			\begin{align*}
				j>0	&\tab	h_j=	g_{j,0}+\sum_{(x,y)\in\S;y\le-2}\sum_{i=1}^{-y-1}f_{j-1,i-1} t^x h_{i+y},\\
				j<0	&\tab	h_j=	g_{0,-j}+\sum_{(x,y)\in\S;y\ge2}\sum_{i=1}^{y-1}f_{y-i-1,-j-1} t^x h_i.
			\end{align*}
			An $h_{j>0}$-walk can either be a walk \tbl directly\tbr\ to the $x$-axis ($g_{j,0}$) or it consists of an arbitrary walk that does not touch the $x$-axis ($f_{j-1,i-1}$), followed by a step across ($t^x$), and another walk ($h_{i+y}$) to return to the $x$-axis. Similarly an $h_{j<0}$-walk can either be a walk \tbl directly\tbr\ to the $x$-axis ($g_{j,0}=g_{0,-j}$) or it consists of an arbitrary walk that does not touch the $x$-axis ($f_{j+1,i-y+1}=f_{y-i-1,-j-1}$), followed by a step across ($t^x$), and another walk ($h_{i}$) to return to the $x$-axis.
			
			A brief lapse in concentration allowed me to write the alternative (and much simpler) relation:
			\beq
				h_0	=	2g_{0,0}	+	\sum_{(x,y)\in\S;y\ne0}t^x h_y.
			\eeq
			The reason for shunning this \tbl simpler\tbr\ description is that the $h_y$ may be irreducible walks directly back to the $x$-axis, and thus double-count a walk from $g_{0,0}$.\\
		
			The system of equations will be closed since the index of $h_i$-walks is bounded by $\max_{(x,y)\in\S}|y|-1$. All of the equations and variables are produced in {\bf EqualUnBoundedScoringPathsEqsVars}. Again, we use the {\bf Basis} procedure in the {\bf Groebner} package to find a polynomial $p$ such that $p(G)=0$ in terms of formal power series. To produce this polynomial, use {\bf EqualUnBoundedScoringPaths}. To produce the minimal polynomial for $h_j$, use {\bf SpecificUnBoundedScoringPaths}.\\
			
			We can try to produce the coefficients of the taylor series using the iteration $G_0\to p(G_0)+G_0$. However, because of convergence issues, this will typically not work. It is also potentially not ideal because it requires finding a Groebner basis, which can be time-consuming. However, iterating a vector of solutions \`a la the semi-bounded Section \ref{sec:SemiBounded} does work. To check these values, use {\bf SpecificUnBoundedScoringPathsNumber} to verify $G$ and {\bf SpecificIrreducibleUnBoundedScoringPathsNumber} to verify $h_{j}$.
		\end{subsection}%

		\begin{subsection}{Alternative Method}
			\label{subsec:AlternativeMethod}
			There is an alternative method for enumerating unbounded walks in the specific case that all of the steps have $x$-step 1 \cite{AZd}. The functions themselves are taken from the Maple package {\bf EKHAD} by Zeilberger. I will start with an example.\\
			
			\begin{example}
				Consider the function $G(t,n)=\left(t^{-2}+t^{-1}+t+t^2\right)^n$. This represents a possible step set of $\S=\{[1,-2],[1,-1],[1,1],[1,2]\}$. A specific monomial in the expansion of $G(t,n)$ comes from all the ways of picking a power of $t$ from each term; this is equivalent to picking which step to take and in which order.
				\ben
					\item	$G(t,0)=1$. There is only 1 walk: the stationary walk.
					\item	$G(t,1)=t^{-2}+t^{-1}+t+t^2$. There are 4 possible steps each leading to a different final altitude.
					\item	$G(t,2)=t^{-4}+2t^{-3}+t^{-2}+2t^{-1}+4+2t+t^2+2t^3+t^4$. There are 4 ways to return to the $x$-axis after 2 steps.
				\een
				$G(t,n)$ enumerates all the ways to change altitude by $c$ in the coefficient of the $t^c$ monomial.
			\end{example}
			In general, for a step set $\S$, let $G(t,n)=\left(\sum_{(x,y)\in\S}t^y\right)^n$. The procedure {\bf AZd(A,t,n,N)} from package {\bf EKHAD} will give a recurrence for the (contour around 0) integral of $A$ (with respect to $t$) under hypergeometric assumptions, i.e., {\bf AZd} returns the residue. And if we represent $A$ as a power series, that means we extract the $t^{-1}$ term. So to find a recurrence for the number of walks that return to the $x$-axis, we will use $A=G(t,n)/t$.\\
			
			\begin{example}
				The input {\bf AZd($\left(t^{-2}+t^{-1}+t+t^2\right)^n$,t,n,N)} yields as output:
				\beq
					18(2n+1)(5n+8)(n+1)+(n+1)(35n^2+91n+54)N-2(n+2)(2n+3)(5n+3)N^2.
				\eeq
				This translates as, if we let $B(n)$ denote the number of walks that return to the $x$-axis after $n$ steps,
				\begin{align}
					18(2n+1)(5n+8)(n+1)B(n)	&+	(n+1)(35n^2+91n+54)B(n+1)\nonumber\\
										&-2(n+2)(2n+3)(5n+3)B(n+2)	=	0.
					\label{eqn:AZdrecurrence}
				\end{align}
			\end{example}
			We can find a recursion for {\it any} change in altitude, not just bridges. If we are interested in walks that change in elevation by $c$, then we should use $A=G(t,n)/t^{c+1}$.
		\end{subsection}%

		\begin{subsection}{Selected Step Sets}
			The following examples are all produced by the one call:
			\begin{center}{\bf EqualUnBoundedScoringPaths($\S$,t,G)}.\end{center}
			\begin{example}[Old-Time Basketball]
				We will follow up the previous example by using this paper's g.f. method of dynamical programming with the same step set. Let \linebreak$\S=\{[1,-2],[1,-1],[1,1],[1,2]\}$. Let $G$ denote the g.f. for the number of unbounded bridges subject to step set $\S$. Then $G$ satisfies
				\beq
					(9t+4)(4t-1)^2G^4	-	2(3t-2)(4t-1)G^2	+	t	=	0.
				\eeq
				And the coefficients satisfy (using {\bf algtorec})
				\begin{align*}
					108(n+1)(2n+1)B(n)&+(78n^2+246n+216)B(n+1)\nonumber\\
					&-(n+2)(17n-9)B(n+2)-2(n+3)(2n+5)B(n+3)	=	0.
				\end{align*}
				Interestingly, we have obtained a different recurrence than what was found with the {\bf AZd} method: Eqn.\ \eqref{eqn:AZdrecurrence}. Both recurrences are correct. Typically, {\bf AZd} will produce a lower order recurrence but with higher degree polynomial coefficients.\\
				
				The simple alteration of the step set to $\S'=\{[2,-2],[1,-1],[1,1],[2,2]\}$ makes the alternative method in Section \ref{subsec:AlternativeMethod} ineffectual. However, the method presented in this paper is not bothered in the slightest. Let $G'$ denote the g.f. for the number of unbounded bridges subject to step set $\S'$. Then $G'$ satisfies
				\beq
					\left( 8{t}^{2}+5 \right)  \left( 4t^4-8t^2+1 \right) ^{2} G'^4	+	2 \left( 4{t}^{2}-3 \right)  \left( 4t^4-8t^2+1 \right)G'^2	+	1	=	0.
				\eeq
				Now we only have the option of conversion or guessing to find a recurrence for the coefficients of $G'$. Using {\bf algtorec} took only 69ms and produced:
				\begin{align*}
					32(n+3)(n+2)B(n)+4(5n^2+110n+356)B(n+2)-8(15n^2+160n+439)B(n+&4)\\
						-(59n+438)(n+6)B(n+6)+10(n+8)(n+6)B(n+8)	&=	0.
				\end{align*}
			\end{example}
			The reason to have $x$-step as all 1s versus having $x$-step equal to $y$-step is that the choice changes what $G$ counts. If all $x$-steps are 1, $G$ will count the number of ways to be tied after $n$ total baskets. If $x$-step equals $|y|$-step, then $G'$ counts the number of ways to be tied at a score of $\frac{n}{2}$ to $\frac{n}{2}$. Note that this interpretation means $G'$ has coefficient 0 for all odd powers of $t$; teams cannot be tied after an odd number of points have been scored. We could make the substitution $t\to\sqrt{t}$ in $G'$ and then $G'$ would count the number of ways to be tied at a score of $n-n$.
			
			\begin{example}[Current Basketball]
				Let $\S=\{[1,-3],[1,-2],[1,-1],[1,1],[1,2],[1,3]\}$. Let $G$ denote the g.f. for the number of unbounded bridges subject to step set $\S$. Then $G$ satisfies
				\begin{align*}
					&\left( 8{t}^{2}-68t-27 \right) ^{2} \left( 6t-1 \right) ^{4} \left( 2t+1 \right) ^{4} G^8\\
					&\tab	+	4 \left( 68{t}^{2}+10t-9 \right) \left( 8{t}^{2}-68t-27 \right) \left( 6t-1 \right) ^{3} \left( 2t+1 \right) ^{3} G^6\\
					&\tab	+	2 \left( 9120{t}^{4}+3744{t}^{3}-1264{t}^{2}-212t+135 \right) \left( 6t-1 \right) ^{2} \left( 2t+1 \right) ^{2} G^4\\
					&\tab	+	4 \left( 1216{t}^{4}+832{t}^{3}+4{t}^{2}-46t+7 \right)  \left( 6t-1 \right)  \left( 2t+1 \right) G^2\\
					&\tab	+	 \left( 16{t}^{2}+8t-1 \right) ^{2}	=	0.
				\end{align*}
				{\bf EqualUnBoundedScoringPaths($\S$,t,G)} originally gave a much larger polynomial that $G$ satisfies. After {\bf FindProperRoot} was added to parse the output, the smaller polynomial above was the result. The coefficients satisfy (using {\bf algtorec})
				\begin{align*}
					&	36864 \left( n+1 \right)  \left( n+2 \right)  \left( n+3 \right)B(n)\\
					&\tab	-3072 \left( n+3 \right)  \left( n+2 \right)  \left( 97n+142 \right) B(n+1)\\
					&\tab	-64 \left( n+3 \right)  \left( 4031{n}^{2}+17601n+19504 \right)B(n+2)\\
					&\tab	-16\left(1684n^3+13491n^2+31178n+15240\right)B(n+3)\\
					&\tab	+ 24\left(663n^3+7222n^2+28628n+41563 \right)B(n+4)\\
					&\tab	+ 4\left(467n^3+7011n^2+35842n+62490 \right)B(n+5)\\
					&\tab	-2 \left( n+6 \right)  \left( 115{n}^{2}+1080n+2273 \right)B(n+6)\\
					&\tab	-3 \left( 3n+20 \right)  \left( 3n+19 \right)  \left( n+7 \right)B(n+7)	=	0.
				\end{align*}
				or, using {\bf AZd},
				\begin{align*}
						96 \left( n+1 \right)  \left( n+2 \right)  \left( n+3 \right) & \left( 2058{n}^{3}+20335{n}^{2}+66857n+73300 \right)B(n)\\
						-8 \left( n+2 \right)  \left( n+3 \right)  \bigg(& 201684{n}^{4}+2295356{n}^{3}\\
					&+9540055{n}^{2}+16998380n+10742400 \bigg)B(n+1)\\
						-4 \left( n+3 \right)  \bigg(& 310758{n}^{5}+4313617{n}^{4}+23611469{n}^{3}\\
					&+63712598{n}^{2}+84804508n+44608800 \bigg)B(n+2)\\
						-2 \left( n+3 \right)  \bigg(&41160{n}^{5}+591920{n}^{4}+3361281{n}^{3}\\
					&+9453847{n}^{2}+13262292n+7512000 \bigg)B(n+3)\\
					+3 \left( 3n+11 \right)  \left( 3n+10 \right) & \left( n+4 \right)  \left( 2058{n}^{3}+14161{n}^{2}+32361n+24720 \right)B(n+4)	=	0.
				\end{align*}
				
				Let $\S'=\{[3,-3],[2,-2],[1,-1],[1,1],[2,2],[3,3]\}$. Let $G'$ denote the g.f. for the number of unbounded bridges subject to step set $\S'$. Then $G'$ satisfies
				\begin{align*}
					&\left( 108{t}^{6}-99{t}^{4}-52{t}^{2}-44 \right) ^{2} \left( 4{t}^{6}+4{t}^{4}+8{t}^{2}-1 \right) ^{4}G'^8\\
					&\hspace{0.5em}	+4\left( 36{t}^{6}+29{t}^{4}+24{t}^{2}-20 \right)  \left( 108{t}^{6}-99{t}^{4}-52{t}^{2}-44 \right)  \left( 4{t}^{6}+4{t}^{4}+8{t}^{2}-1 \right) ^{3}G'^6\\
					&\hspace{0.5em}	+2\left( 2160{t}^{12}+4248{t}^{10}+4347{t}^{8}+992{t}^{6}-16{t}^{4}-1184{t}^{2}+976 \right)  \left( 4{t}^{6}+4{t}^{4}+8{t}^{2}-1 \right) ^{2}G'^4\\
					&\hspace{0.5em}	+4\left( 112{t}^{12}+288{t}^{10}+665{t}^{8}+412{t}^{6}+552{t}^{4}-112{t}^{2}+96 \right)\left( 4{t}^{6}+4{t}^{4}+8{t}^{2}-1 \right)G'^2\\
					&\hspace{0.5em}	+ \left( 4{t}^{6}+3{t}^{4}+12{t}^{2}+4 \right) ^{2}	=	0.
				\end{align*}
				Again, {\bf EqualUnBoundedScoringPaths($\S'$,t,G')} originally gave a more complicated polynomial before {\bf FindProperRoot} was added. We cannot use {\bf AZd} for this case. However, we can still obtain a recurrence using {\bf algtorec}. The conversion only took about 13 seconds but produced a recurrence of order 66 and degree 3 and as such is not included here.
				\label{exa:CurrentBasketball}
			\end{example}

		\end{subsection}%

		\begin{subsection}{Time Comparison}
			To reiterate the benefits of this paper and Maple package over guessing for the polynomial, the following examples illustrate the time savings.
			\begin{example}[Finding a Polynomial]
				Let $\S=\{[1,-2],[1,-1],[1,0],[1,1]\}$. Let $G$ be the g.f. for bridges subject to step set $\S$. Then
				\beq
					(16t^3+8t^2+11t-4) G^3	+	(3-2t) G	+	1	=	0.
				\eeq
				This takes
				\begin{table}[H]\begin{center}
					\caption{Finding Minimal Polynomial $\S=\{[1,-2],[1,-1],[1,0],[1,1]\}$}
					\begin{tabular}{|c|cccc|}
						\hline
						Method	&	Memory Used	&	Memory Allocation	&	CPU Time	&	Real Time\\
						\hline
						{\bf EUSP}	&	2.10MiB	&	0 bytes	&	18.27ms	&	21.53ms\\
						\hline
						{\bf Empir}	&	34.94MiB	&	0.73MiB	&	294.93ms	&	295.63ms\\
						\hline
						{\bf EmpirF}	&	1.44MiB	&	0 bytes	&	14.27ms	&	16.23ms\\
						\hline
					\end{tabular}\end{center}
					\label{tab:UnboundedMinPolynomial1}
				\end{table}
				In this case, {\bf EmpirF} appears to be (25\%) faster than {\bf EqualUnBoundedScoringPaths} ({\bf EUSP}) but the time is so small the difference is likely to go unnoticed. The next example shows why the {\bf ScoringPaths} package is superior to using {\bf EmpirF}.
				\label{exa:UnBoundedPoly1}
			\end{example}
			
			\begin{example}[Finding a Polynomial 2]
				Let $\S=\{[1,-2],[1,-1],[1,0],[1,1],[1,2]\}$. Let $G$ be the g.f. for bridges subject to step set $\S$. Then
				\beq
					(5t+4) (5t-1)^2 (t-1)^2 G^4	+	2 (t-1) (5t-2) (5t-1)	G^2	+	t	=	0.
				\eeq
				This takes
				\begin{table}[H]\begin{center}
					\caption{Finding Minimal Polynomial $\S=\{[1,-2],[1,-1],[1,0],[1,1],[1,2]\}$}
					\begin{tabular}{|c|cccc|}
						\hline
						Method	&	Memory Used	&	Memory Allocation Change	&	CPU Time	&	Real Time\\
						\hline
						{\bf EUSP}	&	7.70MiB	&	0 bytes	&	58.33ms	&	61.43ms\\
						\hline
						{\bf Empir}	&	377.91MiB	&	31.29MiB	&	2.89s	&	2.75s\\
						\hline
						{\bf EmpirF}	&	&\multicolumn{2}{c}{Failed}	&\\
						\hline
					\end{tabular}\end{center}
					\label{tab:UnboundedMinPolynomial2}
				\end{table}
				{\bf EmpirF} happened to fail for this $\S$ while {\bf EqualUnBoundedScoringPaths} ({\bf EUSP}) only tripled in time and memory usage. The only difference between this example and Example \ref{exa:UnBoundedPoly1} is the second step with positive $y$-value.
				\label{exa:UnBoundedPoly2}
			\end{example}
			Eventually guessing should work, but we have no way of knowing what the upper bound of search space is before blindly investigating. Expressing the g.f. equations is a set amount of time that is guaranteed to produce the solution. The time complexity to find the minimal polynomial with Groebner bases is an interesting question we leave to the reader.
	
			We have many different ways to enumerate the number of bridges with step set $\S$. We will compare their speed and memory usage here.
			\begin{example}[Speed Enumeration - Unbounded]
				As with most previous examples, we use the benchmark step set $\S=\{[1,-2],[1,-1],[1,0],[1,1],[1,2]\}$. The following are time and memory requirements for enumerating 1000 terms of the sequence of bridges with step set $\S$. The different methods are
				\ben
					\item	Iterating the single polynomial found with {\bf EqualUnBoundedScoringPaths} does not work due to convergence issues beyond the scope of this paper.
					\item	{\bf taylor} cannot be used. It gives the error: \tbl does not have series solution\tbr. {\bf taylor} tends to work only if the minimal polynomial contains a $G^1$ term. See Example \ref{exa:UnBoundedPoly1}.
					\item	Iterating a vector solution with {\bf EqualUnBoundedScoringPathsSeries}.%
					\item	Brute force recursion using {\bf SpecificUnBoundedScoringPathsNumber}.%
					\item Converting minimal polynomial to recurrence:
						\ben
							\item	Obtain the polynomial with {\bf EqualUnBoundedScoringPaths}.%
							\item	Convert to recurrence with {\bf algtorec}. The coefficients satisfy%
								\begin{align*}
									0	&=	125(n+1)(n+2)B(n)	-	25(n+2)(n-1)B(n+1)\\
										&\tab-	5(21n^2+89n+96)B(n+2)	+	(n^2-43n-140)B(n+3)\\
										&\tab+	2(n+4)(n+7)B(n+4).
								\end{align*}
							\item	Enumerate using {\bf SeqFromRec}.%
						\een
					\item	Using {\bf AZd} and then enumerating with {\bf SeqFromRec}.
						\begin{align*}
							0	&=	25(3n+8)(n+2)(n+1)B(n)	-	5(3n+5)(2n+5)(n+2)B(n+1)\\
								&\tab	-	(3n+8)(19n^2+76n+75)B(n+2)	+	2(3n+5)(2n+5)(n+3)B(n+3).
						\end{align*}
				\een
				\begin{table}[H]\begin{center}
					\caption{1000 term Unbounded Enumeration}
					\begin{tabular}{|c|cccc|}
						\hline
						Method	&	Memory Used	&	Memory Allocation	&	CPU Time	&	Real Time\\
						\hline
						Vector Set-Up	&	20.47KiB	&	0 bytes	&	200$\mu$s	&	200$\mu$s\\
						Vector Iterating	&	2.79TiB	&	30.54MiB	&	2.27h	&	110.54m\\
						\hline
						Polynomial	&	7.76MiB	&	0 bytes	&	59.40ms	&	63.03ms\\
						Conversion	&	5.97MiB	&	8KiB	&	49.03ms	&	51.97ms\\
						Enumeration	&	20.64MiB	&	2.16MiB	&	84.87ms	&	81.50ms\\
						\hdashline
						Total	&	34.37MiB	&	2.17MiB	&	193.30ms	&	196.50ms\\
						\hline
						{\bf AZd}	&	2.45MiB	&	0 bytes	&	25.77ms	&	30.53ms\\
						{\bf SeqFromRec}	&	19.34MiB	&	3.59MiB	&	74.40ms	&	74.43ms\\
						\hdashline
						Total	&	21.79MiB	&	3.59MiB	&	100.17ms	&	104.96ms\\
						\hline
						Brute-Force Recursion	&	1.41GiB	&	4.29MiB	&	12.69s	&	11.91s\\
						\hline
					\end{tabular}
					\label{tab:UnboundedEnumeration1000}
				\end{center}\end{table}
				The alternate method is the fastest in this case. It is about half of the time as converting, but the actual enumeration after the preliminary set-up is very similar; the memory usage is similarly low. However, {\bf AZd} is only appropriate for specific cases with constant $x$-step. The key take-away is that enumerating by combining {\bf ScoringPaths} and {\bf algtorec} is much faster and leaner than brute-force recursion and applicable in a wide-range of cases.
			\end{example}
			
			If we change the step set back to $\S=\{[1,-2],[1,-1],[1,0],[1,1]\}$, then we can use {\bf taylor} on the minimal polynomial, as well as iterating the single solution (after an appropriate change of iteration). The difference is that Example \ref{exa:UnBoundedPoly1} had a $G^1$ term in its minimal polynomial while Example \ref{exa:UnBoundedPoly2} did not. This allows one to write $G$ as copies of itself, which is what the iterative method is effectively accomplishing.
		\end{subsection}%
	
		\label{sec:Unbounded}
	\end{section}%

	\begin{section}{Asymptotics}
		For determining asymptotic behavior, I avoided analyzing the bounded cases as those resulted in explicit rational solutions, which can already be handled very easily. The asymptotic behavior of unbounded general walks is simply finding the number of combinations of $x$-steps to yield a length of $n$; the actual walk altitude does not matter.\\
	
		For semi-bounded and unbounded cases, once we have recurrences for the coefficients, we can derive asymptotic expressions. Because of the nature of these quantities, the number of walks of length $n$ will always follow $b^n$ for some base $b$ \cite{Singularity}. Wimp and Zeilberger \cite{wimp1985resurrecting} created a method for automatically determining asymptotics for a linear recurrence, including the constant of the leading term! Their package {\bf AsyRec} provides {\bf Asy}, which attempts to determine the exponential power with which a recurrence grows. {\bf AsyC}, along with suitable starting values, will also give the correct constant.\\
		
		One task we can use this for is finding what proportion of walks are nonnegative. Let $\S=\{[1/2,-1],[1/2,1]\}$.\footnote{$x$-step is $1/2$ so that excursions and bridges do not all have even length, producing a bunch of extraneous $0$s.} Let $F$ denote the g.f. for nonnegative excursions, and $B(n)$ denote the number of nonnegative excursions of length $n$. Then\footnote{All of this can be produced by the one command {\bf EqualSemiBoundedPaper($\S$,20,true,true)}.}
		\begin{align*}
			0	&=	1	-	F	+	t F^2,\\
			0	&=	(4n+2)B(n)	-	(n+2)B(n+1),\\
			B(n)	&\sim	\frac{4^n}{\sqrt{\pi}n^{3/2}}\cdot\left(1-{\frac {9}{8n}}+{\frac {145}{128{n}^{2}}}-{\frac {1155}{1024{n}^{3}}}+{\frac {36939}{32768{n}^{4}}}-{\frac {295911}{262144{n}^{5}}}.
	\right)
		\end{align*}
		Because the minimal polynomial is simple, we could use the Lagrange Inversion Formula to compute exact values of $B(n)$.\footnote{This is what Duchon used for several simple cases \cite{Duchon}.} But for most step sets with $|\S|>2$, this will not work, so the LIF is not considered here. Let $G$ denote the g.f. for bridges and $C(n)$ the number of bridges of length $n$. Then\footnote{This can be produced by {\bf UnBoundedPaper($\S$,20,true,true)}.}
		\begin{align*}
			0	&=	1	+	(4t-1) G^2,\\
			0	&=	(4n+2)C(n)	-	(n+1)C(n+1),\\
			C(n)	&\sim	\frac{4^n}{\sqrt{\pi}\sqrt{n}}\cdot\left(1-\frac{1}{8n}+{\frac {1}{128{n}^{2}}}+{\frac {5}{1024{n}^{3}}}-{\frac {21}{32768{n}^{4}}}-{\frac {399}{262144{n}^{5}}}\right).
		\end{align*}
		Then the proportion of binary bridges that are Dyck paths is asymptotically
		\beq
			\frac{B(n)}{C(n)}	\sim	\frac{1}{n}.
		\eeq
		This matches with the known exact proportion of $\frac{1}{n+1}$.
		
		The above example was not exactly revolutionary, but the method allows for \tbl quick\tbr\ analysis of walks with any step set.\\
		
		Dyck paths have been studied quite extensively. The typical result is for walks with step set $\{[1,0],[0,1]\}$ from $(0,0)$ to $(n,n)$ that stay below the line $y=x$. This is equivalent to enumerating nonnegative excursions with step set $\{[1,-1],[1,1]\}$. Changing the slope of the upper bound to rational $\frac{a}{b}$ is equivalent to using a step set $\{[1,-a],[1,b]\}$; this is called Rational Catalan Combinatorics. For a general conversion from walks below a line with slope $m$, reflect the walk across the line through the origin with slope $\frac{m}{2}$. This produces an equivalent nonnegative walk in the right hand plane. Andrew Lohr studied the field of paths below rational slope with the goal of obtaining asymptotic constants \cite{Andrew}. For his results, go to {\bf http://sites.math.rutgers.edu/$\sim$ajl213/DrZ/RSP.pdf}.%

		Lohr states a result by Duchon in 2000 \cite{Duchon} that for any slope $\frac{a}{b}$, the number of paths below a line of that slope is asymptotically $\Theta\left(\frac{1}{n}\binom{(a+b)n}{an}\right)$. Are we able to say anything about step sets of size $>2$? Intuition would say the growth rate should be larger, but how much larger?
		
		First we need to tackle the scenario if we have steps with $x$-step 0. Recall we can only have steps directly up OR down. The maximum number of up-steps we could have in an excursion or bridge of length $n$ is
		\beq
			n\frac{-\min\{\frac{y}{x}:(x,y)\in\S,x\ne0\}}{\min\{y:(0,y)\in\S\}}.
		\eeq
		Replace the $\min$s with $\max$s for down-steps.\footnote{We technically do not need the $x\ne0$ since we assumed there were no steps directly down.} We can only take so many steps in one direction before we have to start returning in order to make it to the $x$-axis before length $n$. Let $c$ denote the fraction above ($c=0$ for step sets without vertical steps). The \tbl worst case\tbr\ scenario\footnote{Largest number of walks.} is that the remaining steps all have $x$-step 1. Then we have $\binom{n+cn}{n}$ ways to place the vertical steps. Let $S$ denote the size of $\S$ with the 0-steps removed. Then the number of walks is upper bounded (usually very poorly) by
		\beq
			\binom{n(1+c)}{n}\cdot S^n	\approx	\sqrt{\frac{1+c}{2\pi cn}}\left[S(1+c)(1+1/c)^c\right]^n.
		\eeq
		So the number of every type of walk is $O\left(\frac{b^n}{\sqrt{n}}\right)$ for some $b$.\footnote{Which matches with the earlier mentioned asymptotic.} The bound may seem to contradict the Dyck path example, but remember that we cut the steps in half to avoid extraneous 0s. A good goal would be to improve this bound to something more meaningful.\\
		
		\begin{example}
			Let $\S=\{[1, -1], [1,0], [1, 2]\}$. Then, assuming $F,B,G,C$ are as before,
			\begin{align*}
				0	&=	1	+	(t-1)F	+	t^3F^3,\\
				0	&=	31(n+1)(n+2)B(n)-6(2n+5)(n+2)B(n+1)\\
					&\tab	+2(6n^2+36n+53)B(n+2)-2(2n+9)(n+3)B(n+3),\\
				B(n)	&\sim	0.8001188640\cdot\frac{2.889881575^{n}}{n^{3/2}}\cdot\left(1- \frac{1.7475722}{n}- \frac{2.6532889}{n^2}+ \frac{4.0131981}{n^3}\right),
			\end{align*}
			and
			\begin{align*}
				0	&=	1	-	3(t-1)G	+	(31t^3-12t^2+12t-4)G^3,\\
				0	&=	31(n+1)(n+2)C(n)-6(n+2)(2n+3)C(n+1)\\
					&\tab	+2(6n^2+24n+23)C(n+2)-2(n+3)(2n+3)C(n+3),\\
				C(n)	&\sim	0.3488331868\cdot\frac{2.889881575^{n}}{\sqrt{n}}\cdot\left(1-\frac{0.24757219}{n}-\frac{0.03549572}{n^2}+\frac{0.046925761}{n^3}\right).
			\end{align*}
			The exponential bound is simply $|\S|^n=3^n$: fairly close to the actual asymptotic base. And the proportion of excursions to bridges is
			\beq
				\frac{B(n)}{C(n)}	\sim	\frac{2.293700526}{n}:
			\eeq
			also just a constant times $\frac{1}{n}$.
		\end{example}
		
		Finding the asymptotic behavior does not work in every case, e.g., meanders with step set $\{[0, -1], [1, -1], [2, -1], [2, 0], [2, 1]\}$. If the asymptotic does not match with empirical data, then {\bf AsyC} will notify the user of this.\\
		Is the ratio of nonnegative excursions to bridges always some constant times $\frac{1}{n}$? The following example is more experimental evidence for this conjecture.
		\begin{example}
			Let $\S=\{[1,-2],[3,0],[0,1],[2,1],[2,-2]\}$ and $F,B,G,C$ as before. Then
			\begin{align*}
				0	&=	1+(t-1)(t^2+t+1)F+t(1+t)(t^2+1)^2F^3,\\
				B(n)	&\sim	\frac{4}{15}\cdot \frac{7.898354145^n}{n^{3/2}},\\
				0	&=	1-3(t-1)(t^2+t+1)G+(4t^9+15t^6+27t^5+54t^4+66t^3+27t^2+27t-4)G^3,\\
				C(n)	&\sim	\frac{15\sqrt{3}}{58}\cdot \frac{7.898354145^n}{\sqrt{n}}.
			\end{align*}
			The actual recurrences obtained from conversion ({\bf algtorec}) are not listed here due to size. They are degree 2 order 18 and degree 2 order 16, respectively. The base of the exponential is exactly the same for $B$ and $C$; they are the same root of the same polynomial.\footnote{The polynomial in question is $4x^9-27x^8-27x^7-66x^6-54x^5-27x^4-15x^3-4$ and the root has index=1 in Maple notation.} And once again,
			\beq
				\frac{B(n)}{C(n)}	\sim	\frac{232\sqrt{3}}{675}\cdot\frac{1}{n}.
			\eeq
		\end{example}
		To compute the ratio for any step set, use the shorthand {\bf RatioOfWalks($\S$)}.
		
		We can also try to reason the ratio heuristically. Suppose we start with an excursion, $E$. We can reorder the steps cyclically (or in any order) and still have a bridge. Consider the set of such walks; there are size of $E$ such walks. The reason for reordering cyclically would be to maintain the \tbl unique\tbr\ excursion in the set, while the rest are bridges. Actually, $E$ is only unique if it was an irreducible excursion but potentially there is not much overlap.
		
		Similarly, from a bridge $B$ we could reorder the steps cyclically and we must obtain at least 1 excursion: starting from the step after the minimum altitude of $B$. We actually have the same number of excursions in this cyclic set as points of minimum altitude.
		
		We have a relation between excursions and bridges that appears to be roughly linear. Some issues may be that an excursion of length $n$ does not necessarily have size $n$. We still have a linear relationship between the two bounded by $\max\{x,(x,y)\in\S\}$. And even steps with $x$-step 0 do not mess this up too horribly because we are considering excursions and bridges; the walk returns to the $x$-axis so can only go so far away before it must start returning since the return rate $\frac{y}{x}$ is finite (we do not have steps directly up AND down). This again leads to a linear number of $x$-step 0 allowed steps.

		\begin{subsection}{Discriminant}
			Another interesting result from analyzing asymptotic behavior is what the base of the exponent appears to be. For all three of our examples, the base can be found by taking the reciprocal of the smallest modulus of the roots of the coefficient of the leading term in the minimal polynomial of the bridges g.f., i.e.,
			\begin{align*}
				\left(\min\{|z|:4z-1=0\}\right)^{-1}	&=	4,\\
				\left(\min\{|z|:31z^3-12z^2+12z-4=0\}\right)^{-1}	&=	\frac{3}{2^{2/3}}+1	\approx 2.89,\\
				\left(\min\{|z|:4z^9+15z^6+27z^5+54z^4+66z^3+27z^2+27z-4=0\}\right)^{-1}	&=	7.90.\\
			\end{align*}
			This is the first estimate of the asymptotics from the inverse of the radius of convergence. Actually, the polynomial whose root modulus we need can be taken to be the discriminant of the minimal polynomial, which in the step sets given is the same for excursions and bridges.
			
			In general, one method that appears to find the base $b$ for the exponential asymptotic behavior is
			\ben
				\item	Find the discriminant of the minimal polynomial.
				\item	Take the smallest positive real root.
				\item	Take the reciprocal of that root.
			\een
			This will miss the possible sub-exponential factors of $n^{-3/2}$, etc. The shorthand for this computation is implemented as {\bf AsymptoticBase}. A more detailed analysis of singularities and the associated asymptotics is available in {\it Analytic Combinatorics} by Flajolet and Sedgewick \cite{Singularity}.
		\end{subsection}%

		\begin{subsection}{Meanders}
			We have yet to analyze the asymptotic behavior of meanders. Let $K$ denote the g.f. and $D(n)$ the number of nonnegative meanders of length $n$. For the case of Dyck paths with a step set of $\{[1,1],[1,-1]\}$ (since we can have odd length walks now), meanders satisfy 
			\begin{align*}
				1+(2t-1)K+t(2t-1)K^2	&=	0,\\
				4(n+1)D(n)+2D(n+1)-(n+3)D(n+2)	&=	0,\\
				D(n)	&\sim	\frac{0.797}{\sqrt{n}}2^n.
			\end{align*}
			This says there are roughly $2.283$ nonnegative meanders for each bridge of the same length. All of the asymptotic behavior in this section was found using the conversion to recurrence and extracting the asymptotics from there.\\
			
			It turns out that meanders can follow very different asymptotic behavior depending on the step set.
			\begin{example}
				For $\S=\{[1, -1],[1,0], [1,2]\}$,%
				\begin{align*}
					0	&=	1+(4t-1)K+3t(3t-1) K^2+t(3t-1)^2 K^3,\\
					0	&=	93(n+1)(n+2)(n+3)D(n)-2(80n+293)(n+3)(n+2)D(n+1)\\
						&\tab	+(n+3)(115n^2+793n+1318)D(n+2)\\
						&\tab	-4(18n^3+210n^2+820n+1071)D(n+3)\\
						&\tab	+4(n+4)(7n^2+65n+152)D(n+4)-2(n+5)(n+4)(2n+11)D(n+5),\\
					D(n)	&\sim	\frac{3-\sqrt{5}}{2}3^n.
				\end{align*}
				The base of the exponent happens to match the inverse of the smallest modulus of the roots of the discriminant of the minimal polynomial.
			\end{example}
			So meanders can, by virtue of their endpoint flexibility, greatly outnumber excursions and bridges.\\
			
			Finally, meanders do not always follow $|\S|^n$. Let $\S=\{[1,-1],[1,0],[1,1],[2,2]\}$. Then
			\begin{align*}
				0	&=	1+(3t^2+3t-1)K+t(3t+1)(t^2+3t-1)K^2+t^2(t^2+3t-1)^2 K^3,\\
				D(n)	&\sim	0.307\left(\frac{3+\sqrt{13}}{2}\right)^n	\approx	0.307\cdot3.303^n.
			\end{align*}
			The recurrence found was $11^{th}$ order and $3^{rd}$ degree so is not included here.\footnote{The $0.307$ is actually a root of $169x^4-1014x^3+507x^2-78x+4$. Maple's identify command was used on experimental data.} It was produced in less than one second using {\bf algtorec}. And $\S=\{[1,-2],[2,-1],[1,0],[1,2],[2,1]\}$, yields relations of
			\begin{align*}
				0	&=	1+(2t^2+3t-1) K+t(t+2)(2t^2+3t-1) K^2\\
					&\tab	+t(2t^2+3t-1)^2 K^3+(2t^2+3t-1)^2t^2 K^4,\\
				D(n)	&\sim	\frac{7\sqrt{2}}{13\sqrt{n}}\left(\frac{3+\sqrt{17}}{2}\right)^n	\approx	\frac{0.7615}{\sqrt{n}}\cdot3.562^n.
			\end{align*}
			The recurrence for this step set was $25^{th}$ order and of $4^{th}$ degree.\\
			
			The behavior of meanders appears to be a lot harder to peg down than the anticipated $\frac{b^n}{\sqrt{n}}$ for bridges and $\frac{b^n}{n^{3/2}}$ for excursions. Though it appears that meanders are always at least as large as bridges asymptotically.
		\end{subsection}%
		
		\label{sec:Asymptotics}
	\end{section}%

	\begin{section}{Applications}

		\begin{subsection}{Combining Solutions}
			After obtaining all of these g.f.s, we can produce much more. Sports are always a subject of interest for a good portion of the population. Sports statistics are an integral part of many a fan base. So a question of interest beyond the number of ways to be tied, may be the number of ways to win by at least $X$. The task at hand then is how to describe a walk such as that. We have spent the majority of this paper breaking down walks into smaller components. We can now use those components to build other types of walks. The important step is making sure that we count all of our walks, and do not double count any walks.\\
			
			Suppose we want to win by $\ge2$ and never trail by more than 3. $f_{3,5}$ will count walks that drop by no more than 3, and we will have a 2 point lead at the end. Then any walk that stays above that line ($k_0$) will produce what we want. So does $f_{3,5}\cdot k_0$ count what we want? Not necessarily. $f_{3,5}$ ensures that at some point we are exactly 2 points ahead of where we began. But depending on the step set, we may skip over this lead and never actually hit the altitude 2 steps higher than our beginning. In addition, $f_{3,5}$ may finish with a step up and then down, while $k_0$ could start that way; effectively tracing the same walk in \tbl different\tbr\ ways. Thus, $f_{3,5}\cdot k_0$ double counts some walks and misses others.\\
			
			There is at least one way, though not as elegant, of describing these walks. Let $L$ denote the g.f. of interest. Then $L=k_3-f_{3,0}-f_{3,1}-f_{3,2}-f_{3,3}-f_{3,4}$. We count all meanders that never drop more than 3 points, and then remove those that change in altitude by exactly $-3,-2,-1,0,1$.
			\begin{example}
				Let $\S=\{[0, -1], [1, 0], [1, 1], [1, 2]\}$. Then $L$ satisfies
				\begin{align*}
					0	&=	t \bigg( 3{t}^{12}-63{t}^{11}+555{t}^{10}-2673{t}^{9}+7671{t}^{8}-13371{t}^{7}\\
						&\tab\tab	+13745{t}^{6}-7554{t}^{5}+1615{t}^{4}+179{t}^{3}-138{t}^{2}+21t-1 \bigg)\\
						&\tab	+	\bigg( 9{t}^{14}-216{t}^{13}+2286{t}^{12}-13986{t}^{11}+54558{t}^{10}-141402{t}^{9}+246687{t}^{8}\\
						&\tab\tab	-288270{t}^{7}+221709{t}^{6}-109548{t}^{5}+33981{t}^{4}-6448{t}^{3}+715{t}^{2}-42t+1 \bigg) L\\
						&\tab	-	3{t}^{3} \left( 3{t}^{7}-36{t}^{6}+153{t}^{5}-261{t}^{4}+126{t}^{3}+50{t}^{2}-26t+2 \right) {L}^{2}	+	9{t}^{6}{L}^{3},
				\end{align*}
				and has truncated expansion
				\begin{align*}
					L	&=	t+21t^2+305t^3+4064t^4+52431t^5+666657t^6\\
						&\tab	+8420130t^7+106070229t^8+1335635352t^9+16832212452t^{10}.
				\end{align*}
				The minimal polynomial was derived by finding a Groebner basis. We had to describe all of $\{k_3,f_{3,0},f_{3,1},f_{3,2},f_{3,3},f_{3,4}\}$ as well, which included a chain of describing further walks. But those methods are already set and have been shown to be closed.
			\end{example}
		\end{subsection}%

		\begin{subsection}{Weighted Walks}
			One extension that would be fairly easy to implement is adding a weight to each step; $[x,y]$ has an associated weight $w$. Then, for example, we would write\footnote{Taken from general walks in Bounded Section \ref{sec:Bounded}.}
			\beq
				f_{a,b}	=	1+\sum_{(x,y,w)\in\S}w\cdot t^{x}f_{a-y,b-y}.
			\eeq
			Simply multiply by the weight whenever we take a certain step. Now we can accomplish more with the weights in place. If $\sum_{(x,y,w)\in\S}w=1$ (and all $x=1$), then $w$ represents the probability of taking a specific step. And then $f_{a,b}$ is the g.f. for the probability that a given walk of length $n$ maintains altitude $a\ge y\ge b$. If we want to know the probability of a bounded walk being a bounded bridge, use weights to describe $f'_{a,b}$, the g.f. for probability of a general walk being a bounded bridge, and divide its coefficients by the appropriate coefficient of $f_{a,b}$, the probability of a general walk being bounded. This is an explicit description of using $\Pr[A|B]=\frac{\Pr[A\cap B]}{\Pr[B]}$. 
			If all of the weights are the same, then simply enumerate each one to get the probabilities; it is much easier.
			
			\begin{example}
				Let us try finding the probability of a general walk with step set \linebreak$\S=\{[1,2,1/3],[1,-1,1/6],[1,-2,1/2]\}$ being bounded above by $y=3$ and below by $y=-2$. Then the g.f. is explicitly
				\beq
					3{\frac {3888+3888t-756t^2-972t^3-12t^4-{t}^{5}}{11664-7776t^2-432t^3+1296t^4+{t}^{6}}},
				\eeq
				which has taylor expansion
				\beq
					1+t+{\frac{17}{36}}{t}^{2}+{\frac{49}{108}}{t}^{3}+{\frac{77}{324}}{t}^{4}+{\frac{811}{3888}}{t}^{5}+{\frac{53}{432}}{t}^{6}+{\frac{3407}{34992}}{t}^{7}+{\frac{26483}{419904}}{t}^{8}+{\frac{58247}{1259712}}{t}^{9}+O \left( {t}^{10} \right).
				\eeq
			\end{example}
			
			We need the steps to all have the same $x$-value, otherwise we aren't representing the probability that a given walk of length $n$ has some property. We would somehow be combining the probability that a walk has length $n$ and the probability that it satisfies our desired property(ies) in a way I cannot currently describe.
		\end{subsection}%

		\begin{subsection}{2-step Examples}
			With all of the tools at our disposal, let's now use them to produce more information about many sequences.

			\begin{subsubsection}{Bridges with $\S=\{[1,1],[1,-k]\}$}
				It is well known\footnote{Try proving it for yourself. Or prove it for any finite $k$ with this package.} that the g.f. for nonnegative excursions, denoted $f$, satisfies
				\beq
					f	=	1	+	t^{k+1} f^{k+1}.
				\eeq
				We would like to show something about bridges with this step set. The g.f. for the first $k=0,\ldots,7$ step sets have minimal polynomials
				\begin{align*}
					(t-1)G	+	1,\\
					(4t^2-1)G^2	+	1,\\
					(27t^3-4)G^3+3G+1,\\
					(256t^4-27)G^4+18G^2+8G+1,\\
					(3125t^5-256)G^5+160G^3+80G^2+15G+1,\\
					(46656t^6-3125)G^6+1875G^4+1000G^3+225G^2+24G+1,\\
					(823543t^7-46656)G^7+27216G^5+15120G^4+3780G^3+504G^2+35G+1.%
				\end{align*}
				The general form appears to have leading term $([(k+1)t]^{k+1}-k^{k})G^{k+1}$. %
				The $G^1$ term is fairly easy to see have coefficient $(k+1)(k-1)$. The $G^2$ term, after some manipulations, has coefficients that follow $\frac{k^2}{2}(k+1)(k-2)$. %
				The computer could only obtain up to $k=9$ before memory requirements became too large ($>3GB$ allocated). From that limited data, I found that the $G^3$ term has coefficient that follows $\frac{k^3}{6}(k+1)(k-1)(k-3)$. Does this generalize? And if so, how?\\
				
				These results are not groundbreaking as we can already enumerate the number of bridges of length $(a+b)n$ with step set $\{[1,a],[1,-b]\}$ with the simple $\binom{(a+b)n}{an}$. They do allow a different view of the walks by considering how they can be built from copies of themselves.
			\end{subsubsection}%

			\begin{subsubsection}{Duchon Numbers}
				The Duchon numbers\footnote{OEIS \seqnum{A060941} \cite{A060941}.} are one of the cases that Lohr analyzed for asymptotic behavior \cite{Andrew}. We can derive more information for the g.f., which we will denote $f$, of the Duchon numbers. The sequence is defined as the number of paths of length $5n$ from $(0,0)$ to the line $y=2x/3$ with unit North and East steps that stay below the line or touch it. It is equivalent to enumerating excursions with step set $\{[1/5,2],[1/5,-3]\}$. To find the minimal polynomial for its g.f., all we have to do is type {\bf EqualSemiBoundedScoringPaths($\S$,0,t,f)} and hit return.
				\beq
					0	=	1	-	f	+	2t f^5	-	t f^6	+	t f^7	+	t^2	f^{10}.
				\eeq
				We can also produce the minimal polynomial for the related g.f. of irreducible walks (those that only touch at the endpoints). We almost accomplish this with the call \linebreak{\bf SpecificEqualSemiBoundedScoringPaths($\S$,0,t,f,g,g[0,0])}. But to match OEIS \linebreak \seqnum{A293946} \cite{A293946}, we must allow the stationary walk. Simply substitute $g[0,0]=g-1$ and we are done.
				\begin{align*}
					0	&=	{g}^{10}-19{g}^{9}+162{g}^{8}-816{g}^{7}+2688{g}^{6}-\left(2t+6048 \right) {g}^{5}+ \left( 19t+9408 \right) {g}^{4}\\
						&\tab	-\left(73t+9984 \right) {g}^{3}+ \left( 142t+6912 \right) {g}^{2}-\left(140t+2816 \right) g+t^2+56t+512.
				\end{align*}
				\label{subsubsec:Duchon}
			\end{subsubsection}%

			\begin{subsubsection}{Excursions with $\S=\{[1,2],[1,-k]\}$}
				We have some information for this family in the Duchon numbers: $k=3$. But what do they look like in general? We will assume $gcd(2,k)=1$ otherwise we could reduce the steps for an equivalent problem. The first few g.f. have minimal polynomial (we have made the transformation $t\to t^{1/(2+k)}$ for ease of reading)
				\begin{align*}
					\seqnum{A001764}:\hspace{1em}	k=1	\hspace{1em}	&	\tab	f^{3}t-f+1,\\
					\seqnum{A060941}:\hspace{1em}	k=3	\hspace{1em}	&	\tab	{f}^{10}{t}^{2}+{f}^{7}t-{f}^{6}t+2{f}^{5}t-f+1,\\
					\seqnum{A300386}:\hspace{1em}	k=5	\hspace{1em}	&	\tab	{f}^{21}{t}^{3}+2{f}^{16}{t}^{2}-{f}^{15}{t}^{2}+3{f}^{14}{t}^{2}\\
							&\tab\tab	+{f}^{11}t-{f}^{10}t+2{f}^{9}t-2{f}^{8}t+3{f}^{7}t-f+1,\\
					\seqnum{A300387}:\hspace{1em}	k=7	\hspace{1em}	&	\tab	{f}^{36}{t}^{4}+3{f}^{29}{t}^{3}-{f}^{28}{t}^{3}+4{f}^{27}{t}^{3}+3{f}^{22}{t}^{2}-2{f}^{21}{t}^{2}+6{f}^{20}{t}^{2}\\
							&\tab\tab	-3{f}^{19}{t}^{2}+6{f}^{18}{t}^{2}+{f}^{15}t-{f}^{14}t+2{f}^{13}t-2{f}^{12}t\\
							&\tab\tab	+3{f}^{11}t-3{f}^{10}t+4{f}^{9}t-f+1,\\
					\seqnum{A300388}:\hspace{1em}	k=9	\hspace{1em}	&	\tab	{f}^{55}{t}^{5}+4{f}^{46}{t}^{4}-{f}^{45}{t}^{4}+5{f}^{44}{t}^{4}+6{f}^{37}{t}^{3}-3{f}^{36}{t}^{3}+12{f}^{35}{t}^{3}\\
							&\tab\tab	-4{f}^{34}{t}^{3}+10{f}^{33}{t}^{3}+4{f}^{28}{t}^{2}-3{f}^{27}{t}^{2}+9{f}^{26}{t}^{2}-6{f}^{25}{t}^{2}\\
							&\tab\tab	+12{f}^{24}{t}^{2}-6{f}^{23}{t}^{2}+10{f}^{22}{t}^{2}+{f}^{19}t-{f}^{18}t+2{f}^{17}t\\
							&\tab\tab	-2{f}^{16}t+3{f}^{15}t-3{f}^{14}t+4{f}^{13}t-4{f}^{12}t+5{f}^{11}t-f+1,\\
					\seqnum{A300389}:\hspace{1em}	k=11	\hspace{1em}	&	\tab	{f}^{78}{t}^{6}+5{f}^{67}{t}^{5}-{f}^{66}{t}^{5}+6{f}^{65}{t}^{5}+10{f}^{56}{t}^{4}-4{f}^{55}{t}^{4}+20{f}^{54}{t}^{4}\\
							&\tab\tab	-5{f}^{53}{t}^{4}+15{f}^{52}{t}^{4}+10{f}^{45}{t}^{3}-6{f}^{44}{t}^{3}+24{f}^{43}{t}^{3}-12{f}^{42}{t}^{3}\\
							&\tab\tab	+30{f}^{41}{t}^{3}-10{f}^{40}{t}^{3}+20{f}^{39}{t}^{3}+5{f}^{34}{t}^{2}-4{f}^{33}{t}^{2}+12{f}^{32}{t}^{2}\\
							&\tab\tab	-9{f}^{31}{t}^{2}+18{f}^{30}{t}^{2}-12{f}^{29}{t}^{2}+20{f}^{28}{t}^{2}-10{f}^{27}{t}^{2}\\
							&\tab\tab	+15{f}^{26}{t}^{2}+{f}^{23}t-{f}^{22}t+2{f}^{21}t-2{f}^{20}t+3{f}^{19}t-3{f}^{18}t\\
							&\tab\tab	+4{f}^{17}t-4{f}^{16}t+5{f}^{15}t-5{f}^{14}t+6{f}^{13}t-f+1.
				\end{align*}
				The degree is simply $\frac{1}{2}(k+1)(k+2)$. If one looks closer they may recognize that the degrees then decrease at a consistent rate. The degree drops by $k$ and then by $1$ twice to create a \tbl group\tbr\ of 3 terms. The degree then drops again by $k-2$ and then by $1$ to create a group of 5 terms. It can be seen that each polynomial follows this $1, 3, 5, 7,\ldots$ pattern, with the power of $t$ decreasing by 1 in each successive group. There is always a $-f+1$ included.
				
				We have empirically shown the general form of the minimal polynomial, but were unable to describe what the coefficients themselves are.
				
				\seqnum{A001764} \cite{A001764} and \seqnum{A060941} \cite{A060941} are sequences currently in the OEIS. \seqnum{A300386} \cite{A300386}, \seqnum{A300387} \cite{A300387}, \seqnum{A300388} \cite{A300388}, and \seqnum{A300389} \cite{A300389} are new and have been recently submitted and accepted.
			\end{subsubsection}%

			\begin{subsubsection}{Excursions with $\S=\{[1,3],[1,-k]\}$}
				Let us push our computers further. What does this family look like? We cannot derive a lot of information empirically as the $\{[1,3],[1,-5]\}$ case already takes 28 seconds to run. The case $k=7$ ran for over one day and had used 3GiB of allocated memory before it was terminated. Again, we made the transformation $t\to t^{1/(3+k)}$ for compact reading.
				\begin{align*}
					\seqnum{A002293}:\hspace{1em}	k=1	\hspace{1em}	&	\hspace{1em}	f^4t	-	f	+	1,\\
					\seqnum{A060941}:\hspace{1em}	k=2	\hspace{1em}	&	\hspace{1em}	{f}^{10}{t}^{2}+{f}^{7}t-{f}^{6}t+2{f}^{5}t-f+1,\\
					\seqnum{A300390}:\hspace{1em}	k=4	\hspace{1em}	&	\hspace{1em}	{f}^{35}{t}^{5}-{f}^{31}{t}^{4}+{f}^{30}{t}^{4}-{f}^{29}{t}^{4}+5{f}^{28}{t}^{4}-{f}^{25}{t}^{3}+{f}^{24}{t}^{3}+3{f}^{23}{t}^{3}\\
							&\tab	-4{f}^{22}{t}^{3}+10{f}^{21}{t}^{3}+{f}^{19}{t}^{2}-{f}^{18}{t}^{2}+5{f}^{17}{t}^{2}+3{f}^{16}{t}^{2}-6{f}^{15}{t}^{2}\\
							&\tab	+10{f}^{14}{t}^{2}+{f}^{13}t-{f}^{12}t+3{f}^{10}t+{f}^{9}t-4{f}^{8}t+5{f}^{7}t-f+1,\\
					\seqnum{A300391}:\hspace{1em}	k=5	\hspace{1em}	&	\hspace{1em}	{f}^{56}{t}^{7}-2{f}^{51}{t}^{6}+{f}^{50}{t}^{6}-{f}^{49}{t}^{6}+7{f}^{48}{t}^{6}+{f}^{46}{t}^{5}-{f}^{45}{t}^{5}-3{f}^{43}{t}^{5}\\
							&\tab+5{f}^{42}{t}^{5}-6{f}^{41}{t}^{5}+21{f}^{40}{t}^{5}-3{f}^{37}{t}^{4}-3{f}^{36}{t}^{4}+8{f}^{35}{t}^{4}+10{f}^{34}{t}^{4}\\
							&\tab	-15{f}^{33}{t}^{4}+35{f}^{32}{t}^{4}-2{f}^{31}{t}^{3}+2{f}^{30}{t}^{3}-9{f}^{28}{t}^{3}+22{f}^{27}{t}^{3}\\
							&\tab+10{f}^{26}{t}^{3}-20{f}^{25}{t}^{3}+35{f}^{24}{t}^{3}+3{f}^{22}{t}^{2}+5{f}^{21}{t}^{2}-9{f}^{20}{t}^{2}\\
							&\tab	+18{f}^{19}{t}^{2}+5{f}^{18}{t}^{2}-15{f}^{17}{t}^{2}+ {f}^{16}\left( 21t+1 \right)t-{f}^{15}t+3{f}^{13}t\\
							&\tab	-3{f}^{12}t+5{f}^{11}t+{f}^{10}t-6{f}^{9}t+7{f}^{8}t-f+1.
				\end{align*}
				\seqnum{A002293} \cite{A002293} and \seqnum{A060941} \cite{A060941} are already in the OEIS while \seqnum{A300390} \cite{A300390} and \seqnum{A300391} \cite{A300391} are new. There appears to be some pattern again in the degree of $f$, though it is not discernible at first with so few datum. The degree happens to follow $\binom{3+k}{3}$.
			\end{subsubsection}%
		
			\begin{conjecture}
				Let $\S=\{[1,a],[1,-b]\}$ with $gcd(a,b)=1$. Let $f$ denote the g.f. for excursions (or meanders or bridges) with step set $\S$. Then the minimal polynomial of $f$ has degree $\binom{a+b}{a}$.
			\end{conjecture}
			The above conjecture is supported by all of the examples in this Section \ref{subsec:2stepExamples}, as well as several other quick tests of accuracy. Though what does this mean in terms of deconstructing a walk with step set $\S$?%
			
			This is somewhat intuitive for the unbounded case since there are $\binom{a+b}{a}$ ways to have a bridge of smallest length: $a+b$. And we have provided ample evidence that the asymptotic behavior is only off by $\frac{c}{n}$ for these excursions.

			\label{subsec:2stepExamples}
		\end{subsection}%
		
	\end{section}%

	\begin{section}{Conclusion and Future Work}

		To analyze walks, we began by examining the first step (Bounded Section \ref{sec:Bounded}), the last step (Semi-bounded Section \ref{sec:SemiBounded}), and finally a middle step that crosses the $x$-axis (Unbounded Section \ref{sec:Unbounded}).
		
		The way we dissected the g.f. equations is not unique. You could describe them in slightly different ways that may be more optimal. However, the minimal polynomial is called that for a reason; there is no \tbl better\tbr\ way to describe the g.f. except for an exact solution in special cases. It is important to make sure you do not double count or miss any walks in your descriptions. Define your various types of walks in a very particular manner.\\
		
		For bounded cases we produced the new OEIS sequences \seqnum{A301379} \cite{A301379} and \seqnum{A301380} \cite{A301380} as well as \seqnum{A301381} \cite{A301381} and \seqnum{A300998} \cite{A300998}, which are not included in this thesis, though they were produced in conjunction. Enumerating some semi-bounded and unbounded examples may be much faster by actually solving the polynomial for the g.f. in cases where $\deg(p)\le4$ or $p$ has \tbl nice\tbr\ roots.\\
	
		One interesting note that appears in the semi-bounded case is that, no matter what we have chosen as our step set, the minimal polynomial has had terms $-F+1$. This seems to indicate that there is always a way to write semi-bounded excursions and meanders as some combination (and deductions) of copies of ONLY ITSELF. Though it may be that the self-description is extremely complicated seemingly without (but it must be there) any combinatorial interpretation: see Example \ref{exa:EqualSemiBounded}. We cannot necessarily do that with unbounded cases because there is not always a $-F$ term.
		
		Another trick we can use the minimal polynomial for is a bijective proof. All of the walks counted by terms with positive coefficients are also equally counted by those terms with negative coefficients. The trouble with this is that the two sides are very artificial and rarely something of interest by themselves.\\

		In Section \ref{sec:Asymptotics}, we tried looking at the asymptotics of excursions and bridges. One note that came up was the relationship between the number of excursions and bridges. It is known that the number of Dyck paths is $\frac{1}{n+1}$ of the total number of bridges. With a different looking step set, we still obtained a $\frac{c}{n}$ asymptotic relationship (for a constant $c$). Is this always the case? We provided some heuristics and examples in support of the relationship but no found definitive proof. The relationship to meanders of the same step set appears to be much harder to state exactly.\\
	
		A little further in Section \ref{subsec:2stepExamples}, we examined many 2-step cases. I contributed the new sequences \seqnum{A300386} \cite{A300386}, \seqnum{A300387} \cite{A300387}, \seqnum{A300388} \cite{A300388}, \seqnum{A300389} \cite{A300389}, \seqnum{A300390} \cite{A300390}, and \seqnum{A300391} \cite{A300391} to the OEIS. These are equivalent to walks that stay below certain lines of rational slope. A further question: how does one translate walks with general step sets bounded by lines into bridges, excursions, or meanders? What about bounded by something other than straight lines?\\
		
		An extension that would be fairly easy, though laborious to implement, would be generalizing from 2D walks to 3D walks. Or to n-dimensional walks.\\
	
		Beyond simply enumerating walks, we might want to know more about them: how many peaks or valleys do they contain? what is the area beneath the curve? how many times does the walk hit its maximal/minimal altitude? We can try to answer these questions by adding in a catalytic variable to count this new measurement. The generating function relations are very similar, but the g.f.s themselves are now functions of 2 (or more) variables. This can lead to systems that are not closed (under current descriptions). However, the system can still be iterated to enumerate terms.
		
		Ayyer and Zeilberger analyzed how many times a bounded bridge or excursion hits its boundaries \cite{BoundedCase} in the bounded and semi-bounded cases. The extension to measuring the area under the curve has been started. The bounded case still yields explicit solutions but the semi-bounded and unbounded cases are now left as a system of equations rather than a minimal polynomial.\\

		Thank you for reading this chapter. I hope you have enjoyed it and can make use of this package.
	\end{section}%

	\label{chap:LatticeWalkEnumeration}
\end{chapter}%

\begin{chapter}{Conclusion}
	Let us look back at what we have accomplished.\\
	
	In Chapter \ref{chap:UnimodalPolynomials} we adapted a method of proving the unimodal nature of $q$-binomials to build more families of unimomdal polynomials. The special facet of the constructed families is that the unimodal nature is not immediately obvious. Either an explicit formula or a recurrence with the correct initial conditions produce surprisingly unimodal polynomials.\\
	
	In order to enumerate lattice walks in Chapter \ref{chap:LatticeWalkEnumeration} we generalized the technique of describing generating function relations with irreducible walks. This allowed us to give exact generating functions in the bounded case and minimal polynomials in semi-bounded and unbounded cases. Connected to these minimal polynomials are recurrences for enumeration. Once we take a step back, we realize that they produce surprisingly integral sequences given the correct initial conditions.\\
	
	In both major chapters, we produced new sequences for the OEIS. We have revised results of mathematicians before us to work in more ways than originally intended. We created two Maple packages to accompany the work and provide avenues for quick and easy investigation into more general questions. We contributed to the field of experimental mathematics by devising methods for exploration and using them to gain insights into unimodal polynomial construction and lattice walk generating functions that we otherwise could not have gleaned.
\end{chapter}%

\bibliographystyle{amsalpha}
\bibliography{Ek_thesis_bibliography}

\end{document}